\documentclass[twoside]{article}

\usepackage{geometry}
 \usepackage{fancyhdr}
\usepackage{scrhack}

\usepackage[square,numbers,super,sort&compress]{natbib}

\makeatletter
\renewcommand\NAT@citesuper[3]{\ifNAT@swa
\if*#2*\else#2\NAT@spacechar\fi
\unskip\kern\p@\textsuperscript{\NAT@@open#1\if*#3*\else,\NAT@spacechar#3\fi\NAT@@close}%
   \else #1\fi\endgroup}
\makeatother

\usepackage[utopia]{mathdesign}

\usepackage{amsfonts, dsfont}
\usepackage{amsmath}
\usepackage{mathtools}
\usepackage{amsthm}
\usepackage{stmaryrd}
 
\usepackage{framed}
\usepackage{multicol}

\usepackage{color}
\usepackage{array}
\usepackage{booktabs}
\usepackage[final]{graphicx}
\usepackage{enumitem}
\usepackage{url}
\usepackage{varwidth}
\usepackage{hhline}

\usepackage{mathtools}

\usepackage[
  normal,
  font={small,color=black},
  labelfont=bf,
  margin=2em,
  labelsep=period]{caption}[2008/04/01]

\usepackage{caption}

\renewcommand{\kappa}{\varkappa}

\newtheorem{varform}{Variational Formulation}
 \newcommand{\jump}[1]{\llbracket #1 \rrbracket }
 
 \providecommand{\keywords}[1]{\textbf{Keywords}  #1}

\newtheorem{lemma}{Lemma}[section]
\newtheorem{theorem}{Theorem}[section]
\newtheorem{remark}[lemma]{Remark}
\newtheorem*{acknowledgement}{Acknowledgement}

\begin{document}

\title{A Nitsche-based formulation for fluid-structure interactions with contact}

\author{Erik Burman\thanks{Department of Mathematics, University College London, Gower Street, WC1E 6BT, London, UK
(e.burman@ucl.ac.uk)} \and
Miguel A. Fern\'andez\thanks{Inria Paris, 75012 Paris \& Sorbonne Universit\'e \& CNRS, UMR 7598 LJLL, 75005 Paris, France (miguel.fernandez@inria.fr)}
\and Stefan Frei\thanks{Department of Mathematics, University College London, Gower Street, WC1E 6BT, London, UK (s.frei@ucl.ac.uk)}}

\date{}

\maketitle

% \author[1]{Erik Burman}
% 
% \author[2,3]{Miguel A. Fern\'andez}
% 
% \author[1]{Stefan Frei}
% 
% \authormark{Burman \textsc{et al}}
% 
% \title{} 
% 
% 
% \pagestyle{myheadings}
% \thispagestyle{plain}
% \markboth{}{}
% 
% 
% \address[1]{\orgdiv{Department of Mathematics}, \orgname{University College London}, \orgaddress{Gower Street, WC1E 6BT, London, \country{United Kingdom}}}
% 
% \address[2]{\orgname{Inria Paris, 75012 Paris, France}}
% \address[3]{\orgname{Sorbonne Universit\'e \& CNRS, UMR 7598 LJLL, 75005 Paris, France}}
% 
% \corres{*Stefan Frei \email{s.frei@ucl.ac.uk}}

%\presentaddress{This is sample for present address text this is sample for present address text}

\abstract{We derive a Nitsche-based formulation for 
fluid-structure interaction (FSI) problems with contact. The approach is
based on the work of Chouly and Hild [SIAM Journal on Numerical Analysis. 2013;51(2):1295--1307]
for contact problems in
solid mechanics. We present two numerical approaches,
both of them formulating the FSI interface and the contact conditions simultaneously
in equation form on a joint interface-contact surface $\Gamma(t)$.
The first approach uses a relaxation of the contact conditions to allow for a small 
mesh-dependent gap between solid and wall. The second alternative introduces an artificial fluid below the 
contact surface.
The resulting systems of equations
can be included  {in a consistent fashion} within a monolithic variational formulation,
which prevents the so-called ``chattering'' phenomenon.
To deal with the topology changes in the fluid domain at the time of impact, we 
use a  fully Eulerian approach for the FSI problem. We compare the effect of
slip and no-slip interface conditions and study the performance of the method by means
of numerical examples.\\}


\section{Introduction}

Contact problems have to be considered in many physical processes in engineering, medicine and nature. 
To name only a few consider for example the contact of balls and races in roller bearings, 
closing and opening heart valves or falling and jumping bouncing balls.
While an extensive amount of literature exists for the numerical simulation of contact in a purely mechanical context (see for example Wohlmuth~\cite{Wohlmuth2011}
for an overview), 
i.e. disregarding the gas or liquid that mostly lies between contacting structures, 
much less works can be found considering full fluid-structure interaction with contact. 
The flow between contacting surfaces might however be of great importance for the contact dynamics. 
In the example of heart valves, the pulsating blood flow is even the driving force that enables opening and closure.
In the case of ball bearings, fluid forces in the lubricant between ball and 
bearing may have a significant influence on the performance and wear of the bearing~\cite{Knaufetal, Bruyere2012}.

Contact between different structures is typically formulated by means of variational inequalities~\cite{Wohlmuth2011}. 
In the context of 
full fluid-structure interaction (FSI), first results and algorithms can be found  using either {an artificial} penalty 
force~\cite{TezduyarSathe2007} or Lagrange 
multipliers~\cite{MayerWalletal2013, DosSantosEtAl2008, AstorinoGerbeauetal2009} to obtain a well-posed 
and computationally feasible variational formulation.

However, these {approaches} 
have certain drawbacks: The use of a penalty force prevents 
real contact. The force is typically an artificial force 
and involves the choice of penalty parameters. If they are chosen too small, the structures might 
overlap in a numerical simulation. If they are chosen too large, the contact dynamics 
might be significantly perturbed~\cite{FreiRichter2017Sammelband}. In the case of Lagrange multipliers, 
additional variables are introduced on the 
contacting surfaces and an inf-sup condition is needed to ensure the well-posedness of the system.
To tackle the variational inequality numerically, an extra loop is usually used in each time step (for example within 
an \textit{active-set strategy}), meaning in particular that the system of equations has to be solved 
several times in each time step.

Recently, a new contact formulation using Nitsche's method~\cite{Nitsche70} was derived by 
Chouly and co-workers~\cite{ChoulyHild2013, ChoulySymVsNonsym, Choulyetal2015},
following the ideas of Alart and Curnier~\cite{AlartCurnier91}. Their approach is based on an equivalent re-formulation 
of the contact conditions in equality form. In the case of the contact of an elastic body with a wall, the
mechanical contact conditions on the contact surface $\Gamma_C(t)$ read
\begin{align*}
 d_n \leq 0, \;\; \sigma_{s,n}(d) \leq 0, \;\;\sigma_{s,n}(d) d_n = 0,
\end{align*}
where $d_n$ denotes the solid displacement in normal direction and $\sigma_{s,n}(d)$
is the normal stress component. It can be shown that these conditions are equivalent to the equality
\begin{align}\label{reform}
 \sigma_{s,n}(d)  = -\frac{1}{\gamma} \max\{0,d_n - \gamma \sigma_{s,n}(d)\} \quad \text{on } \Gamma_C(t)
\end{align}
for arbitrary $\gamma>0$ (Chouly \& Hild~\cite{ChoulyHild2013}). The authors incorporate this equality condition weakly in
the variational formulation using Nitsche's method.
This approach has the advantage that it is fully consistent 
and hence the contact dynamics will not be perturbed. Furthermore, no additional variables have to 
be introduced and no additional loop within each time-step is needed. 
The authors were able to prove numerical convergence in a series of papers for friction-free 
and frictional contact~\cite{ChoulySymVsNonsym, ChoulyTrescaFriction, ChoulyCoulombFriction}, 
disregarding however the fluid that usually lies between the structures.

Following these works, Burman \emph{\emph{et al.}}$\,$used the re-formulation (\ref{reform}) to derive a
Galerkin Least Squares formulation in equality form for the obstacle 
problem~\cite{BurmanetalObstacle} and a Galerkin Least Squares and a Lagrange multiplier formulation
for membrane contact~\cite{BurmanHansboLarsonMembrane}. {Annavarapu et al.~\cite{Annavarapuetal2014} used a Nitsche formulation to model frictional 
sliding between two solid bodies. In order to derive a Nitsche-based formulation for FSI with contact, we will introduce Lagrange multipliers on the FSI interface 
first, which can then be eliminated in a second step, following the concepts presented by Burman \& Hansbo~\cite{BurmanHansbo2018}.
}

Modelling of contact in an FSI context brings along a further issue: It is unclear, whether the incompressible 
Navier-Stokes equations are an appropriate model in the fluid part, when it comes to contact. Theoretical studies show, that for a smooth, rigid
solid body, no contact with an exterior wall can happen, when no-slip conditions are used on the interface and the outer boundary of the fluid domain,
see Hillairet~\cite{Hillairet2d} and Hesla~\cite{HeslaPhD} in 2 space dimensions and Hillairet \emph{et al.}~\cite{HillairetTakahashi3d, GerardVaretetal2015}
in 3 space dimensions. This changes, when slip- or Navier-slip conditions are used on both the interface and the wall~\cite{GerardVaretetal2015} or when the 
boundary of the solid is non-smooth~\cite{GerardVaretHillairet, Wang2014}. Gerard-Varet and Hillairet~\cite{GerardVaretHillairet} found in a model example that 
it comes to contact for a solid with a $C^{1,\alpha}$-parametrised boundary for $\alpha<1/2$, while no contact happens for $\alpha\geq1/2$. 
In the context
of fluid-structure interactions, the regularity of the solid boundary depends on the solid displacement $d$, for which such a regularity can usually
not be guaranteed. 

For a full FSI problem with {a thin-walled} structure, 
a no-collision result has been shown by Grandmont and Hillairet~\cite{GrandmontHillairet}
in the no-slip case. For an overview on further results regarding existence of 
fluid-structure interaction problems, we refer to Grandmont et al~\cite{GrandmontEtAlInBook}. 
Recently, Muha \& \v{C}ani\'c~\cite{MuhaCanic} showed the well-posedness of a 
fluid-structure interaction system with 
slip-conditions. 

Motivated by {these} theoretical results, we will 
study both no-slip and slip conditions on the FSI interface $\Gamma(t)$
in this work. It will turn out that the latter transits naturally into a ``no-friction'' 
condition when it comes to contact, while the prior leads to frictional 
contact. {In this work, we will therefore consider friction-free contact, 
when a slip-condition is used on $\Gamma(t)$ and the 
specific frictional contact condition that follows from the transition of the 
interface conditions, when a no-slip condition is used on $\Gamma(t)$.}
For recent works on the incorporation of different 
friction laws (in particular
Coulomb and Tresca friction), we refer to Chouly and co-workers~\cite{ChoulyTrescaFriction, ChoulyCoulombFriction}.
Moreover, we will study only contact of a deformable elastic structure with a fixed 
and straight wall for simplicity. Efficient algorithms 
to treat contact between more 
complex structures can be found, for example, in Puso~\cite{Puso2004}, Yang et al~\cite{Yangetal2005}
and {Chouly, Mlika \& Renard~\cite{ChoulyMlikaRenard18, milka-et-al-17}}.

{Concerning the governing equations, we focus in this work on linear model equations for the fluid and solid part, i.e., the incompressible Stokes equations
in the fluid and linear elasticity in the solid sub-domain. These simplifications must be seen as a first step towards the derivation of Nitsche-based contact 
formulations for complex FSI-contact problems.
We consider, however, the case of a moving interface, which is a major challenge from the numerical point of view and 
leads to a non-linear FSI system, already in absence of contact.}

The change of topology in the fluid domain causes additional numerical difficulties. Standard numerical approaches 
as the Arbitrary Lagrangian Eulerian method are not able to deal with topology changes, as the map from the reference domain to the 
Eulerian domain degenerates necessarily in this situation. The dynamics shortly before the impact can only be {handled robustly}, when a Eulerian
description of the fluid equations is used.

In the last years, several numerical approaches have been developed that 
are able to deal with topology changes.
The methods can be split into {two} categories, according
to the coordinate systems that are used for the solid system: Fully Eulerian approaches, where 
also the structure equations are formulated in Eulerian
coordinates~\cite{DunneRannacher, Cottetetal, Richter2012b, FreiRichter2017Sammelband, Pironneau2016}; {and}
Euler-Lagrangian techniques, where Lagrangian coordinates are used for
the solid equations~\cite{Peskin1972, BoffiGastaldi2003, ZhangGerstenberger,LegayChessaBelytschko2006,DosSantosEtAl2008,GerstenbergerWall,BurmanFernandez2014,alauzet-et-al-15,MassingLarsonetal}, {such as in the}
 \textit{Immersed Boundary} or \textit{Immersed Finite Element} methods.

Regarding the \textit{Euler-Lagrangian} techniques, one can further distinguish between methods using Lagrange multipliers for the coupling of fluid and structure
(Legay \emph{et al.}~\cite{LegayChessaBelytschko2006}, Gerstenberger \& Wall~\cite{GerstenbergerWall}) 
and methods based on Nitsche-techniques
(Hansbo \emph{et al.}~\cite{HansboHermanssonSvedberg2004}, Burman \& Fern\'andez~\cite{BurmanFernandez2014}, {Alauzet \emph{et al.}~\cite{alauzet-et-al-15}}, Massing \emph{et al.}~\cite{MassingLarsonetal}, {Kamensky \emph{et al.}~\cite{kamensky-et-al-15}}). For the latter, a theoretical stability and convergence analysis 
has been derived~\cite{BurmanFernandez2014}. {The reader is referred to {Boilevin-Kayl \emph{et al.}}~\cite{boilevinkayl:hal-01704575} for a comparative study on the accuracy of some of these approaches}. 

The FSI approach we use here is based on the monolithic \textit{Fully Eulerian} approach 
(Dunne \& Rannacher~\cite{DunneRannacher}, Cottet \emph{et al.}~\cite{Cottetetal}, 
Richter~\cite{Richter2012b}, Frei \& Richter~\cite{FreiRichter2017Sammelband}, Hecht \& Pironneau~\cite{Pironneau2016}). As the complete system of equations is 
formulated in Eulerian coordinates, the incorporation of contact
conditions is straight-forward by means of variational principles. While some of the early works 
in this context suffered from stability and accuracy issues, see e.g.{$\,$Dunne}~\cite{DunnePhD}, accurate and 
robust discretisation and stabilisation techniques have been developed recently 
(Frei \& Richter~\cite{FreiRichter2014, FreiRichter2017, FreiRichter2017Sammelband}, Hecht \& Pironneau~\cite{Pironneau2016}).
{We remark, however, that the algorithms we derive to incorporate contact can be combined in a straight-forward way with different FSI coupling techniques, 
e.g.$\,$ Fictitious Domain or Immersed Boundary methods.}

Concerning discretisation, we allow both for unfitted and fitted finite element approaches. 
For the unfitted case, so-called ``ghost penalty'' stabilisations can be used to
guarantee the
coercivity of the system~\cite{Burman2010,BurmanFernandez2014}. In order to simplify the 
presentation, 
we concentrate on fitted discretisations in this work and will comment on the unfitted case in a remark.
In the numerical examples at the end 
of this paper, we will use the 
fitted \textit{locally modified finite element method}~\cite{FreiRichter2014}.

The remainder of this paper is organised as follows: In Section~\ref{sec.model}, we first introduce 
the equations and the contact model. Then, we derive a variational formulation for the problem of an obstacle
within the fluid domain {in Section~\ref{sec.ObstModel}}, where we already have to deal with variational inequalities, but a topology change
in the fluid domain is avoided. In Section~\ref{sec.fullContact}, we study contact with an exterior wall and
discuss in particular the effect of slip- and no-slip interface conditions. Then, we show a stability result in
Section~\ref{sec.discstab}.
We show detailed numerical studies {for a model problem on a simple geometry} in Section~\ref{sec.num}, investigating the influence of 
contact parameters, interface conditions and different contact formulations
as well as convergence under mesh refinement. We conclude in 
Section~\ref{sec.conclusion}.

\section{Model}
\label{sec.model}

We begin by presenting the models for the fluid part, the solid part and the fluid-structure interaction
on one hand in Section~\ref{sec.FSI} and the contact model in Section~\ref{sec.ContactModel} on the other hand.
For both models, Nitsche-based variational formulations are introduced.
{Different possibilities to combine the two models will then 
be presented in Section~\ref{sec.FSIcontact}}.

\subsection{Fluid-structure interaction without contact}
\label{sec.FSI}

\begin{figure}[bt]
\centering
\resizebox*{0.9\textwidth}{!}{
\begin{picture}(0,0)%
\includegraphics{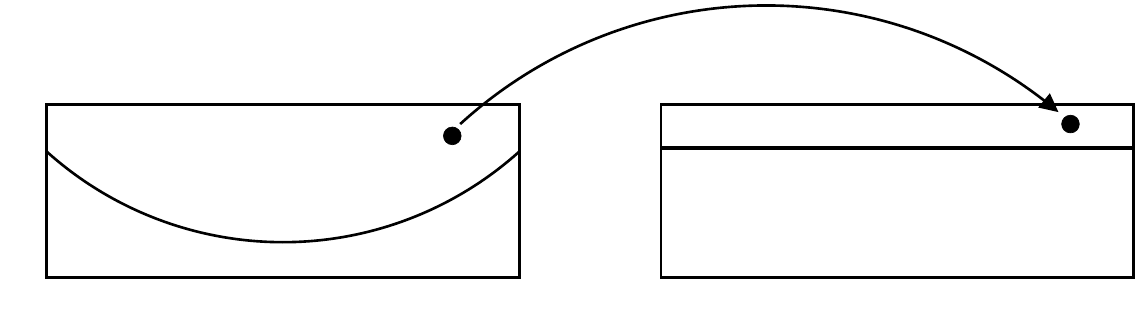}%
\end{picture}%
\setlength{\unitlength}{1657sp}%
\begingroup\makeatletter\ifx\SetFigFont\undefined%
\gdef\SetFigFont#1#2{%
  \fontsize{#1}{#2pt}%
  \selectfont}%
\fi\endgroup%
\begin{picture}(12996,3533)(1269,-3339)
\put(2296,-2671){\makebox(0,0)[lb]{\smash{{\SetFigFont{8}{9.6}{\color[rgb]{0,0,0}$\Omega_f(t)$}%
}}}}
\put(7291,-1321){\makebox(0,0)[lb]{\smash{{\SetFigFont{8}{9.6}{\color[rgb]{0,0,0}$\Gamma_s^D$}%
}}}}
\put(7291,-2581){\makebox(0,0)[lb]{\smash{{\SetFigFont{8}{9.6}{\color[rgb]{0,0,0}$\Gamma_f^N$}%
}}}}
\put(3826,-2356){\makebox(0,0)[lb]{\smash{{\SetFigFont{8}{9.6}{\color[rgb]{0,0,0}$\Gamma(t)=\Gamma_{\text{fsi}}(t)$}%
}}}}
\put(2701,-1636){\makebox(0,0)[lb]{\smash{{\SetFigFont{8}{9.6}{\color[rgb]{0,0,0}$\Omega_s(t)$}%
}}}}
\put(9901,-2536){\makebox(0,0)[lb]{\smash{{\SetFigFont{8}{9.6}{\color[rgb]{0,0,0}$\Omega_f(0)$}%
}}}}
\put(9901,-1276){\makebox(0,0)[lb]{\smash{{\SetFigFont{8}{9.6}{\color[rgb]{0,0,0}$\Omega_s(0)$}%
}}}}
\put(6616,-1411){\makebox(0,0)[lb]{\smash{{\SetFigFont{8}{9.6}{\color[rgb]{0,0,0}$x$}%
}}}}
\put(9586,-151){\makebox(0,0)[lb]{\smash{{\SetFigFont{8}{9.6}{\color[rgb]{0,0,0}$T(x,t)$}%
}}}}
\put(13681,-1321){\makebox(0,0)[lb]{\smash{{\SetFigFont{8}{9.6}{\color[rgb]{0,0,0}$x_0$}%
}}}}
\put(4231,-826){\makebox(0,0)[lb]{\smash{{\SetFigFont{8}{9.6}{\color[rgb]{0,0,0}$\Gamma_s^N$}%
}}}}
\put(11161,-1816){\makebox(0,0)[lb]{\smash{{\SetFigFont{8}{9.6}{\color[rgb]{0,0,0}$\Gamma(0)=\Gamma_{\text{fsi}}(0)$}%
}}}}
\put(4231,-3279){\makebox(0,0)[lb]{\smash{{\SetFigFont{8}{9.6}{\color[rgb]{0,0,0}$\Gamma_f^D$}%
}}}}
\put(1284,-1321){\makebox(0,0)[lb]{\smash{{\SetFigFont{8}{9.6}{\color[rgb]{0,0,0}$\Gamma_s^D$}%
}}}}
\put(1284,-2581){\makebox(0,0)[lb]{\smash{{\SetFigFont{8}{9.6}{\color[rgb]{0,0,0}$\Gamma_f^N$}%
}}}}
\end{picture}%
} 
 \caption{\label{fig.wocontact} Illustration of an FSI problem without contact (left sketch). The domain affiliation in the 
current state can be determined by mapping back to the initial configuration, which is shown on on the right.}
\end{figure}
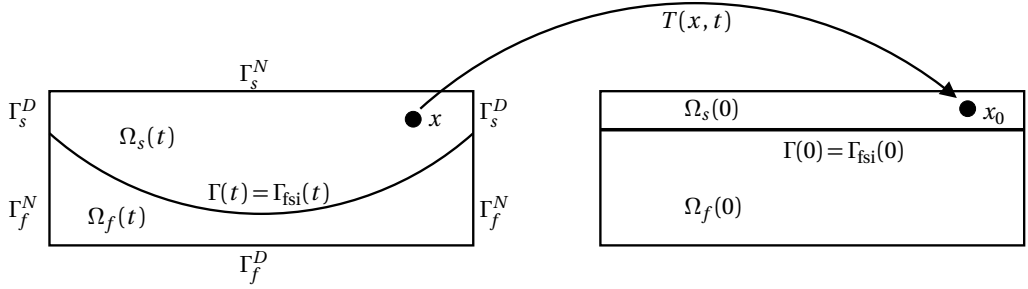

We consider a fluid-structure interaction problem that
is given on an overall domain $\Omega\subset \mathbb{R}^2$ which is split into 
a (variable) sub-domain $\Omega_f(t)$ occupied by a viscous fluid, a
sub-domain $\Omega_s(t)$ occupied by an elastic solid and a
lower-dimensional interface $\Gamma(t)$ separating them, such that
\begin{align*}
 \Omega=\Omega_f(t)\cup\Gamma(t)\cup\Omega_s(t).
\end{align*}
The boundary of the fluid domain is partitioned as follows $\partial \Omega_f(t) = \Gamma_{\rm fsi}(t)\cup \Gamma_{f}^{D}\cup \Gamma_{f}^{N}$, where 
$\Gamma_{\rm fsi}(t)$ stands for the fluid-solid interface. 
As regards the solid boundary, we assume that $\partial \Omega_s(t) = \Gamma(t) \cup \Gamma_{s}^{D} \cup \Gamma_{s}^{N}$, 
where the boundary part $\Gamma(t) = \Gamma_{\rm fsi}(t)\cup \Gamma_C(t)$ is
decomposed in the terms of $\Gamma_{\rm fsi}(t)$ and the contact zone $\Gamma_C(t)$ (see {Figure~\ref{fig.wocontact} for a configuration without contact and the left sketch of Figure~\ref{fig.contact} for a configuration with contact}).
The restriction to two dimensions is made only to simplify the presentation. 
{The models and the methods derived in this paper can be generalised conceptually in a straight-forward way 
to three space dimensions. 
}

{In this work we will use a Eulerian description for the complete FSI problem. As already mentioned in the introduction this is not necessary for 
the contact algorithms derived below, but one convenient way to deal with (possible) topology changes in the fluid domain $\Omega_f(t)$.
In an Eulerian description, the solid sub-domain and the interface are implicitly defined by the (unknown) solid displacement $d$
\begin{align}\label{MovingDomains}
 \Omega_s(t) = \big\{x\in\Omega \, \big| \, T(x,t) \in \Omega_s(0)\big\},\quad
 \Gamma(t) = \big\{x\in\Omega \, \big| \, T(x,t) \in \Gamma(0)\big\},
\end{align}
where $T:\Omega(t)\to\Omega$ is a bijective map, that is given by $T(x,t) =x-d(x,t)$
in the solid domain $\Omega_s(t)$ and by an arbitrary (smooth) extension in $\Omega_f(t)=\Omega\setminus (\Omega_s(t)\cup \Gamma(t))$.
For the details, we refer to the textbook of Richter~\cite{RichterBuch} or Frei~\cite{FreiPhD}.
}

In the variable fluid sub-domain $\Omega_f(t)$, we consider the linear incompressible Stokes equations
\begin{align*}
 \partial_t u - \text{div } \sigma_f(u,p) = f_f,\quad
 \text{div } u =0,
\end{align*}
where the Cauchy stress tensor $\sigma_f$ is defined by
\begin{align*}
 \sigma_f(u,p) = \nu_f {\big( \nabla u^T + \nabla u\big)} - p I,
\end{align*}
$u$ denotes the fluid velocity, $p$ stands for pressure and $\nu_f>0$ is a constant viscosity.
In the solid-subdomain $\Omega_s(t)$, we assume a linear elastic material
\begin{align*}
 \partial_t \dot{d} - \text{div } \sigma_s(d) = f_s,\quad
 \partial_t d = \dot{d},
\end{align*}
where the Cauchy stress tensor $\sigma_s$ is given by
\begin{align*}
 \sigma_s(d) = 2\mu_s E(d) + \lambda_s \text{tr} (E(d)) I, 
 \quad E(d)=\frac{1}{2} \left( \nabla d + \nabla d^T\right),
 \end{align*}
 $\dot{d}=\partial_t d$ denotes the solid velocity
and $\lambda_s,\mu_s>0$ are positive constants.

For the coupling across the {fluid-solid} 
interface $\Gamma_{\text{fsi}}(t)$, the continuity of velocities and
normal stresses 
\begin{align*}
 \dot{d} = u, \quad \sigma_f(u,p) n = \sigma_s(d) n \quad\text{on}\quad\Gamma_{\text{fsi}}(t)
\end{align*}
is typically considered for viscous fluids,
where $n=n_s$ denotes the outer normal vector of the solid domain. 
{We recall that, since in this section there is no contact in the solid (i.e., $\Gamma_C(t)=\emptyset$), we have $\Gamma(t) = \Gamma_{\rm fsi}(t)$}. 
When it comes to 
contact, it is however questionable, whether this condition is still a reasonable approximation
of the underlying physics. Theoretical works~\cite{GerardVaretHillairet} show, that 
the Navier-Stokes equations in combination with these ``no-slip'' conditions do not allow for contact.
Therefore, we will study slip-conditions in this work as well (see Section~\ref{sec.slip}), where the continuity across 
$\Gamma_{\text{fsi}}(t)$ is
only imposed for the normal velocity
\begin{align}
 u \cdot n = \dot{d} \cdot n, \quad \tau^T \sigma_f n = 0 , \quad
 \sigma_f n = \sigma_s n \quad \text{on}\quad\Gamma_{\text{fsi}}(t).\label{slip}
\end{align}

In order to close the system of equations, we define exterior boundary conditions for the fluid 
{and solid}
\begin{align*}
 u&=0 \quad \text{on}\quad \Gamma_f^D,\quad \sigma_f(u,p) n = 0 \quad\text{on}\quad \Gamma_f^N \\
 d&=0 {\quad \text{on}\quad \Gamma_s^D,\qquad \sigma_s(d) n = 0 \quad\text{on}\quad \Gamma_s^N}
\end{align*}
and the initial conditions
\begin{align*}
 u(x,0) = u^0(x) \text{in } \Omega_f(0), 
 \quad d(x,0) = d^0(x), \; \dot{d}(x,0) = \dot{d}^0(x)  \text{in } \Omega_s(0).
\end{align*}

{We introduce the finite element spaces ${\cal V}_h, {\cal Q}_h$ and ${\cal W}_h$ on a quasi-uniform family 
of triangulations $({\cal T}_h)_{h>0}$ and}
use Nitsche's method to combine both equations and interface conditions into 
a monolithic variational formulation
(see e.g. Hansbo \emph{et al.}~\cite{HansboHermanssonSvedberg2004}, Burman \& Fernandez~\cite{BurmanFernandez2014}).
Denoting by $n=n_s$ the outer normal vector of the solid domain $\Omega_s(t)$, the monolithic 
system of equations reads in the no-slip case:
%\begin{varform}{(No-slip conditions)}
\textit{Find $u{(t)}\in {\cal V}_h, p(t) \in {\cal Q}_h, d(t) \in {\cal W}_h$, such that $\dot{d}= \partial_t d$ and}
\begin{align}
 {\cal A}_{\text{fsi, no-slip}} (u,p,d,\dot{d})(v,q,w) = \left(f_f,v\right)_{\Omega_f(t)} + \left(f_s,w\right)_{\Omega_s(t)}\quad 
 \forall v,q,w \in {\cal V}_h \times {\cal Q}_h \times {\cal W}_h,\label{FSI_varForm}
\end{align}
where
\begin{eqnarray}\label{AFSInoslip}
\begin{aligned}
{\cal A}_{\text{fsi, no-slip}} (u,p,&d,\dot{d})(v,q,w)\\
:=&\big(\partial_t u,v\big)_{\Omega_f(t)} + \left(\sigma_f(u,p),\nabla v\right)_{\Omega_f(t)} 
 + \left({\rm div }\, u, q\right)_{\Omega_f(t)}
 + S(p,q)
 +\left(\partial_t \dot{d},w\right)_{\Omega_s(t)}\\
 &\qquad+ \left(\sigma_s(d),\nabla w\right)_{\Omega_s(t)} 
-\left(T_f(u,p,{\dot{d}}) , w-v\right)_{\Gamma(t)} - \left(\dot{d}-u,\sigma_f(v,-q) n\right)_{\Gamma(t)}
 \end{aligned}
\end{eqnarray}
%\end{varform}
{where the numerical fluid traction on the interface is defined by $ T_f(u,p,\dot{d})  :=  \sigma_f(u,p) n - \gamma_{\text{fsi}} ( \dot{d}-u) $. }
{The Nitsche parameter is chosen as $\gamma_{\text{fsi}} := \gamma_{\text{fsi}}^0 \nu_f h^{-1}$.
The term $S(p,q)$ stands for a pressure stabilisation term that is non-zero in case that the discrete fluid spaces do not
fulfil a discrete \textit{inf-sup condition}.}
Note that at the FSI interface, we have used the following relation for the
interface terms arising from integration by parts
\begin{align}
(\sigma_s n,w)_{\Gamma(t)} - (\sigma_f n, v)_{\Gamma(t)} 
= (\sigma_f n, w-v)_{\Gamma(t)} + (\jump{\sigma n},w)_{\Gamma(t)},\label{interfaceterms}
\end{align}
where we have dropped the dependencies of $\sigma_f$ and $\sigma_s$ for better readability.
In the absence of contact, the jump of the stresses defined by 
\begin{equation}\label{eq:jump-stress}
\jump{\sigma n} := \sigma_s n -\sigma_f n,
\end{equation}
vanishes everywhere on $\Gamma(t)$. 
Furthermore, we
have added the term
\begin{align*}
- \left(\dot{d}-u,\sigma_f(v,-q) n\right)_{\Gamma(t)}
\end{align*}
as in Burman \& Fernandez~\cite{BurmanFernandez2014} for stability reasons.

{Using a slip-condition at the FSI interface, the variational formulation reads:
\textit{Find $u\in {\cal V}_h, p \in {\cal Q}_h, d \in {\cal W}_h$, such that $\dot{d}= \partial_t d$ and}
\begin{align*}
 {\cal A}_{\text{fsi, slip}} (u,p,d,\dot{d})(v,q,w) = \left(f_f,v\right)_{\Omega_f(t)} + \left(f_s,w\right)_{\Omega_s(t)}\qquad
 \forall v,q,w \in {\cal V}_h \times {\cal Q}_h \times {\cal W}_h,%\label{FSI_varForm_slip}
\end{align*}
where
\begin{eqnarray}\label{AFSIslip}
\begin{aligned}
&{\cal A}_{\text{fsi, slip}}(u,p,d,\dot{d})(v,q,w)\\
&\quad:=\big(\partial_t u,v\big)_{\Omega_f(t)} + \left(\sigma_f(u,p),\nabla v\right)_{\Omega_f(t)} 
 + \left({\rm div }\, u, q\right)_{\Omega_f(t)} + S(p,q)
 +\left(\partial_t \dot{d},w\right)_{\Omega_s(t)}
 + \left(\sigma_s(d),\nabla w\right)_{\Omega_s(t)} \\ 
 &\qquad\qquad-\left(T_f(u,p,\dot{d}) \cdot n, (w-v)\cdot n\right)_{\Gamma(t)} 
 - \left((\dot{d}-u)\cdot n,n^T \sigma_f(v,-q) n\right)_{\Gamma(t)}.
 \end{aligned}
\end{eqnarray}}

\subsection{Contact model {without fluid}}
\label{sec.ContactModel}

\begin{figure}[bt]
\centering
\begin{minipage}{0.48\textwidth}
\vspace{-0.1cm}
\centering

\resizebox*{0.8\textwidth}{!}{
\begin{picture}(0,0)%
\includegraphics{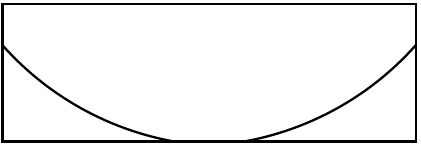}%
\end{picture}%
\setlength{\unitlength}{1450sp}%
\begingroup\makeatletter\ifx\SetFigFont\undefined%
\gdef\SetFigFont#1#2{%
  \fontsize{#1}{#2pt}%
  \selectfont}%
\fi\endgroup%
\begin{picture}(5466,2163)(1768,-3091)
\put(4321,-1321){\makebox(0,0)[lb]{\smash{{\SetFigFont{7}{8.4}{\color[rgb]{0,0,0}$\Omega_s$}%
}}}}
\put(1981,-2581){\makebox(0,0)[lb]{\smash{{\SetFigFont{7}{8.4}{\color[rgb]{0,0,0}$\Omega_f$}%
}}}}
\put(4366,-3031){\makebox(0,0)[lb]{\smash{{\SetFigFont{7}{8.4}{\color[rgb]{0,0,0}$\Gamma_C$}%
}}}}
\put(6211,-3031){\makebox(0,0)[lb]{\smash{{\SetFigFont{7}{8.4}{\color[rgb]{0,0,0}$\Gamma_w$}%
}}}}
\put(5266,-2446){\makebox(0,0)[lb]{\smash{{\SetFigFont{7}{8.4}{\color[rgb]{0,0,0}$\Gamma_{\rm fsi}$}%
}}}}
\end{picture}%
}
\end{minipage}
\hfil
\begin{minipage}{0.48\textwidth}
\centering 
\resizebox*{0.8\textwidth}{!}{
\begin{picture}(0,0)%
\includegraphics{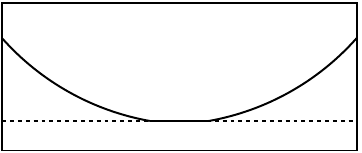}%
\end{picture}%
\setlength{\unitlength}{1243sp}%
\begingroup\makeatletter\ifx\SetFigFont\undefined%
\gdef\SetFigFont#1#2{%
  \fontsize{#1}{#2pt}%
  \selectfont}%
\fi\endgroup%
\begin{picture}(5466,2316)(1768,-3244)
\put(2431,-3076){\makebox(0,0)[lb]{\smash{{\SetFigFont{6}{7.2}{\color[rgb]{0,0,0}$\Omega_f$}%
}}}}
\put(6571,-2671){\makebox(0,0)[lb]{\smash{{\SetFigFont{6}{7.2}{\color[rgb]{0,0,0}$\Gamma_o$}%
}}}}
\put(4321,-1321){\makebox(0,0)[lb]{\smash{{\SetFigFont{6}{7.2}{\color[rgb]{0,0,0}$\Omega_s$}%
}}}}
\put(5491,-2401){\makebox(0,0)[lb]{\smash{{\SetFigFont{6}{7.2}{\color[rgb]{0,0,0}$\Gamma$}%
}}}}
\end{picture}%
}
\end{minipage}
\caption{\label{fig.contact} \textit{Left}: Body in contact with the wall $\Gamma_w$. 
\textit{Right}: Virtual obstacle line $\Gamma_o$ within the fluid domain $\Omega_f$.}
\end{figure}
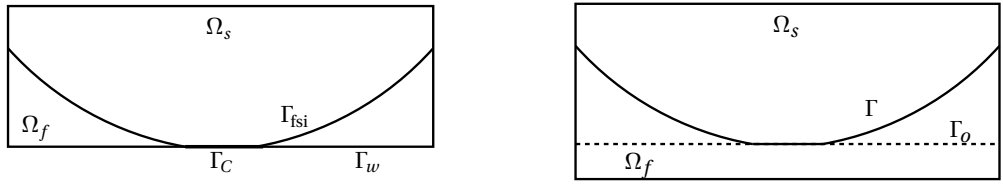

We assume that the solid is at a positive distance to the boundary at initial time and
that contact can only happen with the lower {wall} (see Figure~\ref{fig.contact}, left sketch, where 
the situation at contact is shown)
\begin{align*}
 \Gamma_w=\big\{(x_1,x_2) \in \partial\Omega\, \big| \; x_2=0\big\}.
\end{align*}
We denote the outer normal vector of the fluid domain $\Omega_f(t)$ at $\Gamma_w$
by $n_w= -e_2$. Moreover, let $g_0(x_2)>0$ be the
function describing the initial distance of a point $(x_1,x_2) \in \Gamma(t)$ to the {wall} $\Gamma_w$.

{When contact with $\Gamma_{\rm w}$ occurs on a part $\Gamma_C(t) \subset \Gamma(t)$}, suitable contact conditions are
(Alart \& Curnier~\cite{AlartCurnier91}, 
Chouly \& Hild~\cite{ChoulyHild2013}) 
\begin{align}
 d\cdot n_w\leq g_0, \quad \sigma_{s,n} &:= n_w^T \sigma_s n \leq 0, 
 \quad (d\cdot n_w-g_0) \sigma_{s,n} = 0 \quad \text{on}\quad \Gamma(t).\label{normalContact}
\end{align}

The first inequality in (\ref{normalContact}) ensures that the solid can not pass though $\Gamma_w$,
the second inequality describes that the normal stress is zero (in the absence of contact) or negative (during contact) and the third
condition is a complementarity condition that guarantees that at least one of the inequalities 
is ``active''.

In the model setting, that is considered here, the normal and tangential 
vectors $n_w$ and $\tau_w$ are constant. For later purposes we distinguish between the normal and tangential vectors of 
the wall $n_w, \tau_w$ and the normal and tangential vector of the FSI interface $n,\tau$. The vectors are 
equal on $\Gamma_C(t)$, 
where the body is in contact with the wall, but different before contact.

For arbitrary $\gamma_C>0$ the first line in (\ref{normalContact}) is 
equivalent to~\cite{ChoulyHild2013, AlartCurnier91} 
\begin{align}
  \sigma_{s,n}(d) = -\gamma_C [d\cdot n_w-g_0 -\frac{1}{\gamma_C} \sigma_{s,n}(d)]_+ 
  =: - \gamma_C [P_{\gamma,{s}}(d)]_+\label{ContactFormula},
\end{align} 
where $[f]_+:= \max\{f,0\}$.
The equivalence can be shown by simple calculations for each of the 
cases $P_{\gamma,s}(d)\leq0$ and $P_{\gamma,s}(d)>0$, see Chouly \& Hild~\cite{ChoulyHild2013}.

As tangential contact condition, 
{Chouly \& Hild\cite{ChoulyHild2013} used} the ``no-friction`` condition
\begin{align}
 \tau_w^T \sigma_s n &= 0 \quad \text{on}\quad \Gamma_C(t) \label{tangentialContact}
\end{align}
for the pure solid problem. Choosing $\gamma_C=\gamma_C^0\mu_s h^{-1}$, the variational formulation reads:
\textit{Find $d(t) \in {\cal W}_h$ such that $\dot{d} = \partial_t d$ and}
\begin{equation}\label{theta0}
 \left(\partial_t \dot{d},w\right)_{\Omega_s(t)}
 + \left(\sigma_s(d),\nabla w\right)_{\Omega_s(t)} 
 +\gamma_C \left([P_{\gamma,s}(d)]_+, w\cdot n_w\right)_{{\Gamma(t)}}
 = \left(f_s,w\right)_{\Omega_s(t)}\quad 
 \forall w \in {\cal W}_h.
\end{equation}

{We will discuss the tangential contact conditions for the case of fluid-structure interaction with contact below.}

We close this section by mentioning that Chouly, Hild \& Renard~\cite{ChoulySymVsNonsym} proposed
a more general contact formulation that makes use 
of the consistency of the term
\begin{align*}
\left(\gamma_C [P_{\gamma,s}(d)]_+ + \sigma_{s,n}(d), \sigma_{s,n}(w)\right)_{\Gamma(t)}.
\end{align*}
For $\theta \in [-1,1]$, the contact term can be generalised to
\begin{eqnarray}
 \begin{aligned}\label{thetasolid}
 \gamma_C \big(&[P_{\gamma,s}(d)]_+, w\cdot n_w\big)_{\Gamma(t)}
 -\theta \left(\gamma_C [P_{\gamma,s}(d)]_+ + \sigma_{s,n}(d), \sigma_{s,n}(w)\right)_{\Gamma(t)}\\
 &\qquad= (1-\theta) \gamma_C \left([P_{\gamma,s}(d)]_+, w\cdot n_w\right)_{\Gamma(t)}
 +\theta \left(\gamma_C [P_{\gamma,s}(d)]_+, P_{\gamma,s}(w)\right)_{\Gamma(t)}
 - \theta \left(\sigma_{s,n}(d), \sigma_{s,n}(w)\right)_{\Gamma(t)}.
 \end{aligned}
 \end{eqnarray}

For $\theta=0$, we recover the formulation (\ref{theta0}). Besides that, the case $\theta=1$ 
is of particular interest, as it yields a symmetric formulation, for which
a stability result has been shown~\cite{ChoulySymVsNonsym}.

\section{Fluid-structure interaction with contact} 
\label{sec.FSIcontact}
{
We now address the case where both fluid-structure interaction and contact occur, by combining 
the ideas described in Sections~\ref{sec.FSI}-\ref{sec.ContactModel}. In the next section, we consider first a simplified setting 
in which the solid can enter into contact with a virtual obstacle invisible to the fluid. This simplification allows us to consider the numerical treatment of 
the interface conditions without accounting for the issues related to topology changes within the fluid. Then, in Section~\ref{sec.fullContact}, we show 
how this numerical setting can be extended to model physically realistic contact and fluid-structure interaction. 
}
\subsection{Virtual obstacle within the fluid domain}
\label{sec.ObstModel}

{In order to simplify the presentation we introduce the combined 
FSI-contact formulation first for a model problem with a horizontal obstacle line $\Gamma_o$ within the fluid domain
(see the right sketch of Figure~\ref{fig.contact}), 
before we discuss the full FSI-contact problem including topology changes.
The obstacle is invisible to the fluid, but an obstacle to the solid. One may consider 
for example a membrane that is perfectly permeable for the fluid, 
but not for the solid or a magnetic field below the obstacle that prevents the solid from crossing the line. We assume for simplicity that 
the obstacle line is parallel to the fluid boundary and denote its distance by $\alpha>0$. The initial distance of the FSI interface $\Gamma(0)$
to the obstacle line is then given by $g_{\alpha}(x_2) = g_0(x_2) - \alpha$. We assume that $g_\alpha \geq 0$. 
}

When considering fluid-structure interactions, 
the body is pre-loaded before 
it reaches the obstacle by means of the balance of normal forces
\begin{align}
 {\sigma_{s}n = \sigma_{f}n}\quad \text{on } \Gamma(t).\label{BalNormalForces}
\end{align}
{If the interface $\Gamma(t)$ reaches the obstacle $\Gamma_o$, 
the additional constraint $d\cdot n_w \leq g_{\alpha}$ has to be fulfilled. 
{This gives rise to an additional 
surface force, that acts in the direction $-n_w=e_2$ (normal to $\Gamma_o$)
\begin{align} 
  \jump{\sigma n}  {-} \lambda n_w = 0 \quad \text{on}\quad \Gamma(t).\label{LagrComplete}
\end{align}
The variable $\lambda$ defined through (\ref{LagrComplete}) can be seen as a Lagrange multiplier, which is zero in absence of contact with the obstacle line
(due to \eqref{BalNormalForces}) and can become negative during contact. This is described by the complementarity conditions}
\begin{align}\label{combCont}
 d\cdot n_w\leq g_{\alpha}, \quad \lambda &\leq 0, 
 \quad (d\cdot n_w-g_{\alpha}) \lambda = 0.
\end{align}

\noindent Splitting into normal and tangential contributions, \eqref{LagrComplete} reads
\begin{align}
n_w^T \jump{\sigma n}  - \lambda &= 0, \quad\tau_w^T \jump{\sigma n}  = 0 \quad \text{on}\quad \Gamma(t).\label{DefLagr}
\end{align}
} 

{
\noindent Using the same trick as in Section~\ref{sec.ContactModel}, (\ref{combCont}) is equivalent to 
\begin{align}
 \lambda = -\gamma_C \left[d\cdot n_w -g_{\alpha} -\frac{1}{\gamma_C} \lambda\right]_+ 
 =: -\gamma_C [P_{\gamma}(\lambda,d)]_+ \quad \text{on}\quad \Gamma(t).\label{ContactFormula2}
\end{align}
Eliminating the Lagrange multiplier by using \eqref{DefLagr}, this reads 
\begin{align} 
 {\lambda}= -\gamma_C \left[d\cdot n_w -g_{\alpha} -\frac{1}{\gamma_C} \jump{\sigma_n(u,p,d)}\right]_+ 
 =: -\gamma_C [P_{\gamma}(\jump{\sigma_n(u,p,d)},d)]_+ \quad \text{on} \quad \Gamma(t),\label{ContForce1}
\end{align}
where we have used the abbreviation $\jump{\sigma_n} := n_w^T \jump{\sigma n} $.

{The natural formulation in the discrete setting is to consider the numerical stress $T_f$ in \eqref{eq:jump-stress} and in \eqref{LagrComplete}}. Let us derive these conditions first for the no-slip case. Adding
the additional surface force to the FSI-Nitsche formulation \eqref{FSI_varForm}, the discrete variational formulation reads
\begin{align}
 {\cal A}_{\text{FSI, no-slip}}(u,p,d,\dot{d})(v,q,w) -(\lambda n_w,w)_{\Gamma(t)} = \left(f_f,v\right)_{\Omega_f(t)} + \left(f_s,w\right)_{\Omega_s(t)}\quad
 \forall v,q,w \in {\cal V}_h \times {\cal Q}_h \times {\cal W}_h.\label{FSILagrMult}
\end{align}
Due to the additional Nitsche interface terms in \eqref{LagrComplete}, this formulation 
includes the interface condition
\begin{align}
 \jump{\widetilde{\sigma}_n(u,p,d)}  - \lambda =0 \quad \text{on} \quad \Gamma(t),\label{LagrFluxes}
\end{align}
where the numerical stress jump in the FSI-Nitsche formulation
across the interface is given by 
\begin{equation}\label{eq:jump-numerical-stress}
\jump{\widetilde{\sigma}_n(u,p,d)} := \sigma_{s, n}(d) - T_f(u,p,\dot{d})\cdot n_{\rm w} .
\end{equation}
 
\noindent Eliminating the Lagrange multiplier by means of \eqref{LagrFluxes}, the identity \eqref{ContactFormula2} reads
\begin{align}
 %\jump{\widetilde{\sigma}_n(u,p,d)}
 {\lambda }
 = -\gamma_C \left[d\cdot n_w -g_{\alpha} -\frac{1}{\gamma_C} \jump{\widetilde{\sigma}_n }\right]_+ 
 =: -\gamma_C [P_{\gamma}(\jump{\widetilde{\sigma}_n(u,p,d)},d)]_+ \quad \text{on} \quad \Gamma(t).\label{ContForce2}
\end{align}
For the definition of numerical stress in the slightly more complicated case of slip-interface conditions, we refer to Section~\ref{sec.slip}. 
We will in the following analyse both contact formulations \eqref{ContForce1} and \eqref{ContForce2}
and in particular what the different possibilities 
mean with respect to the weakly imposed interface conditions. In order to avoid too much repetition, we use a general formulation
that includes the Lagrange multiplier $\lambda=\lambda(u,p,d)$, keeping in mind that $\lambda(u,p,d)$ will be either chosen 
as the jump of normal stresses $\jump{\sigma_n(u,p,d)}$ or the jump of 
normal discrete stresses $\jump{\widetilde{\sigma}_n(u,p,d)}$.}

Using \eqref{ContactFormula2}, the variational formulation reads: \begin{varform}\label{varForm_implicit}
\textit{Find $u(t)\in {\cal V}_h, p(t) \in {\cal Q}_h, d(t)\in {\cal W}_h$ such that $\dot{d}=\partial_t d$ and }
\begin{eqnarray}\label{vF_implicit}
 \begin{aligned}
 {\cal A}_{\text{fsi},*} (u,p,d,\dot{d})(v,q,w)
+\gamma_C &\left([P_{\gamma}(\lambda,d)]_+, w\cdot n_w\right)_{\Gamma(t)}\\
 &\qquad\qquad= \left(f_f,v\right)_{\Omega_f(t)} + \left(f_s,w\right)_{\Omega_s(t)}\quad 
 \forall v,q,w \in {\cal V}_h \times {\cal Q}_h \times {\cal W}_h,
\end{aligned}
\end{eqnarray}
\end{varform}
\noindent where the bilinear form is one of the forms ${\cal A}_{\text{fsi, no-slip}}$ or ${\cal A}_{\text{fsi, slip}}$ defined in
\eqref{AFSInoslip} and \eqref{AFSIslip}, respectively, and the contact parameter is chosen~\cite{ChoulyHild2013} $\gamma_C= \gamma_C^0 \mu_s h^{-1}$.

{\subsubsection{Weakly imposed interface conditions}
\label{sec.ImplCond}

Let us now analyse which interface conditions on $\Gamma(t)$ are implicitly included in Variational Formulation~\ref{varForm_implicit}. 
For simplicity, we analyse the formulation with the bilinear form ${\cal A}_{\text{fsi, no-slip}}$ corresponding to no-slip conditions. }
Therefore, we integrate by parts in \eqref{vF_implicit} and consider only the interface terms {by formally neglecting the bulk and inter-element terms}.
For better readability, we drop all the dependencies of $\sigma_f$ and $\sigma_s$.
Collecting all terms with the fluid test function $v$, we obtain
\begin{align}\label{intCondFluid}
 -\sigma_f n + \sigma_f n -\gamma_{\text{fsi}} (\dot{d} - u) = -\gamma_{\text{fsi}} (\dot{d} - u) = 0 \quad \mbox{on}\quad \Gamma(t),
\end{align}
i.e.$\,$the kinematic condition $\dot{d}=u$. Next, we collect the interface terms for the solid part $w$ and split into
a normal part ($w\cdot n_w$) and a tangential part ($w\cdot \tau_w$). {We recall that since the boundary $\Gamma_{\rm w}$ is 
flat, the extension to $\Gamma(t)$ of its tangential and normal vectors are trivial}. For the tangential part, we obtain as usual 
for Nitsche-based FSI
\begin{align}
 \tau_w^T \sigma_s n- \tau_w^T\sigma_f n +\gamma_{\text{fsi}} (\dot{d} - u)\cdot \tau_w = 0\quad {\mbox{on}\quad \Gamma(t)}.\label{ImplTang}
\end{align}

\noindent For the normal part, we have
\begin{align*}
 \sigma_{s,n} - \sigma_{f,n} +\gamma_{\text{fsi}} (\dot{d} - u)\cdot n_w + \gamma_C [P_{\gamma}(d,\lambda)]_+ = 0 \quad \mbox{on}\quad \Gamma(t).
\end{align*}
Let us first consider the case that the contact force is not active. We obtain, as in the standard
FSI-Nitsche formulation
\begin{align*}
 \jump{\sigma_n} + \gamma_{\text{fsi}} (\dot{d}-u)\cdot n_w = 0. 
\end{align*}
If the contact force is active, we get
\begin{eqnarray}
\begin{aligned}
 0 =\jump{\sigma_n} + \gamma_{\text{fsi}} (\dot{d}-u)\cdot n_w &+ \gamma_C (d\cdot n_w-g_{\alpha}) -\lambda 
 = \begin{cases}
 \gamma_C (d\cdot n_w-g_{\alpha}) + \gamma_{\text{fsi}} (\dot{d}-u)\cdot n_w, \quad &\lambda = \jump{\sigma_n}\\
 \gamma_C (d\cdot n_w-g_{\alpha}) , \quad &\lambda = \jump{\widetilde{\sigma}_n}
 \end{cases}
 \label{mixture_wo_dGfluxes}
\end{aligned}
\end{eqnarray}
In the first case, this is a combination of the ''active`` contact condition $d\cdot n_w=g$ and the continuity of velocities. As the continuity of velocities 
is imposed from the fluid side~\eqref{intCondFluid}, this is not an issue for the model problem considered here. It will however lead to problems, when we consider
contact of the solid with the lower wall $\Gamma_w$ in Section~\ref{sec.fullContact}. There, the second formulation, based on the {discrete stresses~\eqref{LagrFluxes}} will be needed. As can be 
seen in \eqref{mixture_wo_dGfluxes}, the pure contact condition $d\cdot n_w=g$ is valid from the solid side during contact.

\subsection{Contact with the boundary of the fluid domain}
\label{sec.fullContact}

In this section, we will {derive} two numerical approaches for the full FSI-contact problem, considering contact with the lower fluid boundary $\Gamma_w$.
First, we note that the formulation derived in the previous section has a simple extension to the case of contact with the fluid boundary, if we relax the contact
condition by a small $\epsilon>0$. see Section~\ref{sec.relaxed}. {Then, we introduce an artificial fluid below the contact line in Section~\ref{sec.artf},
where we drive the fluid velocity to zero by means of a volume penalty approach. This enables us to include the FSI interface and contact conditions fully implicitly
without relaxation.}

{
 \subsubsection{Relaxed contact formulation}
\label{sec.relaxed}

 {The idea of the relaxed formulation is to {place} the virtual obstacle of 
 Section~\ref{sec.ObstModel} at a distance $\alpha =\epsilon(h)$ from $\Gamma_w$, 
 see Figure~\ref{fig.contact1vs2} on the left.}
We assume that $\epsilon(h) \to 0$, as the mesh size $h$ tends to zero. Now we can use the contact formulation 
derived for an obstacle within the fluid domain above using the 
function $g_{\epsilon(h)}(x_2) = g_0(x_2) - \epsilon(h)$. 
{The corresponding variational formulation is Variational Formulation~\ref{varForm_implicit} using
$g_{\epsilon(h)}$ in the definition of $P_{\gamma}$ \eqref{ContForce1} instead of $g_{\alpha}$.
}

Besides its simplicity, the main advantages of this contact formulation are:
\begin{itemize}
 \item The numerical difficulties related with a topology change of the fluid domain are avoided.
 \item No-slip conditions can be used on both $\Gamma(t)$ and the lower wall $\Gamma_w$.
\end{itemize}

On the other hand, from a modelling point of view the contact conditions {\eqref{LagrComplete}-\eqref{combCont} are first of all 
only justified}, if we assume that an infinitesimal fluid layer remains 
between the fluid-structure-interface $\Gamma(t)$ and the 
lower wall $\Gamma_w$.
Then the fluid stresses $\sigma_f$ appearing in the jump terms $\jump{\sigma n}$ remain 
well-defined during contact. {This assumption} might be justified, if no-slip conditions are used on
interface $\Gamma(t)$ and boundary $\Gamma_w$, see the theoretical works 
of e.g.$\,$Gerart-Varet \& Hillairet~\cite{GerardVaretHillairet}. {In this framework, the choice 
of $\lambda$ in terms of the physical or numerical stresses appears not to be essential.}

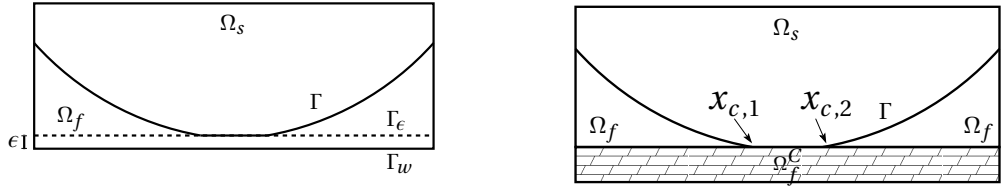
\begin{figure}[bt]
\centering

\begin{minipage}{0.48\textwidth}
\vspace{-0.2cm}
\centering

\resizebox*{0.8\textwidth}{!}{
\begin{picture}(0,0)%
\includegraphics{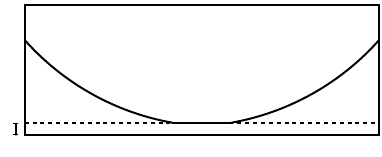}%
\end{picture}%
\setlength{\unitlength}{1243sp}%
\begingroup\makeatletter\ifx\SetFigFont\undefined%
\gdef\SetFigFont#1#2{%
  \fontsize{#1}{#2pt}%
  \selectfont}%
\fi\endgroup%
\begin{picture}(5808,2343)(1426,-3271)
\put(4321,-1321){\makebox(0,0)[lb]{\smash{{\SetFigFont{6}{7.2}{\color[rgb]{0,0,0}$\Omega_s$}%
}}}}
\put(2116,-2581){\makebox(0,0)[lb]{\smash{{\SetFigFont{6}{7.2}{\color[rgb]{0,0,0}$\Omega_f$}%
}}}}
\put(5536,-2401){\makebox(0,0)[lb]{\smash{{\SetFigFont{6}{7.2}{\color[rgb]{0,0,0}$\Gamma$}%
}}}}
\put(1441,-2919){\makebox(0,0)[lb]{\smash{{\SetFigFont{6}{7.2}{\color[rgb]{0,0,0}$\epsilon$}%
}}}}
\put(6571,-2626){\makebox(0,0)[lb]{\smash{{\SetFigFont{6}{7.2}{\color[rgb]{0,0,0}$\Gamma_{\epsilon}$}%
}}}}
\put(6571,-3211){\makebox(0,0)[lb]{\smash{{\SetFigFont{6}{7.2}{\color[rgb]{0,0,0}$\Gamma_w$}%
}}}}
\end{picture}%
}
\end{minipage}
\hfil
\begin{minipage}{0.48\textwidth}
\centering
\resizebox*{0.8\textwidth}{!}{
\begin{picture}(0,0)%
\includegraphics{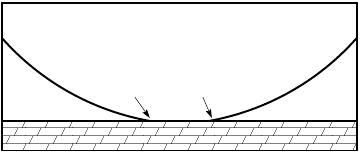}%
\end{picture}%
\setlength{\unitlength}{1243sp}%
\begingroup\makeatletter\ifx\SetFigFont\undefined%
\gdef\SetFigFont#1#2{%
  \fontsize{#1}{#2pt}%
  \selectfont}%
\fi\endgroup%
\begin{picture}(5466,2316)(1768,-3244)
\put(4321,-1321){\makebox(0,0)[lb]{\smash{{\SetFigFont{6}{7.2}{\color[rgb]{0,0,0}$\Omega_s$}%
}}}}
\put(1981,-2581){\makebox(0,0)[lb]{\smash{{\SetFigFont{6}{7.2}{\color[rgb]{0,0,0}$\Omega_f$}%
}}}}
\put(6751,-2581){\makebox(0,0)[lb]{\smash{{\SetFigFont{6}{7.2}{\color[rgb]{0,0,0}$\Omega_f$}%
}}}}
\put(5671,-2401){\makebox(0,0)[lb]{\smash{{\SetFigFont{6}{7.2}{\color[rgb]{0,0,0}$\Gamma$}%
}}}}
\put(4321,-3054){\makebox(0,0)[lb]{\smash{{\SetFigFont{5}{6.0}{\color[rgb]{0,0,0}$\Omega_f^C$}%
}}}}
\put(3466,-2266){\makebox(0,0)[lb]{\smash{{\SetFigFont{5}{6.0}{\color[rgb]{0,0,0}{\small $x_{c,1}$}}%
}}}}
\put(4636,-2266){\makebox(0,0)[lb]{\smash{{\SetFigFont{5}{6.0}{\color[rgb]{0,0,0}{\small $x_{c,2}$}}%
}}}}
\end{picture}%
}
\end{minipage}
\caption{\label{fig.contact1vs2} Illustration of the two possible approaches to include contact. \textit{Left}: Relaxed contact formulation: Body in contact with an 
obstacle line $\Gamma_{\epsilon}$ close to $\Gamma_w$. 
\textit{Right}: Introduction of an artificial fluid domain $\Omega_f^C$ below $\Gamma_w$.}
\end{figure}
}

\subsubsection{Contact formulation using an artificial fluid}
\label{sec.artf}

The idea of this second approach is to add an artificial fluid domain $\Omega_f^C$ below $\Gamma_w$, see 
Figure~\ref{fig.contact1vs2}, right sketch. In $\Omega_f^C$, we will drive the velocity $u$ to zero by means of a penalty term.

The variational formulation reads: \begin{varform}\label{varForm_artF}
\textit{Find $u\in {\cal V}_h, p \in {\cal Q}_h, d \in {\cal W}_h$ such that $\dot{d}=\partial_t d$ and }
\begin{eqnarray}\label{vF_artF}
 \begin{aligned}
 {\cal A}_{FSI,*}^C (u,p,d,\dot{d})(v,q,w)
+\gamma_C &\left([P_\gamma(\lambda,d)]_+, w\cdot n_w\right)_{\Gamma(t)}\\
 &= \left(f_f,v\right)_{\Omega_f(t)} + \left(f_s,w\right)_{\Omega_s(t)}\quad 
 \forall v,q,w \in {\cal V}_h \times {\cal Q}_h \times {\cal W}_h,
\end{aligned}
\end{eqnarray}
\end{varform}
\noindent where in the no-slip case
\begin{eqnarray}
 \begin{aligned}
 {\cal A}_{\text{fsi,no-slip}}^C (u,p,d,\dot{d})(v,q,w)&:=
 \big(\partial_t u,v\big)_{\Omega_f(t)\cup \Omega_f^C} + \left(\sigma_f(u,p),\nabla v\right)_{\Omega_f(t)\cup\Omega_f^C} 
 + \left(\text{div } u, q\right)_{\Omega_f(t)\cup\Omega_f^C} \\
 &\qquad+ S(p,q) 
+ \gamma_a (u,v)_{\Omega_f^C}+\left(\partial_t \dot{d},w\right)_{\Omega_s(t)}
 + \left(\sigma_s(d),\nabla w\right)_{\Omega_s(t)}\\
 &\qquad\qquad
 %+ \left(\partial_t d,z\right)_{\Omega_s(t)} - \left(\dot{d},z\right)_{\Omega_s(t)}\\
 - \left(T_f(u,p,\dot{d}) , w-v\right)_{\Gamma(t)} 
 %+\gamma_{FSI} \left(\dot{d}-u, w-v\right)_{\Gamma(t)} -
 -\left(\dot{d}-u,\sigma_f(v,-q) n\right)_{\Gamma(t)}\\
 \end{aligned}
\end{eqnarray}
and for slip interface conditions
\begin{eqnarray}
 \begin{aligned}
 {\cal A}_{\text{fsi,slip}}^C (u,p,d,\dot{d})(v,q,w)&:=
 \big(\partial_t u,v\big)_{\Omega_f(t)\cup \Omega_f^C} + \left(\sigma_f(u,p),\nabla v\right)_{\Omega_f(t)\cup\Omega_f^C} 
 + \left(\text{div } u, q\right)_{\Omega_f(t)\cup\Omega_f^C}\\ 
 &\quad+ S(p,q)
 + \gamma_a (u,v)_{\Omega_f^C}
 +\left(\partial_t \dot{d},w\right)_{\Omega_s(t)}
 + \left(\sigma_s(d),\nabla w\right)_{\Omega_s(t)}\\
 %+ \left(\partial_t d,z\right)_{\Omega_s(t)} - \left(\dot{d},z\right)_{\Omega_s(t)}\\
 &\qquad- \left(T_f(u,p,\dot{d}) \cdot n, (w-v)\cdot n\right)_{\Gamma(t)}
 - \left((\dot{d}-u)\cdot n,n^T \sigma_f(v,-q) n\right)_{\Gamma(t)}.\\
 \end{aligned}
\end{eqnarray}
{Note the presence of the penalty term $\gamma_a (u,v)_{\Omega_f^C}$ within the artificial fluid, where $\gamma_a:=\gamma_a^0 h^{-2}$.

{\begin{remark}[(Porous medium analogy)]
 The penalisation used in Variational Formulation~\ref{varForm_artF} corresponds to the so-called penalty approach 
 that is sometimes used for the coupling of free flow and flow through porous medium~\cite{CimolinDiscacciati, IlievLaptev}. There,
 the Stokes equations in the fluid part $\Omega_f$ and the Darcy equations in the porous medium $\Omega_p$ (which corresponds to $\Omega_f^C$) are formulated simultaneously 
 in the whole domain $\Omega=\Omega_f \cup \Omega_p$ in the spirit of the volume penalty approach
 \begin{align*}
  \partial_t u -\mu \Delta u + \nabla p + \frac{\mu}{K} u \chi_{\Omega_p} = 0, \qquad \nabla \cdot u = 0 \quad \text{ in } \Omega,
 \end{align*}
 where $K$ denotes the permeability of the porous medium and $\chi_{\Omega_p}$ is the characteristic function of the domain $\Omega_p$. 
 In this sense 
 the artificial fluid in our approach can be seen as a porous medium with asymptotically vanishing permeability $K={\cal O}(h^2)$. A mathematical
 justification of this penalisation has been given by Angot~\cite{Angot1999}.
\end{remark}}

\paragraph{Discussion of the weakly imposed interface conditions: The no-slip case}

We have already derived the weakly imposed interface conditions for the no-slip case in Section~\ref{sec.ImplCond}. Here, however, we have to 
consider that the fluid below the contact line is artificial, and we should in particular make sure that there is no feedback from the
artificial fluid to the solid. {In other words, we want that the artificial fluid acts as a \textit{slave} to the solid during contact.}
Owing to \eqref{mixture_wo_dGfluxes}, this naturally motivates the choice of $\lambda$ in terms of the numerical stress. 

First, we note that the continuity of velocities~\eqref{intCondFluid} is imposed from the fluid side, such that no feedback to 
the solid is included from this equation. Considering the contact condition for the normal contact~\eqref{mixture_wo_dGfluxes}, 
we obtain
\begin{align}
 \gamma_C (d\cdot n_w-g_0) + \gamma_{\text{fsi}} (\dot{d}-u)\cdot n_w =0,
\end{align}
when choosing $\lambda =\jump{\sigma_n }$. Instead of the condition $d\cdot n_w=g_0$, this induces an influence 
from the artificial velocity $u$ from $\Omega_f^C$
onto the solid displacement. Moreover, if $u$ is driven to zero in $\Omega_f^C$, 
$\dot{d} \cdot n_w$ goes to zero as well, which might prevent the body from releasing from contact. 
On the other hand, using 
the jump of fluxes $\lambda = \jump{\widetilde{\sigma}_n} $, we obtain the ''pure`` contact condition 
\begin{align}
 d\cdot n_w=g_0 \label{dGflux_argumentation}
\end{align}
 as desired. 

The {weakly imposed} tangential contact condition \eqref{ImplTang} reads
\begin{align}
 \tau_w^T \sigma_s n +\gamma_{\text{fsi}} \dot{d}\cdot \tau_w = 0,\label{penalty}
\end{align}
when considering that the fluid velocity $u$ is driven to zero asymptotically. As $\gamma_{\text{fsi}}\to\infty$ for $h\to 0$,
this means (asymptotically) that the solid is not allowed to slide along the line $\Gamma_w$. While this might seem restrictive at first sight,
this condition is in fact in some sense inherited from
the no-slip condition at $\Gamma(t)$ before contact. This is due to the continuity of velocities on $\Gamma_{\text{fsi}}(t)$
and the fact that the velocity is driven to zero in $\Omega_f^C$ (and hence on $\Gamma_w$). Moreover, the no-slip condition on the fluid part
of $\Gamma_w$ implies that the solid can not slide on the end points $x_{C,1}$ and $x_{C,2}$ of the contact interval (see Figure~\ref{fig.contact1vs2}, right sketch).
Altogether, this shows in agreement with a number of theoretical works (e.g. Gerard-Varet et al~\cite{GerardVaretetal2015})
that the no-slip interface conditions are not an appropriate model for the case that it comes to contact with an exterior wall.
}

{
\begin{remark}[(Relation of the two contact formulations)]
 The relaxed contact formulation derived in Section~\ref{sec.relaxed} can also be seen as an extension of the fluid forces $\sigma_f$
 to a region below the contact line (here $\Gamma_\epsilon$), namely by using the Stokes equations in the extended domain. We will see in the numerical examples 
 below that the two approaches yield similar results. 
 In this way the use of the relaxed contact formulation might be justified also in cases, where real contact with the wall is expected. Moreover, it would 
 be enough to use the extension in the artificial fluid approach only in a small layer of size ${\cal O}(\epsilon)$.
\end{remark}

\begin{remark}[(Lagrange multiplier formulation)]
 A further possibility would be to keep the Lagrange multiplier $\lambda$ in the variational formulation \eqref{FSILagrMult}-\eqref{LagrFluxes} and 
 to discretise additionally the Lagrange multiplier space,
 as in~\cite{LegayChessaBelytschko2006, GerstenbergerWall, Kamenskyetal2017} for the pure FSI case. Due to the difficulties concerning the discrete \textit{inf-sup stability}
 and the additional computational effort, we will, however, not consider this alternative in the remainder of this work. 
\end{remark}
}

\subsection{Slip conditions}
\label{sec.slip}

Motivated by the above considerations, we study slip conditions in this section. {The arguments that follow can be applied to the 
formulations proposed in Sections~\ref{sec.ObstModel} and \ref{sec.fullContact}. }
As mentioned in the introduction,
Gerard-Varet, Hillairet and Wang showed that contact
can not happen for a smooth rigid body falling down to the ground by means of gravity,
if the no-slip boundary conditions are used for the velocities on both the fluid-structure 
interface $\Gamma_{\text{fsi}}(t)$ and
the contact {wall} $\Gamma_w$, but that contact 
happens in this example, if slip boundary conditions are used on both $\Gamma_{\text{fsi}}(t)$ and
$\Gamma_w$~\cite{GerardVaretetal2015}.

{Observe that in the continuous case the discussion of \eqref{BalNormalForces}-\eqref{LagrComplete} remains valid in the case of contact with slip conditions (noting that the tangential 
stresses vanish on both sides of the fluid-solid interfaces).
However, at the discrete level, care has to be taken to use only the normal 
component of the numerical stress. Indeed, in this case the relation \eqref{LagrComplete}
becomes 
 \begin{align} 
  {(n^T\jump{\sigma n}) n}  - \lambda n_w = 0 \quad \text{on}\quad \Gamma(t),\label{LagrCompleteSlip}
\end{align}
which yields {$\lambda = (n^T\jump{\sigma n}) (n \cdot n_{\rm w})=:\jump{\sigma_{n,\text{slip}}}$}. 

At the discrete level, using the numerical stress this expression translates 
to
\begin{align}
\lambda = {(n^T\jump{\widetilde \sigma n})} (n \cdot n_{\rm w}) =: \jump{\widetilde\sigma_{n,\text{slip}}},\label{lambdaslip}
\end{align}
where the jump operator is given by \eqref{eq:jump-numerical-stress}. 

The resulting discrete {formulations are given by Variational 
Formulation \ref{varForm_implicit} or Variational 
Formulation \ref{varForm_artF}, respectively, with the respective choice} of $\lambda$.
}

\paragraph{Weakly imposed interface conditions: The slip case}
%\label{subsec.InterfacedG}

Let us consider again which interface conditions are implicitly included in the variational
formulation,{ when using the numerical stresses \eqref{lambdaslip} for $\lambda$}. 
Considering the interface terms with fluid test function $v$
yields as usual the interface conditions $\tau^T \sigma_f n=0$ and $(\dot{d}-u)\cdot n=0$. 
Let us therefore concentrate on
the terms with the solid test function $w$
\begin{align*}
 (\sigma_s n, w)_{\Gamma(t)} - (T_f \cdot n, w\cdot n)_{\Gamma(t)} 
 +\gamma_C([\widetilde{P}_{\gamma}]_+, w\cdot n_w)_{\Gamma(t)}=0.
\end{align*}
{In the case without contact, the last term vanishes and hence we retrieve the standard 
consistency of Nitsche's method for fluid-structure interaction with slip conditions. 
On the other hand, by developing the solid test function $w$ in the local basis of $\Gamma(t)$ we have 
\begin{align*}
( \tau^T \sigma_{s}n, w\cdot \tau)_{\Gamma(t)} +
(  (n^T\jump{\widetilde \sigma n}) , w\cdot n)_{\Gamma(t)} 
 +\gamma_C([P_{\gamma}(d,\lambda)]_+, w\cdot n_w)_{\Gamma(t)}=0.
\end{align*}
Now, we use the identity $ n = (n\cdot n_w)n_w + (n\cdot \tau_w)\tau_w$ to get 
\begin{align*}
( \tau^T \sigma_{s}n, w\cdot \tau)_{\Gamma(t)} +
(  (n^T\jump{\widetilde \sigma n}) , (w\cdot \tau_w) ( n\cdot \tau_w))_{\Gamma(t)} +
(  (n^T\jump{\widetilde \sigma n}) , (&w\cdot n_w) ( n\cdot n_w))_{\Gamma(t)} \\
 &+\gamma_C([P_{\gamma}(d,\lambda)]_+, w\cdot n_w)_{\Gamma(t)}=0.
\end{align*}
Hence, in the case of contact by using the definition of $P_{\gamma}$, we get 
\begin{align}\label{eq:final-chier-tout-le-monde}
( \tau^T \sigma_{s}n, w\cdot \tau)_{\Gamma(t)} +
(  (n^T\jump{\widetilde \sigma n}) , (w\cdot \tau_w) ( n\cdot \tau_w))_{\Gamma(t)} +
\gamma_C (d\cdot n_w - g_0 , w\cdot n_w)_{\Gamma(t)}=0.
\end{align}
When it comes to contact, we have in the asymptotic limit $\tau_w=\tau$ and $ n\cdot\tau_w =0$, 
so that \eqref{eq:final-chier-tout-le-monde} reduces to the ``no-friction'' condition 
$\tau^T \sigma_s n =0$ and the non-penetration condition $d\cdot n_w =
g_0$. This is the desired contact condition, as in the no-slip case, see \eqref{dGflux_argumentation}.

{To summarise we have shown that the following conditions are weakly imposed on $\Gamma(t)$ in the asymptotic limit
\begin{align*}
 \tau^T \sigma_f n = 0, \quad (\dot{d} -u) \cdot n = 0, \quad \tau^T \sigma_s n = 0.
\end{align*}
Moreover, we have on the part $\Gamma_C(t) \subset \Gamma(t)$ that is in contact with $\Gamma_w$
\begin{align*}
 d\cdot n_w - g_0 = 0,
\end{align*}
and on the part $\Gamma_{\text{fsi}}(t) \subset \Gamma(t)$ that is not in contact
\begin{align*}
 \sigma_{f,n} - \sigma_{s,n} = 0.
\end{align*}
}

We observe that in \eqref{eq:final-chier-tout-le-monde} both the
tangential and the normal components (with respect to $\tau_w$ and
$n_w$) have asymptotically vanishing perturbations. It is possible to eliminate the
perturbation in the normal component, from the solid stress term $( \tau^T \sigma_{s}n, w\cdot
\tau)_{\Gamma(t)}$, by adding the corresponding
term to the definition of $\lambda$, {i.e.}
\begin{equation}\label{eq:slip-perturbed}
\lambda = %(n^T\jump{\widetilde \sigma n}) (n \cdot n_{\rm w})
{\jump{\widetilde \sigma_{n,\text{slip}}}}+ (\tau^T \sigma_s n) ( \tau \cdot n_w).
\end{equation}
In this case the relation \eqref{eq:final-chier-tout-le-monde} takes
the form
\begin{align*}
( \tau^T \sigma_{s}n, w\cdot \tau_w (\tau \cdot \tau_w))_{\Gamma(t)} +
(  (n^T\jump{\widetilde \sigma n}) , {(w\cdot \tau_w) (n\cdot \tau_w)})_{\Gamma(t)} +
\gamma_C (d\cdot n_w - g_0 , w\cdot n_w)_{\Gamma(t)}=0.
\end{align*}
Testing with $w=n_w$ we see that here the non-penetration condition $d\cdot n_w =
g_0$ is imposed unperturbed for all $h>0$. By moving the perturbation
to $\lambda$ {as in~\eqref{eq:slip-perturbed}} it is instead the form $P_\gamma$ that is perturbed,
which implies a (weakly consistent) perturbation of the contact zone.
}

\section{Stability}
\label{sec.discstab}

 In this work, we will use fitted finite elements, i.e. we assume that both 
the interface $\Gamma(t)$ and the wall $\Gamma_w$ are resolved by mesh lines. The 
approaches presented, however, can be combined with unfitted finite elements as well, with the only difference that some more stabilisation 
terms have to be added to the variational formulation.
In order to simplify the presentation, we will concentrate on the fitted case 
first and discuss the extension to
unfitted finite elements in a remark afterwards.

We will use equal-order finite elements in combination with a pressure stabilisation 
term $S_p(p,q)$ for the fluid equations. 
For the stability analysis, 
the only requirement on $S_p$ is that it leads to a well-posed discrete fluid problem. Possibilities
include the Brezzi-Pitk\"aranta stabilisation~\cite{BrezziPitkaeranta1984}, local projections 
(LPS)~\cite{BeckerBraack2001}, the pressure-stabilised Petrov-Galerkin approach (PSPG)~\cite{HughesFrancaBalestra1986}
or the Continuous Interior Penalty method (CIP)~\cite{BurmanHansbo2006}).

{In order to present the stability analysis in a general setting, we introduce some 
further notation splitting the {contact force} variable $\lambda$ into 
a fluid part $\lambda_f$ and a solid part $\lambda_s$, {such that $\lambda = \lambda_s - \lambda_f$}. For the no-slip case, we define 
$\lambda_s(w)=\sigma_{s,n}(w)$ and
\begin{align}
 \lambda_f(v,q,w) &= \begin{cases} \sigma_{f,n}(v,q), \qquad\quad\qquad \text{if } \lambda = \jump{\sigma_n},\\
 T_f(v,q,w)\cdot n, \quad\qquad\, \text{if } \lambda = \jump{\tilde{\sigma}_n},
 \end{cases}
\end{align}
see \eqref{DefLagr} and \eqref{LagrFluxes}. For the slip case, we define $\lambda_s(w)=\sigma_{s,n}(w)(n\cdot n_w)$ and
\begin{align}
 \lambda_f(v,q,w) &= \begin{cases} \sigma_{f,n}(v,q)(n\cdot n_w), \qquad\quad\;\; \text{if } \lambda = \jump{\sigma_{n,\text{slip}}},\\
 (T_f(v,q,w) \cdot n) (n\cdot n_w), \quad \text{if } \lambda = \jump{\widetilde \sigma_{n,\text{slip}}},
 \end{cases}
\end{align}
see \eqref{LagrCompleteSlip} and \eqref{lambdaslip}.
}

\subsection{Generalised contact formulation}
\label{sec.symmetric}

 Before we conduct the stability analysis, let us introduce a generalised
contact formulation,
following the ideas of Chouly \emph{et al.}~\cite{ChoulySymVsNonsym}. We have already briefly 
discussed their ideas for the case of 
a pure solid with contact in \eqref{thetasolid}.
The generalisation of (\ref{thetasolid}) to the FSI-contact system
(Variational Formulation \ref{varForm_artF})
would be to add the terms
\begin{align}\label{impractical}
{-} \left(\gamma_C [P_{\gamma}(\lambda,d)]_+ + \lambda , 
\lambda(\partial_t v, \partial_t q, w)\right)_{\Gamma(t)}.
\end{align}
The time derivatives on the test functions $v$ and $q$ are motivated by the stability analysis
below, where we have to test the variational form
with $v=u$, $q=p$ and $w=\partial_t d$, in order to show stability (see also Burman \& Fern\'andez~\cite{BurmanFernandez2014}).

On the other hand, the term (\ref{impractical}) is not usable within a time-stepping 
scheme due to the time derivatives on the test functions. A remedy is to shift the time derivatives
to the first integrand (ignoring the boundary terms), i.e.$\,$adding the consistent terms
\begin{align*}
 &{-}\left(\gamma_C [P_{\gamma}(\lambda,d)]_+ + \lambda , 
{\lambda_s(w)}\right)_{\Gamma(t)} {-} \left(\partial_t \left(\gamma_C [P_{\gamma}(\lambda,d)]_+ 
+ \lambda \right), {\lambda_f(v,q,w)}\right)_{\Gamma(t)}.
\end{align*}

\noindent This yields the variational formulation: \begin{varform}
\textit{Find $u(t)\in {\cal V}_h, p(t)\in {\cal Q}_h, d(t)\in {\cal W}_h$, such that $\dot{d}=\partial_t d$ and }
\begin{eqnarray}\label{DiscSysGen}
 \begin{aligned}
 {\cal A} (u,p,d)(v,q,w)
 &:={\cal A}_{*,FSI}^*(u,p,d)(v,q,w)+\gamma_{C} \left([P_\gamma(\lambda,d)]_+, w\cdot n_w\right)_{\Gamma(t)}
 \\
 &\qquad
 {-}\theta \left(\gamma_C [P_{\gamma}(\lambda,d)]_+ + \lambda , 
{\lambda_s(w)}\right)_{\Gamma(t)}
-\theta \left(\partial_t \left(\gamma_C [P_{\gamma}(\lambda,d)]_+ + \lambda \right), 
{\lambda_f(v,q,w)}\right)_{\Gamma(t)} \\
 &\qquad\qquad\qquad\qquad\qquad= \left(f_f,v\right)_{\Omega_f(t)} + \left(f_s,w\right)_{\Omega_s(t)}\quad 
 \forall v,q,w \in {\cal V}_h \times {\cal Q}_h \times {\cal W}_h,
\end{aligned}
\end{eqnarray}
where ${\cal A}_{*,FSI}^*$ is one of the bilinear forms ${\cal A}_{\text{no-slip,FSI}},{\cal A}_{\text{slip,FSI}}, 
{\cal A}_{\text{no-slip,FSI}}^C$ or ${\cal A}_{\text{slip,FSI}}^C$.
\end{varform}

\subsection{A stability result}
\label{sec.stability}

In this section, we will investigate stability of the discrete formulation for different values of $\theta$.
In particular, we will show a stability result for the symmetric formulation ($\theta=1$) and stability up to a specific
term in the general case (including $\theta=0$). These results correspond to the results that 
have been obtained by Chouly \emph{et al.}~\cite{ChoulySymVsNonsym} for the pure solid case. For the stability 
analysis, we will assume infinitesimal displacements, i.e. the sub-domains $\Omega_f$ and $\Omega_s$ as 
well as the interface $\Gamma$ are fixed. We introduce the notation $\widetilde{\Omega}_f$ for the combined fluid 
and artificial fluid domain in the case of ${\cal A}_{\text{*,FSI}}^C$ and for $\Omega_f$ else.

\begin{theorem}\label{theo.stability}
Let $u,p,d \in V_h$ and $\dot{d}=\partial_t d$. We have the 
following stability result for the form ${\cal A}$ defined in \eqref{DiscSysGen},
where $\theta \in [0,1]$ and $\gamma_C^0$ sufficiently large
\begin{align*}
 &\| u(T) \|^2_{\widetilde{\Omega}_f} + \| \dot{d}(T)\|^2_{\Omega_s} +\| d(T) \|_{H^1(\Omega_s)}^2
 +\int_0^T \left(\nu_f \|\nabla u\|_{\widetilde{\Omega}_f}^2 + S_p(p,p)
 +\gamma_a \|u\|_{\Omega_f^C}^2 \right) dt \\
 &\qquad\qquad+\int_0^T \gamma_{\text{fsi}} \| (\dot{d}-u)\cdot n\|^2_{\Gamma}\,dt 
 + \theta \left\| \gamma_C^{-1/2} \lambda (T) 
 +\gamma_C^{1/2} [P_\gamma(\lambda,d)]_+ (T) \right\|_\Gamma^2\\
 &\qquad\qquad\qquad\qquad\leq C\Bigg( \int_0^T {\cal A}(u,p,d; u,p,\dot{d})\, 
 - (1-\theta)\gamma_C \left([P_\gamma(\lambda,d)]_+, \dot{d}\cdot n_w\right)_{\Gamma} dt 
 + \| u_0 \|^2_{\widetilde{\Omega}_f}\\
 &\qquad\qquad\qquad\qquad\qquad\qquad + \| \dot{d}_0 \|^2_{\Omega_s} 
 + \| d_0 \|_{H^1(\Omega_s)}^2 +\theta 
\left\| \gamma_C^{-1/2}\lambda_0
 + \gamma_C^{1/2} [P_{\gamma}(\lambda_0,d_0)]_+ \right\|_\Gamma^2\Bigg),
\end{align*}
where we have used the abbreviation $\lambda_0:=\lambda(u_0,p_0,d_0)$.
For the no-slip case, the term $\| (\dot{d}-u)\cdot n\|_{\Gamma}$ can be replaced by $\| (\dot{d}-u)\|_{\Gamma}$.

\end{theorem}
\begin{remark}{(Contact terms)}
 The second line gives us control over the satisfaction of the FSI-contact condition for $\theta>0$.
{In contrast to the work by Chouly  \emph{et al.}\cite{ChoulySymVsNonsym} for a pure solid problem, here we obtain 
discrete stability for $\theta=1$, for the following positive discrete energy}
\[
{E(T):=\| u(T) \|^2_{\widetilde{\Omega}_f} + \| \dot{d}(T)\|^2_{\Omega_s} +\| d(T) \|_{H^1(\Omega_s)}^2+\left\| \gamma_C^{-1/2} \lambda (T) 
 +\gamma_C^{1/2} [P_\gamma(\lambda,d)]_+ (T) \right\|_\Gamma^2.}
\]
For $\theta\neq 1$ on the other hand, the contact term 
 $(1-\theta)\gamma_C \left([P_\gamma(\lambda,d)]_+, \dot{d}\cdot n_w\right)_{\Gamma}$
 appears on the right-hand side. The last term
 on the right-hand side vanishes, if we assume that the contact conditions are fulfilled at initial time,
 for example if the solid is not in contact with at $t=0$.
\end{remark}

\begin{proof}
We test (\ref{DiscSysGen}) with $w=\dot{d}=\partial_t d, v=u$ and $q=p$ and integrate in time.
We start by deriving a lower bound for ${\cal A}_{*,\text{FSI}}^*$.
For the fluid part, we use the techniques from Burman \& Fernandez~\cite{BurmanFernandez2014}, 
to show coercivity of the Stokes part including the coupling terms. For the no-slip case, the authors have shown
\begin{align*}
  &\left(\sigma_f(u,p),\nabla u\right)_{\widetilde{\Omega}_f} 
 + \left(\text{div } u, p\right)_{\widetilde{\Omega}_f} + S_p(p,p) + \gamma_a(u,u)_{\Omega_f^C}\\
 &\qquad\quad- (\sigma_f(u,p) n, \dot{d}-u)_{\Gamma} - (\dot{d}-u,\sigma_f(u,-p) n)_{\Gamma} 
 + \gamma_{\text{fsi}} \| \dot{d}-u\|^2_{\Gamma}\\
 &\qquad\qquad\qquad\geq c \left(\nu_f \|\nabla u\|_{\widetilde{\Omega}_f}^2 
 + \gamma_{\text{fsi}} \| \dot{d}-u\|^2_{\Gamma} + S_p(p,p) + \gamma_a\|u\|_{\Omega_f^C}^2\right).
\end{align*}
Analogously, one can show in the slip-case that
\begin{align*}
  &\left(\sigma_f(u,p),\nabla u\right)_{\widetilde{\Omega}_f} 
 + \left(\text{div } u, p\right)_{\widetilde{\Omega}_f} + S_p(p,p) + \gamma_a(u,u)_{\Omega_f^C}\\
 &\qquad\quad- (n^T \sigma_f(u,p) n, (\dot{d}-u)\cdot n)_{\Gamma} - ((\dot{d}-u)\cdot n,n^T\sigma_f(u,-p) n)_{\Gamma} 
 + \gamma_{\text{fsi}} \| (\dot{d}-u)\cdot n\|^2_{\Gamma}\\
 &\qquad\qquad\qquad\geq c \left(\nu_f \|\nabla u\|_{\widetilde{\Omega}_f}^2 
 + \gamma_{\text{fsi}} \| (\dot{d}-u)\cdot n\|^2_{\Gamma} + S_p(p,p) + \gamma_a\|u\|_{\Omega_f^C}^2\right).
\end{align*}

\noindent Using the symmetry of $\sigma_s$, integration in time and a Korn's inequality, we obtain for 
the solid part
 \begin{align*}
 \int_0^T (\sigma_s(d), \nabla \dot{d})_{\Omega_s} \, dt
 = \frac{1}{2}\int_0^T \partial_t (\sigma_s(d), \nabla d)_{\Omega_s} \, dt
 &= \frac{1}{2} \left( (\sigma_s(d(T)), \nabla d(T))_{\Omega_s} - 
 (\sigma_s(d(0)), \nabla d(0))_{\Omega_s} \right)\\
 &\geq c_1 \|\nabla d(T)\|_{\Omega_s}^2 -c_2 \|\nabla d_0\|_{\Omega_s}^2.
 \end{align*}

\noindent Moreover, we have 
\begin{align*}
 \int_0^T (\partial_t u, u)_{\widetilde{\Omega}_f} + (\partial_t \dot{d}, \dot{d})_{\Omega_s} \, dt
 = \frac{1}{2} \left(\| u(T)\|_{\widetilde{\Omega}_f}^2 + \| \dot{d}(T)\|_{\Omega_s}^2 
 - \| u_0\|_{\widetilde{\Omega}_f}^2 - \|\dot{d}_0\|_{\Omega_s}^2\right).
\end{align*}
Together, we have shown that 
\begin{align*}
  \| u(T) \|^2_{\widetilde{\Omega}_f} + \| \dot{d}(T)\|^2_{\Omega_s} +\| d(T) \|_{H^1(\Omega_s)}^2
 &+\int_0^T \nu_f \|\nabla u\|_{\widetilde{\Omega}_f}^2 + S_p(p,p)
 +\gamma_a \|u\|_{\Omega_f^C}^2 + \gamma_{\text{fsi}} \| \dot{d}-u\|^2_{\Gamma}\, dt\\
 &\leq C\Bigg( \int_0^T {\cal A}_{*,FSI}^*(u,p,d; u,p,\dot{d})\,dt
 + \| u_0 \|^2_{\widetilde{\Omega}_f} + \| \dot{d}_0 \|^2_{\Omega_s} + \| d_0 \|_{H^1(\Omega_s)}^2\Bigg)
\end{align*}
 
\noindent Let us now estimate the contact terms. We split the principal contact term into
\begin{align*}
\int_0^T\gamma_C \left([P_\gamma(\lambda,d)]_+, \dot{d}\cdot n_w\right)_{\Gamma} \,dt
&=\int_0^T\theta\gamma_C \left([P_\gamma(\lambda,d)]_+, \dot{d}\cdot n_w\right)_{\Gamma}
 +(1-\theta)\gamma_C \left([P_\gamma(\lambda,d)]_+, \dot{d}\cdot n_w\right)_{\Gamma}\, dt.
\end{align*}
We have to estimate the terms
\begin{eqnarray}
\begin{aligned}\label{contTerms}
\theta \int_0^T \underbrace{\gamma_C \left([P_\gamma(\lambda,d)]_+, \dot{d}\cdot n_w\right)_{\Gamma}}_{I_1}
&-\underbrace{\left(\gamma_C [P_\gamma(\lambda,d)]_+ + \lambda , 
{\lambda_s(\dot{d})}\right)_{\Gamma}}_{I_2}\\
&-\underbrace{\left(\partial_t \left(\gamma_C [P_\gamma(\lambda,d)]_+ 
+ \lambda \right), 
{\lambda_f(u,p,\dot{d})}\right)_{\Gamma}}_{I_3}\, dt
\end{aligned}
\end{eqnarray}

From the definition of $P_{\gamma}$ {we can write $ d\cdot n_w = P_\gamma(\lambda,d) + \gamma_C^{-1} \lambda$. Hence, since 
 the lower wall is assumed to be time independent, we have} 
\begin{align*}
 \int_0^T I_1 \, dt= &\gamma_C\int_0^T\big([P_\gamma(\lambda,d)]_
 + ,\partial_t (P_\gamma(\lambda,d) + \gamma_C^{-1} \lambda )\big)_\Gamma\, dt \\
  = &\gamma_C\int_0^T\big([P_\gamma(\lambda,d)]_+ ,\partial_t ([P_\gamma(\lambda,d)]_+ 
  + \gamma_C^{-1} \lambda )\big)_\Gamma\, dt. \\
 \end{align*}  
 In the second line, we have used that~\cite{ChoulySymVsNonsym}
\begin{align*}
 \frac{1}{2} \partial_t [\phi]_+^2 =  [\phi]_+ \partial_t [\phi]_+ = [\phi]_+ H(\phi) \partial_t [\phi]_+ = [\phi]_+ \partial_t \phi 
\end{align*}
where $H$ is the Heaviside function. We insert $\pm \gamma_C^{-1} \lambda $ 
and integrate by parts  
 \begin{align*}  
 \int_0^T I_1 \,dt
 &=\gamma_C\int_0^T\big([P_\gamma(\lambda,d)]_+
 +\gamma_C^{-1} \lambda ,
 \partial_t ([P_\gamma(\lambda,d)]_+ + \gamma_C^{-1} \lambda )\big)_\Gamma\,dt \\
  &\qquad\qquad- \int_0^T\underbrace{\big(\lambda ,\partial_t ([P_\gamma(\lambda,d)]_+ 
  + \gamma_C^{-1} \lambda )\big)_\Gamma}_{I_4} \,dt \\
  %= & \frac{\gamma_C}{2}\int_0^T\frac{\rm d}{{\rm d}t}\| [P_\gamma(\lambda,d)]_+ + \gamma_C^{-1} \lambda(u,p,d)  \|_{0,\Gamma}^2 \\ 
    %&- \int_0^T\underbrace{\big(\lambda(u,p,d),\partial_t ([P_\gamma(\lambda,d)]_+  + \gamma_C^{-1} \lambda(u,p,d) )\big)_\Gamma}_{I_4}\\
    &= \underbrace{\frac{\gamma_C}{2}\|  [P_\gamma(\lambda,d)(T)]_+  + \gamma_C^{-1} \lambda (T) \|_{0,\Gamma}^2}_{D_T}
   - \underbrace{\frac{\gamma_C}{2}\|  [P_{\gamma}(\lambda_0,d_0)]_+  + \gamma_C^{-1} \lambda_0 \|_{0,\Gamma}^2}_{D_0}   
       - \int_0^T  I_4 \,dt.
\end{align*}
Now the idea is to split the contribution from the term $I_4$ into fluid and solid stresses and to apply integrating by parts in time (only)
in the solid stress contribution. By definition, 
we have  
$$
\begin{aligned}
 -\int_0^T I_4\,dt  
 =&  - \int_0^T \big(\lambda_s(d) , \partial_t([P_\gamma(\lambda,d)]_+  
 + \gamma_C^{-1} \lambda  )\big)_\Gamma\,+\, \big(\lambda_f(u,p,\dot{d}) , \partial_t([P_\gamma(\lambda,d)]_+  
  + \gamma_C^{-1} \lambda  )\big)_\Gamma\,dt
\end{aligned}
$$
and, by integrating by parts in the first term of the right-hand side, we have 
$$
\begin{aligned}
 -\int_0^T I_4\,dt  = & \underbrace{\int_0^T \big(\lambda_s(\dot{d}),
 [P_\gamma(\lambda,d)]_+  + \gamma_C^{-1} \lambda \big)_\Gamma }_{I_5}\,dt
  \underbrace{- \big(\lambda_s( d(T)) , [P_\gamma(\lambda,d)(T)]_+  + \gamma_C^{-1} \lambda (T) \big)_\Gamma}_{I_6} \\
  &+ \underbrace{\big(\lambda_s( d_0) , [P_{\gamma}(\lambda_0,d_0)]_+  
  + \gamma_C^{-1} \lambda_0\big)_\Gamma}_{I_7}
  +  \underbrace{\int_0^T \big(\lambda_f(u,p,\dot{d}),
  \partial_t([P_\gamma(\lambda,d)]_+  + \gamma_C^{-1} \lambda  )\big)_\Gamma}_{I_8}\,dt
\end{aligned}
$$
The terms $I_5$ and $I_8$ cancel with the terms $I_2$ and $I_3$ in \eqref{contTerms}. 
The term $I_6$ is treated in a standard fashion using Young's inequality,
an inverse inequality and the dissipation provided by $D_T$ and 
the elastic energy $\frac12a^{\rm s}(d(T),d(T))$ for $\gamma_C^0$ sufficiently large:
\begin{align*}
I_6 &\geq - \gamma_C^{-1}\| \lambda_s(d(T)) \|_{0,\Gamma}^2 
-\frac{\gamma_C}{4}  \| [P_{\gamma}(\lambda(T),d(T))]_+  
+ \gamma_C^{-1} \lambda (T) \|_{0,\Gamma}^2\\ 
&\geq  - \frac{1}{4} \|d(T)\|_{H^1(\Omega_s)}^2 - \frac{1}{2} D_T.
\end{align*}
In the same way, we obtain for $I_7$
\begin{align*}
 I_7 \leq \frac{1}{4} \|d_0\|_{H^1(\Omega_s)}^2 + \frac{1}{2} D_0,
\end{align*}
which completes the proof.
\end{proof}

 \begin{remark}{(Unfitted finite elements)}
When using unfitted finite elements~\cite{BurmanHansbo2012, BurmanFernandez2014},
additional stabilisation terms $S_u$ and $S_d$ are needed, if
the interface $\Gamma(t)$ is not resolved by mesh lines. Their purpose is to extend 
the coercivity of the fluid
system from $\Omega_f$ (resp. $\Omega_s$) to the extended domains $\Omega_h^f$ (resp. $\Omega_h^s$)
that consists of all element $T\in {\cal T}_h$,
with a non-empty intersection with the respective sub-domain ($T\cap \Omega_i\neq \emptyset$). 
Suitable ``ghost penalty'' operators have been defined in Burman~\cite{Burman2010}. The same stability result
as in Theorem~\ref{theo.stability} can then be shown with an analogous argumentation.
\end{remark}

\begin{remark}{(Newton convergence)}\label{rem.generalContact}
While the symmetric formulation ($\theta=1$) seems beneficial from the theoretical point of view, 
the additional terms in (\ref{DiscSysGen}) can cause severe difficulties for the non-linear solver. 
The reason is that the
additional contact terms are not only highly non-linear, but also
non-smooth, especially due to the time derivative acting on the maximum operator $\partial_t [\cdot]_+$.
In our numerical
tests, we were not able to obtain numerical convergence for different versions of generalised Newton 
methods~\cite{TureketalNewton}. The 
investigation of the case $\theta\neq 1$ by means of numerical tests and in particular the 
construction of a robust non-linear solver
are subject to future research.
\end{remark}

\section{Numerical Results}
\label{sec.num}

In this section, we show some numerical results to analyse and to compare the different contact formulations.
As mentioned in Remark~\ref{rem.generalContact}, we 
were not able to obtain results for the generalised contact formulation with $\theta \neq 0$,
due to divergence of the generalised Newton-type methods we have tried. Therefore, 
we only show results for $\theta=0$, where Newton convergence was not an issue, at least when the 
time step $\delta t$ was chosen reasonably small. 
{Following the standard approach for contact in solid mechanics, we could in this case simply ignore the 
non-differentiability of the maximum operator when computing the Newton
derivatives, as the term $P_\gamma$ inside the bracket $[\cdot]_+$ is in practice typically never exactly zero.
For all other values of $P_\gamma$ the derivatives are well-defined.} 
In the computations made for this paper, 
the Newton algorithm needed 1-2 iterations per time step to 
reduce the initial residual by a factor of $10^{-7}$, if the contact force was not getting active during the iteration,
and 1-5 iterations per time step in and around the interval of contact.
{This makes the method highly competitive in terms of computational costs compared to approaches using 
Lagrange multipliers and/or active-sets.}

We first give some details in Section~\ref{sec.numdisc} on the fitted, equal-order finite element discretisation 
and the stabilisations we use. 
Then, in Section~\ref{sec.obstacle}, we study {the problem of a virtual obstacle within the fluid domain introduced
in Section~\ref{sec.ObstModel}.}
The purpose of this example is to isolate the effect of the contact
terms from issues related to discretisation during contact and the 
topology change in the fluid domain $\Omega_f(t)$.
Then, we study in Section~\ref{sec.contact} a model
problem with {contact with the boundary of the fluid domain,
where we compare among other aspects the two contact formulations introduced in Section~\ref{sec.FSIcontact}, 
the different possibilities to choose the fluxes $\lambda$ and the effect of slip and no-slip boundary 
and interface conditions.} 

\subsection{Details on discretisation and stabilisation}
\label{sec.numdisc}

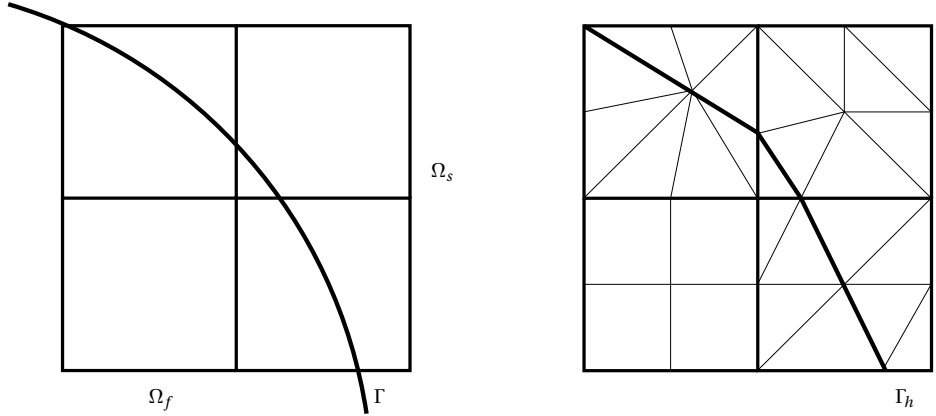
\begin{figure}[bt]
  \centering
  \begin{picture}(0,0)%
\includegraphics{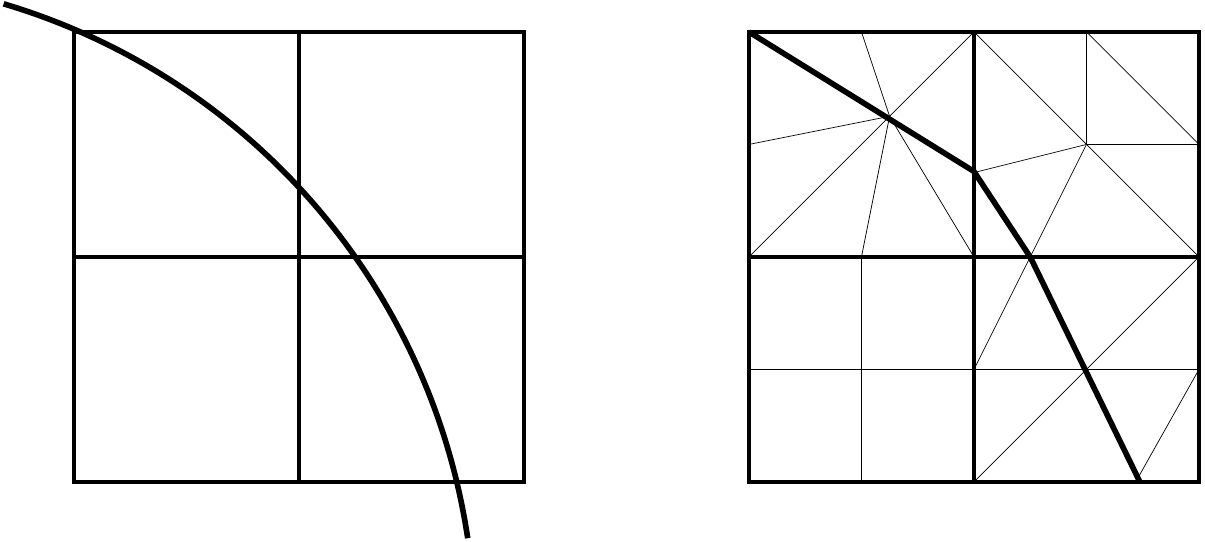}%
\end{picture}%
\setlength{\unitlength}{1776sp}%
\begingroup\makeatletter\ifx\SetFigFont\undefined%
\gdef\SetFigFont#1#2{%
  \fontsize{#1}{#2pt}%
  \selectfont}%
\fi\endgroup%
\begin{picture}(12836,5776)(413,-5799)
\put(12710,-5613){\makebox(0,0)[lb]{\smash{{\SetFigFont{8}{9.6}{\color[rgb]{0,0,0}$\Gamma_h$}%
}}}}
\put(5506,-5611){\makebox(0,0)[lb]{\smash{{\SetFigFont{8}{9.6}{\color[rgb]{0,0,0}$\Gamma$}%
}}}}
\put(2401,-5611){\makebox(0,0)[lb]{\smash{{\SetFigFont{8}{9.6}{\color[rgb]{0,0,0}$\Omega_f$}%
}}}}
\put(6301,-2461){\makebox(0,0)[lb]{\smash{{\SetFigFont{8}{9.6}{\color[rgb]{0,0,0}$\Omega_s$}%
}}}}
\end{picture}%  
 
  \caption{\textit{Left:} Fixed triangulation ${\cal T}_{2h}$ of the domain $\Omega$. \textit{Right:}
    Subdivision of the patches $P\in{\cal T}_{2h}$ such that the interface $\Gamma(t)$
    is resolved in a linear approximation by the discrete interface $\Gamma_h$.\label{fig.lmfem}} 
  \label{fig:mesh}
\end{figure}

{For the numerical results in this paper, we will use a monolithic \textit{Fully Eulerian} approach
on a global mesh ${\cal T}_h$ covering $\Omega(t)$. In order to resolve the interface $\Gamma(t)$ within the discretisation,} we
use the \textit{locally modified finite element method}
introduced by Frei \& Richter~\cite{FreiRichter2014}. The idea of this approach is to use a fixed
coarse triangulation ${\cal T}_{2h}$ of the overall domain $\Omega=\widetilde{\Omega}_f(t)
\cup\Gamma(t)\cup\Omega_s(t)$ that is independent of the position of the interface $\Gamma(t)$. Then,
in each time step, this coarse grid is refined once by splitting each so-called ``patch'' element
in either
eight triangles or four quadrilaterals to resolve the interface in at least a linear approximation, see
Figure~\ref{fig.lmfem} for an illustration.

The finite element space $V_h$ is then defined as a combination of piece-wise linear and 
piece-wise bi-linear finite elements on the patches. It can be
guaranteed that a maximum angle condition is fulfilled in each of the sub-cells, leading to optimal-order 
interpolation and error estimates~\cite{FreiRichter2014}.

For temporal discretisation, we split the time interval $I$ into $m$ equidistant-distant time 
intervals $I_j=(t_{j-1},t_{j-1}+k]$ and
use a time-stepping scheme that is based on a modified discontinuous 
Galerkin time discretisation of lowest order (dG(0)), see Frei \& Richter~\cite{FreiRichter2017}. 
{The {displacement-velocity relation} $\partial_t d = \dot{d}$ is included by means of the $L^2$-projection
\begin{align*}
 (\partial_t d, z)_{\Omega_s(t)} - (\dot{d}, z)_{\Omega_s(t)} &= 0 \quad \forall z\in V_h^s,
\end{align*}}
where $V_h^s$ denotes the (modified) finite element space that is spanned by the degrees of freedom 
of the elements in the solid part $\Omega_h^s$.

{The domain affiliation of a point $x\in\Omega(t)$ is determined by means of the 
\textit{Initial Point Set/Backward Characteristics method}~\cite{DunneRannacher, Cottetetal},
that uses the displacement $d(t)$ in the solid domain and an extension to $\Omega_f(t)$ in order to trace back points 
to their initial position in $\Omega(0)$, following the definition~\eqref{MovingDomains}.}

For pressure stabilisation, we use an anisotropic variant of the Continuous Interior Penalty method,
see Frei~\cite{FreiToAppear, FreiPhD}. In addition, we add the temporal pressure stabilisation term
\begin{align*}
 S_{pt}(p,q) = \gamma_{pt} h (p^m-p^{m-1},q)_{\Gamma(t)}
\end{align*}
in each time interval $I_m$. This additional stabilisation is needed, as the mesh ${\cal T}_h(t_m)$, and
hence the finite element spaces, change from time-step to time-step. The solution $u^{m-1}$ 
{from the previous time-step $t_{m-1}$} is therefore not discrete divergence-free with respect to the new 
mesh ${\cal T}_h(t_m)$, which
gives rise to pressure oscillations, see for example Besier \& Wollner~\cite{BesierWollner}.

All the following results have been obtained using the finite element library \textit{Gascoigne 3d}~\cite{Gascoigne}.

\subsection{{Virtual} obstacle within the fluid domain}
\label{sec.obstacle}

In order to isolate the effect of the contact model from issues related to the Navier-Stokes-contact paradox and the topology changes 
in the fluid domain,
let us first study the simplified problem introduced in Section~\ref{sec.ObstModel}. We define the initial fluid and solid domains as
\begin{align*}
\Omega_f(0)= (0,1) \times (0,0.5), \qquad \Omega_s(0)= (0,1) \times (0.5,0.6) 
\end{align*}
and a lower-dimensional obstacle $\Gamma_w = (0,1) \times {0.25}$ within the fluid domain.
We consider a moving interface $\Gamma(t)$, 
which is resolved using the locally modified finite element method. The sub-domains $\Omega_s(t)$ 
and $\Omega_f(t)$ and the interface $\Gamma(t)$ depend on the solid displacement $d(t)$, 
see (\ref{MovingDomains}).

The constraint for the solid displacement is given by
\begin{align}
 d\cdot n_w \leq 0.25 \quad (=:g_{\alpha}).\label{ContObstacle}
\end{align}

We use the elasticity parameters $\lambda_s=\mu_s=2\cdot 10^6$ and the fluid 
viscosity $\nu_f=1$. The structure is pulled towards the bottom
by fluid forces due to a prescribed pressure mean value at the left and right boundary of the fluid domain 
\begin{align*}
  \int_{\Gamma_{f,\text{left}}} p \, ds 
  = \int_{\Gamma_{f,\text{left}}} \overline{P}\, ds,\qquad 
  \int_{\Gamma_{f,\text{right}}} p \, ds = \int_{\Gamma_{f,\text{right}}} \overline{P}\, ds
\end{align*}
where $\overline{P}:=1.3\cdot 10^5$.
We consider the {Variational Formulation~\ref{varForm_implicit} with no-slip conditions
and $\lambda=\jump{\tilde{\sigma}_n}$, see \eqref{LagrFluxes},
on a Cartesian mesh that
consists of 5120 elements and with a small time step $\delta t = 10^{-5}$.}
The Nitsche constant at the 
FSI interface is chosen $\gamma_{\text{fsi}}^0 = 10^3$
and temporal pressure stabilisation with $\gamma_{pt}=10^{-2}$ is used. 

To analyse the results, we define the ``minimal distance''
\begin{align*}
 d_{\min} := \min_{x\in \Gamma(t)} x_2 - 0.25
\end{align*}
 of the interface $\Gamma(t)$ to $\Gamma_w$. To be precise the term ``minimal distance'' is only
 correct before contact, as $d_{\min}$ gets negative in case of an overlap. Moreover, we define 
 the following functionals in order to analyse the pressure $p$ and 
 the contact force at the interface $\Gamma(t)$
 \begin{align*}
  J_p &:= \big| \int_{\Gamma(t)} p \, ds \big|,\qquad
  J_{P_{\gamma}} := \gamma_C \int_{\Gamma(t)} [P_\gamma(\lambda,d)]_+ \, ds.
 \end{align*}

\begin{figure}[h!]
\centering
 \begin{minipage}{0.51\textwidth}
  \includegraphics[width=\textwidth]{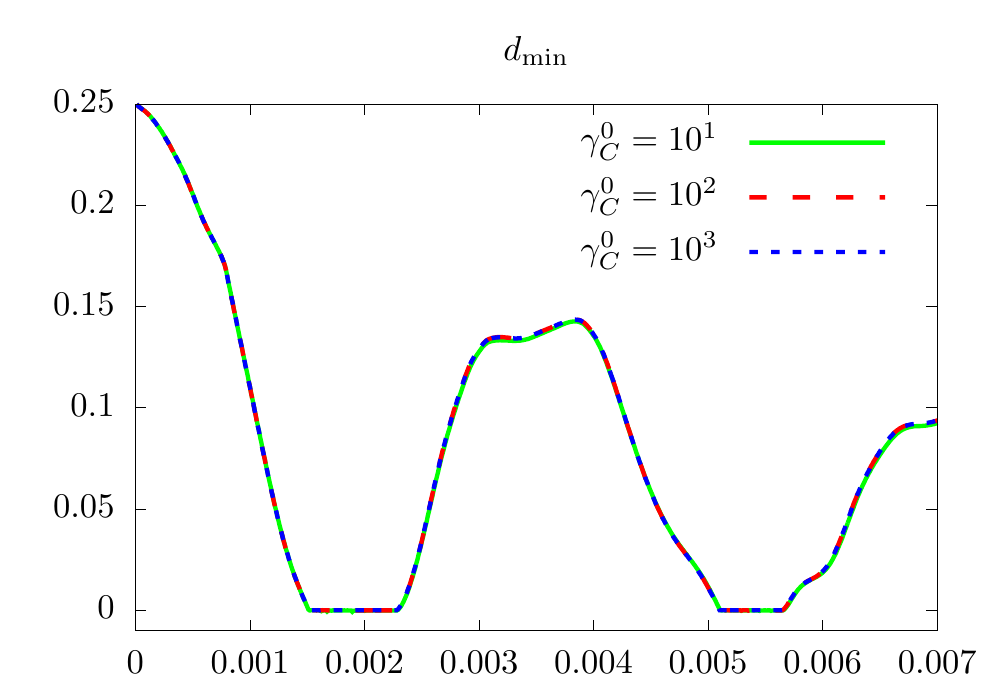}
 \end{minipage}
 \hspace{-0.6cm}
  \begin{minipage}{0.51\textwidth}
    \includegraphics[width=\textwidth]{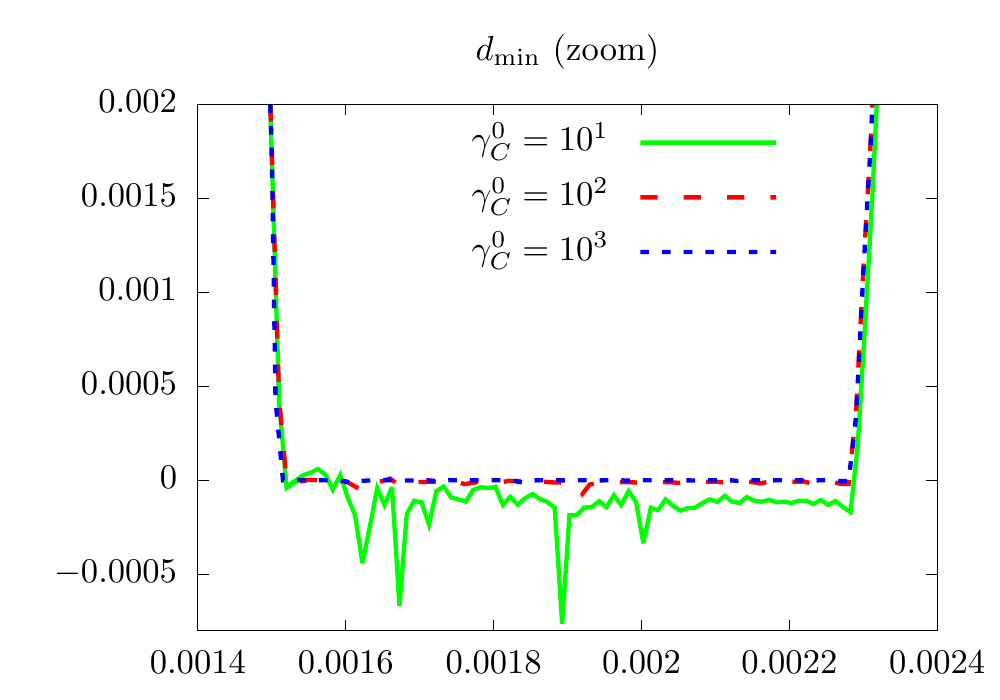}
 \end{minipage}
 
 \begin{minipage}{0.51\textwidth}
  \includegraphics[width=\textwidth]{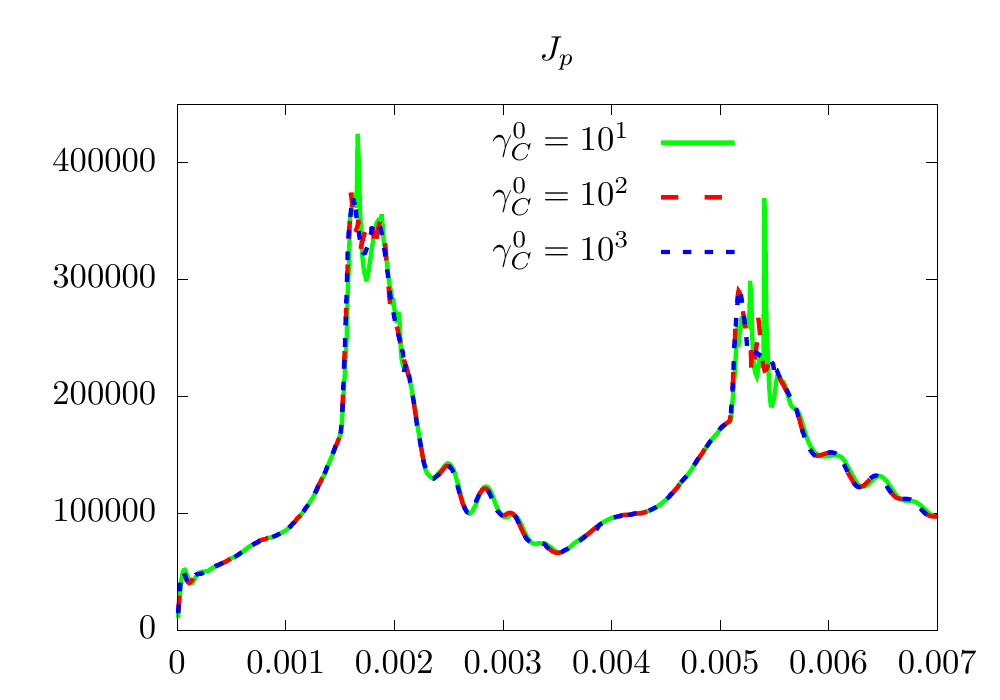}
 \end{minipage}
 \hspace{-0.6cm}
    \begin{minipage}{0.51\textwidth}
     \includegraphics[width=\textwidth]{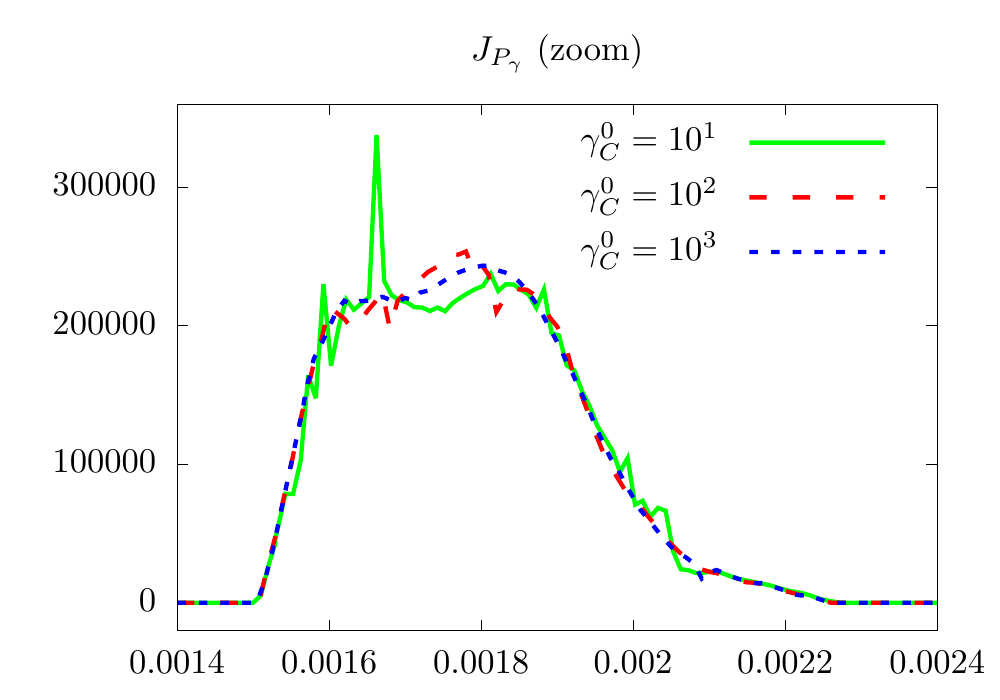}
  \end{minipage}
 \caption{\label{fig.vertDefandpress} \textit{Top row}: 'Minimal distance' $d_{\min}$ to $\Gamma_w$
 over time with two contact periods. \textit{Left}: Total time interval. 
 \textit{Right}: Zoom-in at the first contact interval. \textit{Bottom left}: 
 Pressure functional $J_p$ over time. 
 \textit{Bottom right}:
 Functional $J_{P_{\gamma}}$ measuring the contact force around the contact interval 
 over time.}
\end{figure}

In Figure~\ref{fig.vertDefandpress}, we plot these three functionals  
over time for two contact periods and the 
contact parameters $\gamma_C^0 = 10, 10^2, 10^3$. 
In the top left plot, we observe that the solid is pulled down until it reaches $\Gamma_w$
at $t\approx 1.5\cdot 10^{-3}$. After a short \emph{contact} period, it is released again due to its
elastic properties
before it reaches the obstacle for a second time at $t\approx 5\cdot 10^{-3}$. 

The contact condition $d\cdot n_w\leq 0.25$ is only significantly violated
for the smallest contact parameter, where $d_{\min}$ reaches a minimum value of around
$-8\cdot 10^{-4}$, {see the zoom-in around the contact interval on the top left}. 
This value is more than an order of magnitude smaller than the 
mesh size in vertical direction $h\approx 1.4\cdot 10^{-2}$. For the larger values of $\gamma_C^0$, 
the minimal value of $d_{\min}$ is even much closer to zero. 
On the other hand, we observe that even for the largest value $\gamma_C^0=10^3$, the contact condition
is slightly relaxed, allowing for very small overlaps of solid and contact line.

In the second row of Figure~\ref{fig.vertDefandpress}, we observe that the pressure shows 
a peak at the beginning ($t\approx 1.5\cdot10^{-3}$ and $t\approx 5\cdot10^{-3}$)
of the contact periods, followed by some small oscillations. The peak is caused by the 
fluid dynamics and will be discussed below. The oscillations get smoother for larger values of 
$\gamma_C^0$ and are barely visible for $\gamma_C^0=10^3$.

Similarly, the contact force $J_{P_\gamma}$ shows oscillations for $\gamma_C^0=10$ and a much smoother 
behaviour for $\gamma_C^0\geq10^2$. Note that this does not contradict the stability result in 
Theorem~\ref{theo.stability}, where we have assumed that $\gamma_C^0$ is large enough.
{The relatively large value for $\gamma_C^0$ that is needed here is due to the anisotropic 
cells that appear in some of the time-steps, when using the \textit{locally modified finite element method}. 
In the absence of extreme anisotropies a value of $\gamma_C^0 \approx 1$ seems to be enough to obtain stable
numerical results.
The optimal choice of the 
contact parameter $\gamma_C$ in the context of anisotropic cells is subject to future research.

If $\gamma_C^0$ is chosen large enough, the contact force is roughly 
of the same size for different $\gamma_C^0$. This is in agreement with the observations of Chouly et al. for the 
case of a pure solid problem~\cite{ChoulyMlikaRenard18}, who showed that the consistency
of the method makes the choice of the contact parameter much less sensitive compared to a pure penalty method.}

\paragraph{Investigation of the pressure peak}

The pressure peak at the beginning of the contact interval can be explained as follows. 
As the fluid does not ''see`` the obstacle before reaching it, the solid is pulled down 
towards it without reducing
its velocity. At the moment when the obstacle is reached, 
its vertical velocity $\dot{d}\cdot n_w$ has to decrease to zero in an instant. Due to the continuity of 
velocities, the same happens for the fluid velocity $u\cdot n_w$ at the interface, and due to the 
incompressibility constraint the velocity has to change globally in the fluid domain $\Omega_f(t)$. 
The pressure can be seen as a
Lagrange multiplier and more specifically as sensitivity of a (jumping) energy functional 
with respect to the incompressibility constraint, which explains the peak.

To substantiate this explanation numerically, we add an artificial penalty for the velocity 
on the sub-domain $\Omega_f^0$
below the contact line
\begin{align*}
 S_a(u,v) :=\gamma_a^0 h^{-2} (u, v)_{\Omega_f^0}.
\end{align*}
For $\gamma_a:= \gamma_a^0 h^{-2} \to \infty$, the fluid velocity is driven to zero below the obstacle. As 
this is already the case before contact, no abrupt changes in the fluid velocity are expected at the moment of
the impact. {Note that the problem {with $\gamma_a^0 > 0$} is purely artificial, as the pressure 
mean values are still applied on 
the whole fluid boundary, including the boundary of $\Omega_f^0$.}

\begin{figure}[bt]
  \begin{minipage}{0.53\textwidth}
    \includegraphics[width=\textwidth]{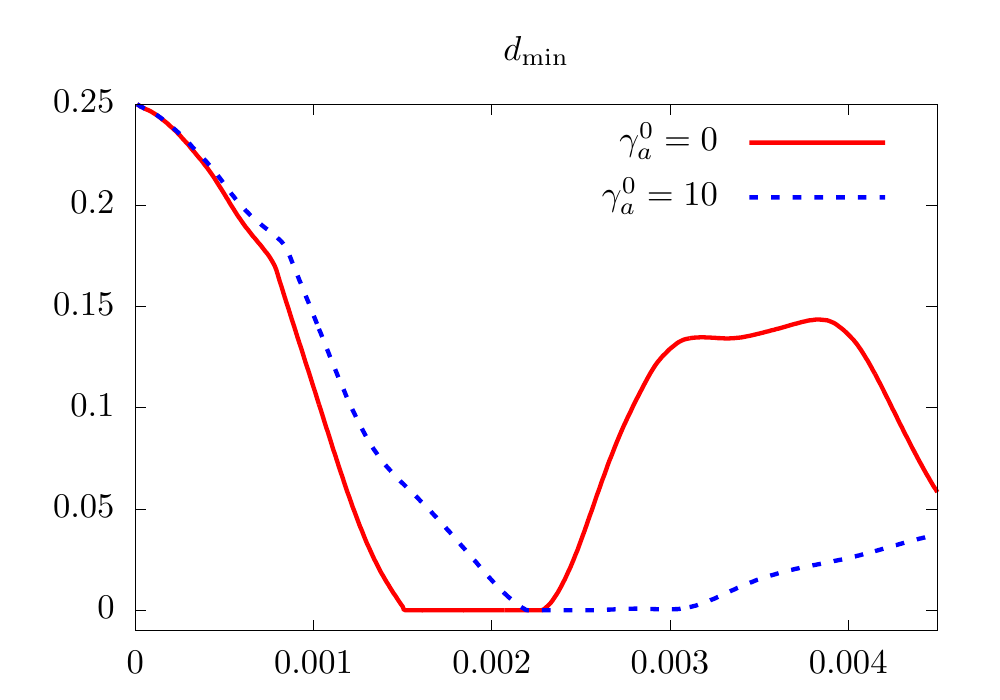}
    \end{minipage}
     \hspace{-0.6cm}
   \begin{minipage}{0.53\textwidth}
  \includegraphics[width=\textwidth]{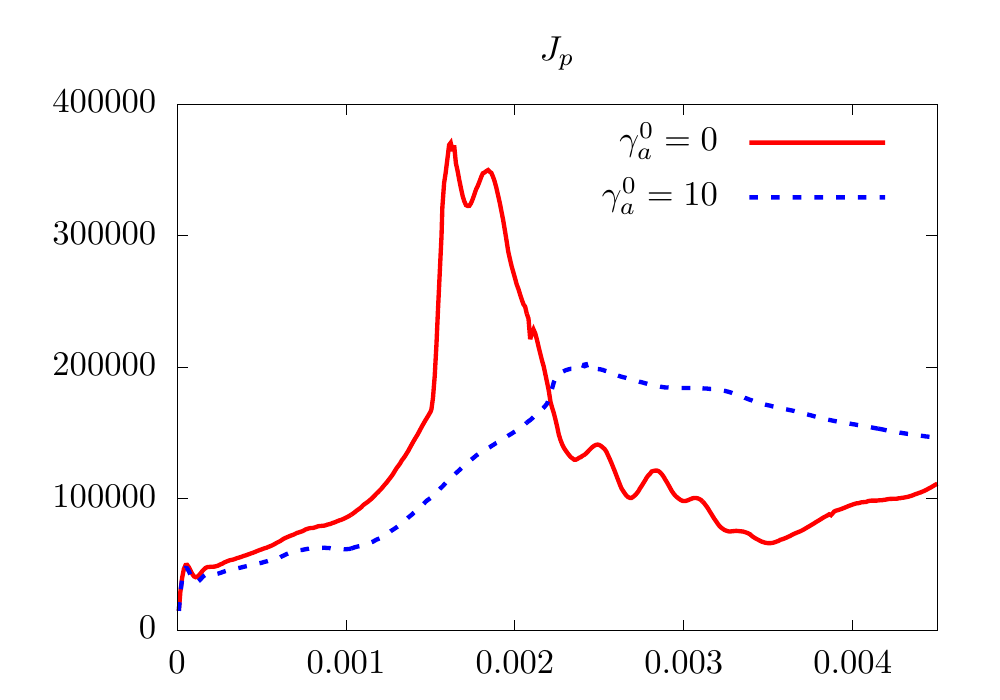}
 \end{minipage}
 \caption{\label{fig.gamma_a} Minimal distance $d_{\min}$ to the contact line $\Gamma_w$ (left)
 and mean pressure over the interface $J_p$ (right), plotted
 over time with an artificial penalty for the velocity 
 below the {virtual obstacle}.}
 \end{figure}

{
In Figure~\ref{fig.gamma_a} we compare the minimal distance and pressure functionals for computations 
without penalty ($\gamma_a^0=0$) to results
for $\gamma_a^0=10$.} First, we note that the vertical displacement is significantly 
influenced by the penalty, which has to be expected
 as the fluid dynamics are altered. {Contact happens later at $t\approx 0.0022$ with the artificial penalty}.
Moreover, we observe indeed that the initial pressure peak {at the time of impact} is significantly
reduced for $\gamma_a^0=10$.

\subsection{Contact problem}
\label{sec.contact}

Next, we study a problem, where it comes to real contact with 
{the wall} $\Gamma_w=\{(x,y)\in \Omega, y=0.25\}$. At time $t=0$, we define
\begin{align*}
 \Omega_f(0) := (0,1)\times (0.25,0.5), 
 \quad \Gamma(0) := (0,1) \times {0.5}, \quad \Omega_s(0) := (0,1)\times (0.5,0.6).
\end{align*}
Below $\Gamma_w$, we define a fixed artificial fluid domain $\Omega_f^C:= (0,1)\times (0,0.25)$.

We apply again a pressure mean value $\overline{P}$ on the lateral boundaries $\Gamma_{f,\text{left}}$ 
and $\Gamma_{f,\text{right}}$ of the 
fluid domain $\Omega_f(t)$. As the size of $\Gamma_{\text{left}}$ and $\Gamma_{\text{right}}$ is smaller 
and the viscous fluid forces acting against the ``real'' contact 
are stronger than in the previous example, we have to set 
a larger pressure force $\overline{P}=3\cdot 10^5$ in order to obtain contact. On the other hand, 
the contact was never released again in our numerical experiments, when we used
this constant boundary force for all times.
Therefore, we decrease $\overline{P}$ linearly
from $t=10^{-3}$ on until it reaches zero at $t=1.2\cdot 10^{-3}$.
{In order to avoid the issues related to no-slip conditions and contact, we use slip-interface conditions
first, i.e.$\,$the Variational Formulation~\ref{varForm_artF} with ${\cal A}_{\text{slip,FSI}}^C$. Unless stated 
differently,
$\lambda$ is chosen as the jump of {numerical stresses} $\jump{\widetilde{\sigma}_{n,\text{slip}}}$.}
Moreover, we use again a Cartesian mesh that
consists of 5120 elements, a time step $\delta t = 10^{-5}$ and temporal pressure stabilisation 
with $\gamma_{pt}=10^{-2}$.
Unless explicitly stated, the Nitsche parameters are chosen
as $\gamma_{\text{fsi}}^0 = \gamma_C^0 = 10^3$ and the penalty in the artificial fluid as $\gamma_a^0=10^2$.

\begin{figure}[bt]
 \begin{minipage}{0.53\textwidth}
  \includegraphics[width=\textwidth]{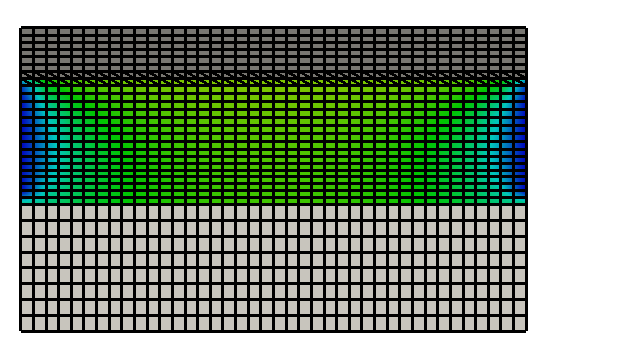}
 \end{minipage}
 \hspace{-1cm}
  \begin{minipage}{0.53\textwidth}
    \includegraphics[width=\textwidth]{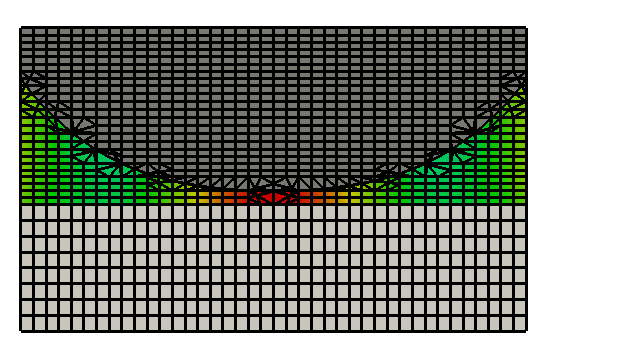}
 \end{minipage}
 \begin{minipage}{0.53\textwidth}
  \includegraphics[width=\textwidth]{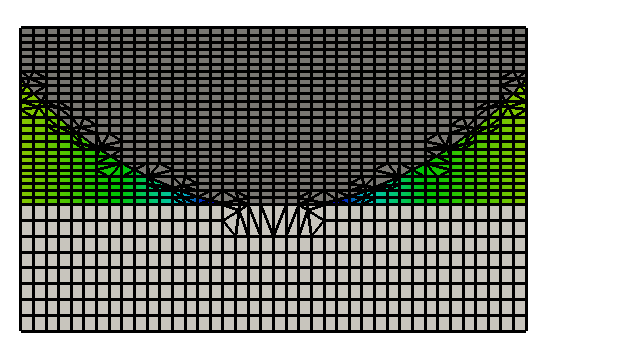}
 \end{minipage}
 \hspace{-1cm}
  \begin{minipage}{0.53\textwidth}
    \includegraphics[width=\textwidth]{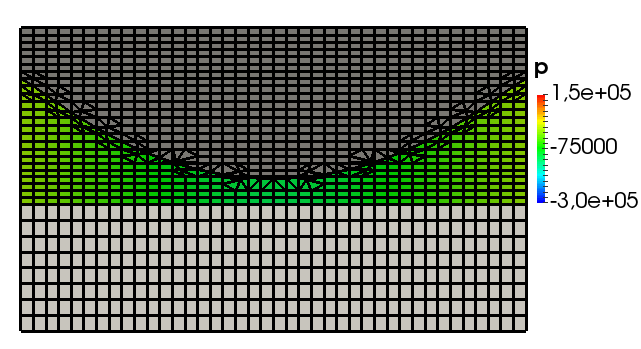}
 \end{minipage}
 \caption{\label{fig.vtks} Illustration of the contact problem at four time instances on a coarse mesh: 
 $t=0$ (top left), $t=1.2\cdot10^{-3}$ (top right), $t=2\cdot10^{-3}$ (bottom left) and $t=2.5\cdot10^{-3}$ (bottom right). 
 The grey part corresponds to the structure $\Omega_s(t)$, the white part is the artificial fluid 
 $\Omega_f^C$. In the fluid domain $\Omega_f(t)$, values of the pressure $p$ are visualised.
 }
\end{figure}

The results on a coarser
mesh are illustrated in Figure~\ref{fig.vtks} at four time {instants}. Contact happens after 
the pressure on the lateral boundaries is released, as the solid
continues moving downwards for some time. During contact, there is a very 
small overlap of the solid with the 
artificial fluid $\Omega_f^C$. As the overlap is of order $10^{-5}$, 
it can barely be seen in the bottom left picture. 
Notice however the triangular
cells in $\Omega_f^C$ that are used only, when a patch is cut by $\Gamma(t)$.

 \paragraph{Comparison of the two contact formulations}

\begin{figure}[bt]
 \centering
  \begin{minipage}{0.5\textwidth}
  \includegraphics[width=\textwidth]{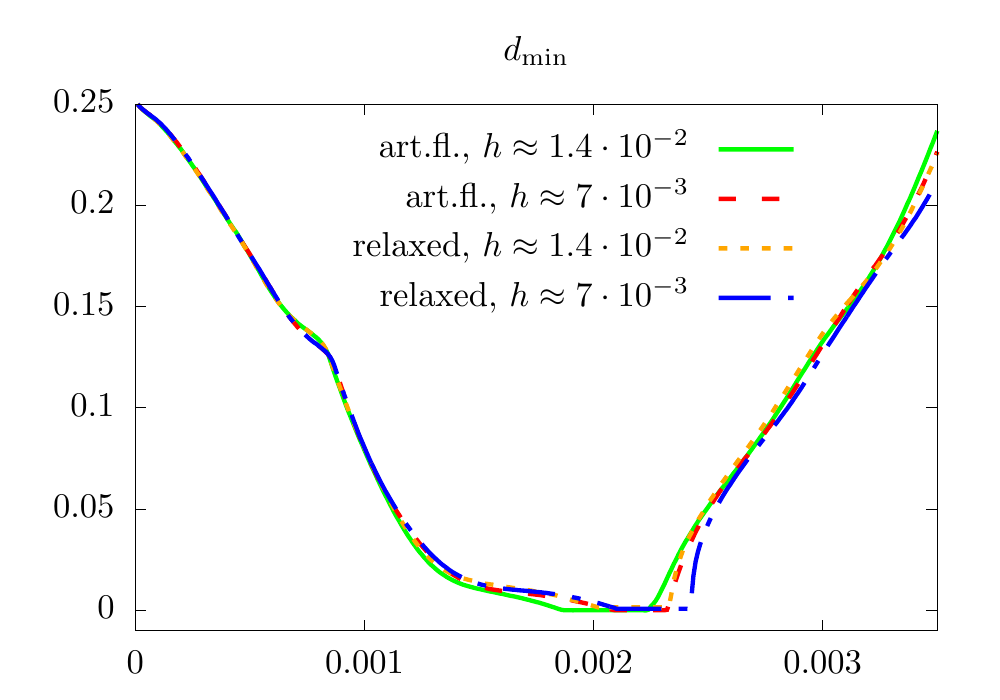}
 \end{minipage}
 \hspace{-0.5cm}
  \begin{minipage}{0.5\textwidth}
    \includegraphics[width=\textwidth]{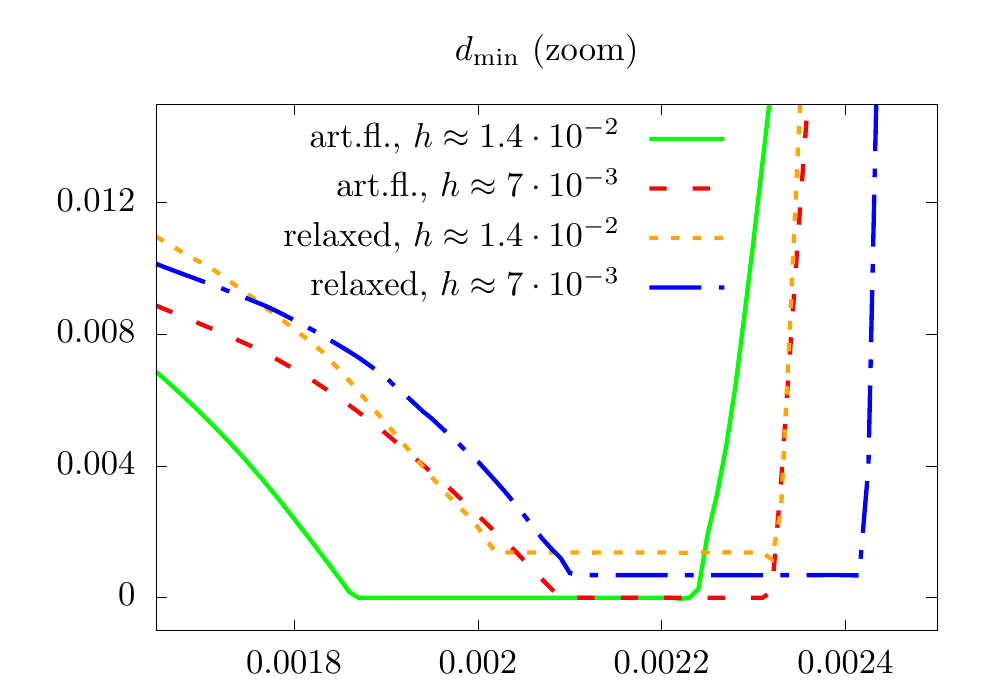}
 \end{minipage}
    \begin{minipage}{0.5\textwidth}
  \includegraphics[width=\textwidth]{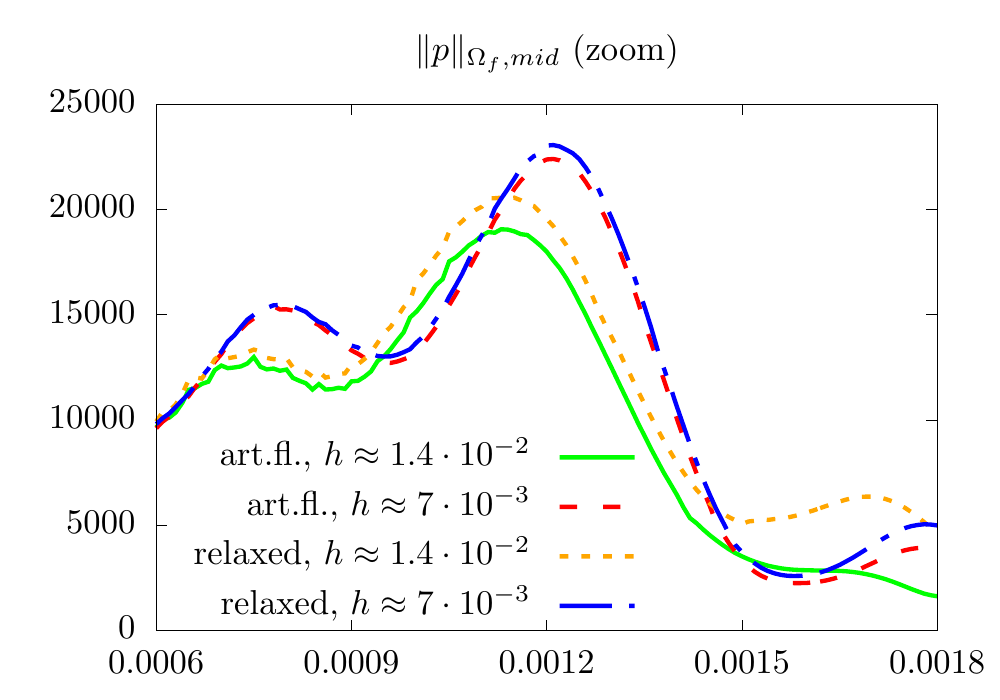}
 \end{minipage}
 \caption{\label{fig.gap} Comparison of the relaxed and the artificial fluid contact formulation. 
 \textit{Top}: Minimal distance $d_{\min}$ to $\Gamma_w$. 
 \textit{Right:} Zoom-in around the contact interval. \textit{Bottom:} Pressure norm $\|p\|_{L^2(\Omega_f,\text{mid})}$
 over the central part of the fluid domain before contact over time.}
 \end{figure}

{First, we compare the two contact strategies derived in Section~\ref{sec.fullContact}, i.e.$\,$the 
relaxed contact formulation introduced in Section~\ref{sec.relaxed}
with a small gap of size $\epsilon(h)=h/10$ between the solid and $\Gamma_w$}
and the strategy using an artificial fluid derived in Section~\ref{sec.artf}.
{A comparison of {the} results for $d_{\min}$, the pressure norm
$
\|p\|_{L^2(\Omega_{f,\text{mid}}(t))}$, where 
$$
\Omega_{f,\text{mid}}(t) := \{x\in \Omega_f(t), 0.4\leq x_1\leq 0.6\}
$$
denotes the central part of the fluid domain, 
and the contact force $J_{P_\gamma}$ are shown in Figure~\ref{fig.gap} on two different meshes
with 5120 and 20480 elements, respectively}.

{First, we observe from the plots in the top row that the
interface stays at a distance {to $\Gamma_w$} of 
about $\epsilon\approx h/10\approx 1.4\cdot 10^{-3}$ on the coarser and 
$\epsilon\approx h/10\approx 7\cdot 10^{-4}$ on the finer mesh for the relaxed formulation. The
much smaller overlap with $\Omega_f^C$ in the artificial fluid formulation is not visible, not 
even in the zoom-in on the right. 

While the curves for $d_{\min}$ look similar in the global picture (left), the zoom-in shows
significant differences already before the impact.
The contact happens earlier for the artificial fluid formulation: on the coarser mesh 
at time $t_{C,a}= 1.87\cdot 10^{-3}$ compared to $t_{C,r}= 2.02\cdot 10^{-3}$ for the relaxed 
formulation. This deviation is already much smaller on the finer mesh, where $t_{C,t}-t_{C,a} = 2\cdot 10^{-5}$.
The reason for this deviation is that in the 
artificial fluid formulation the wall $\Gamma_w$ is only asymptotically for $\gamma_a\to\infty$ impermeable for 
the fluid. Therefore, in practice, the fluid forces acting against the contact, in particular the pressure $p$,
are smaller for this formulation. As $\gamma_a=\gamma_a^0 h^{-2}\to\infty$ for $h\to 0$ 
the difference is significantly reduced on the finer mesh.

To substantiate this explanation, we plot the pressure norm $\|p\|_{L^2(\Omega_{f,\text{mid}}(t))}$ over 
the central part 
of the fluid domain $\Omega_{f,\text{mid}}(t)$ and the time period before the impact in the bottom left figure.
The maximum value of the norm on the coarser mesh at time $t=1.1\cdot 10^{-3}$ is approximately $20.563$ for the relaxed formulation 
and about $19.065$ for the artificial fluid version. {On the finer mesh, the pressure values are much closer.}
After that time
the functional values decrease because the domain $\Omega_{f,\text{mid}}(t)$ gets smaller.

In the next paragraph, we will study the performance of both contact formulations under mesh refinement.
}

 \paragraph{Convergence under mesh refinement}

\begin{figure}[bt]
   \begin{minipage}{0.5\textwidth}
  \includegraphics[width=\textwidth]{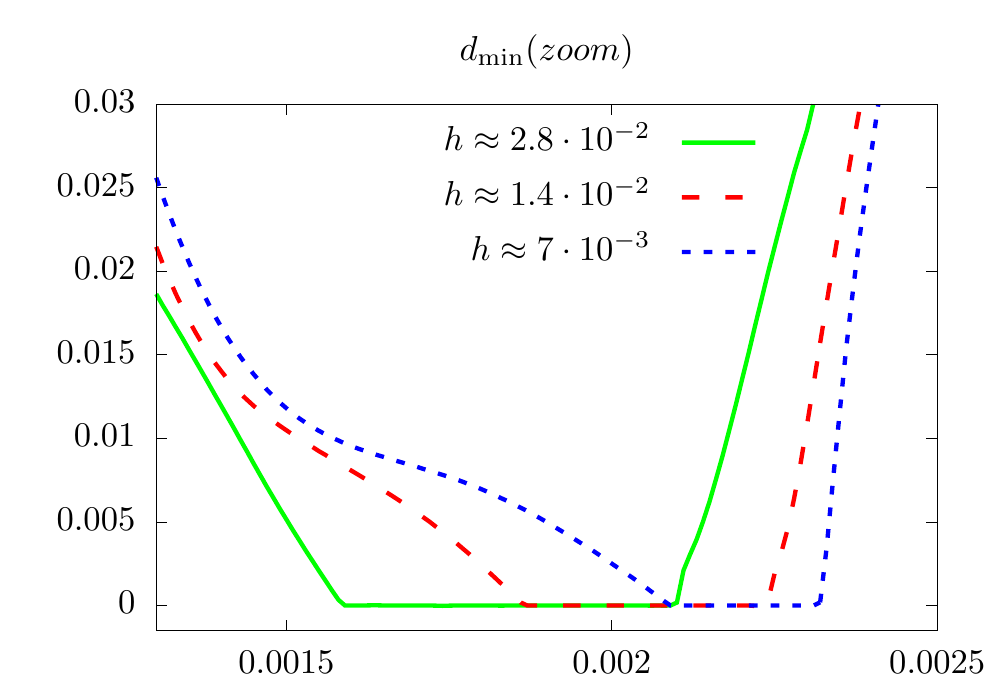}
 \end{minipage}
 \hspace{-0.5cm}
  \begin{minipage}{0.5\textwidth}
    \includegraphics[width=\textwidth]{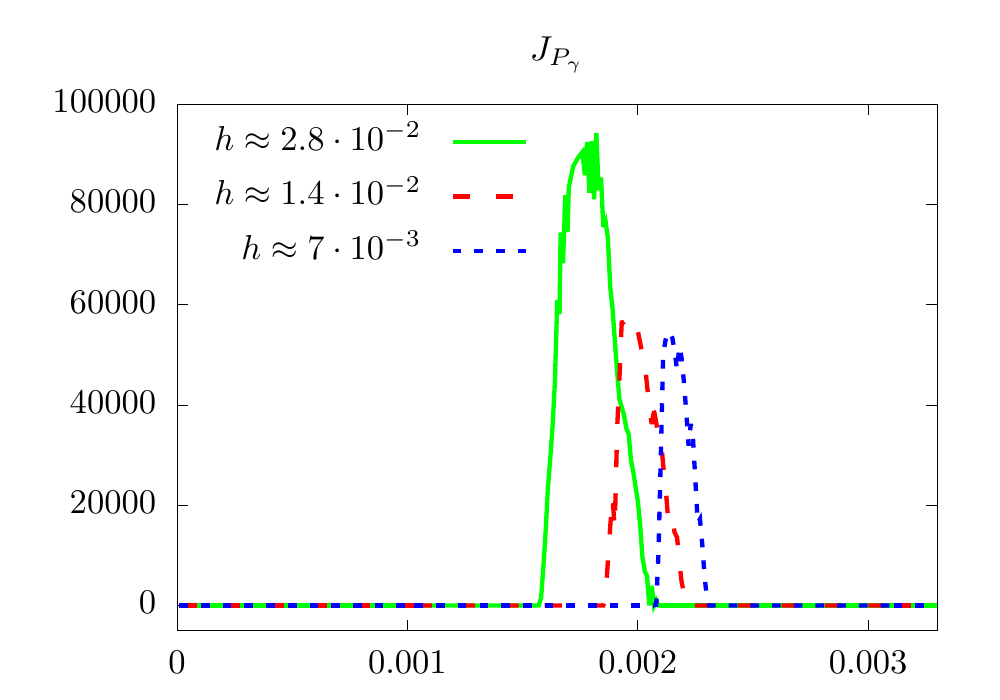}
 \end{minipage}
 
    \begin{minipage}{0.5\textwidth}
  \includegraphics[width=\textwidth]{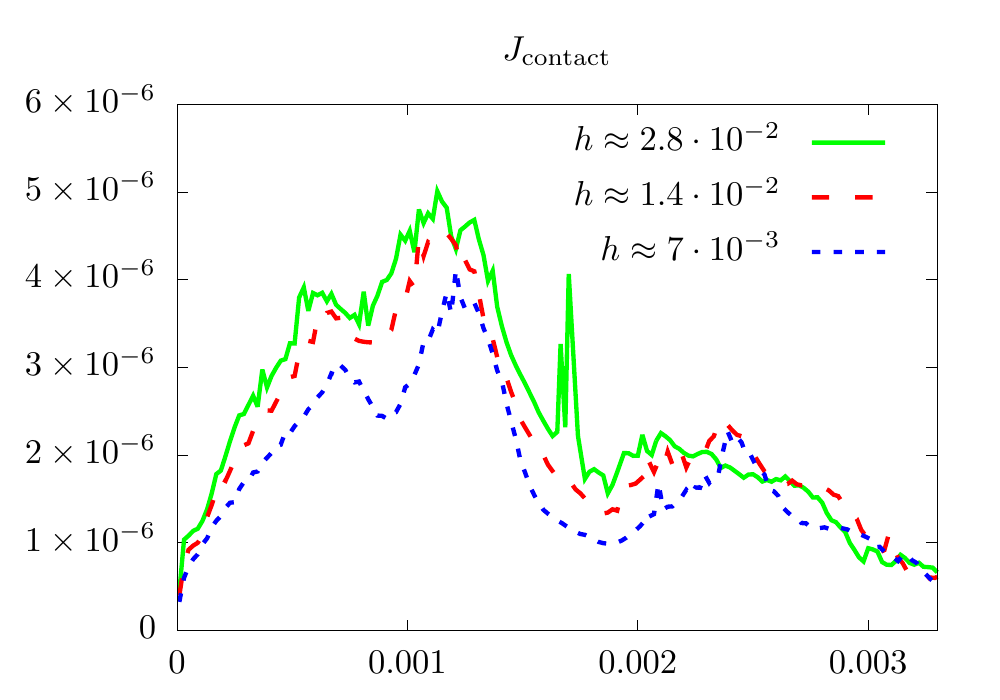}
 \end{minipage}
 \hspace{-0.5cm}
  \begin{minipage}{0.5\textwidth}
    \includegraphics[width=\textwidth]{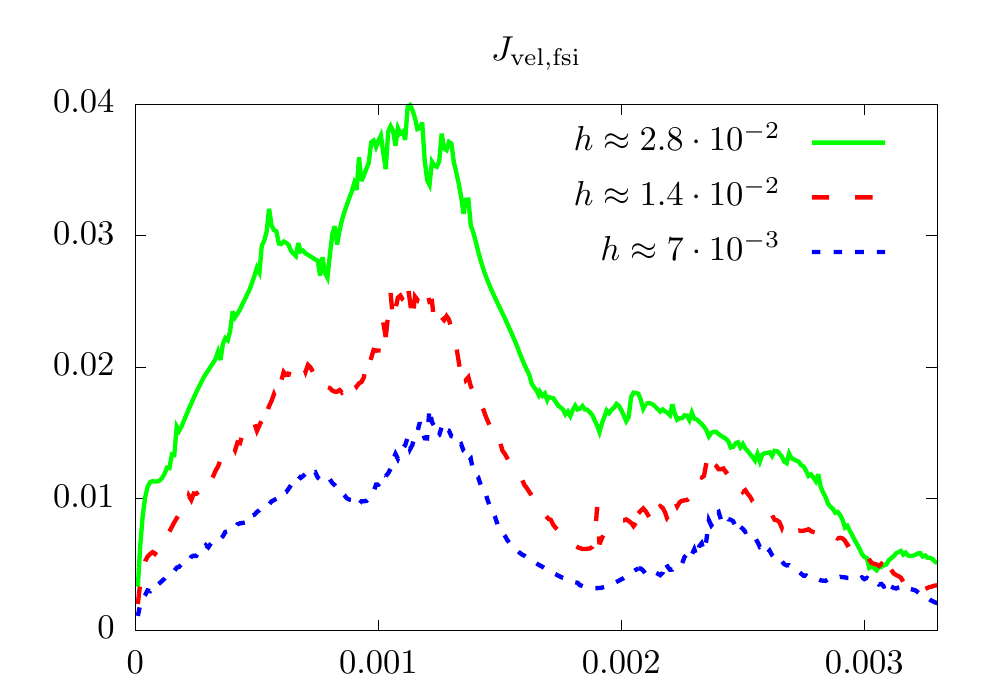}
 \end{minipage}
 \caption{\label{fig.convSpace} Convergence studies under mesh refinement for the \textit{artificial fluid formulation} 
 by means of the following 
 functionals over time:
 \textit{Top left:} Minimal distance $d_{\min}$ of $\Gamma(t)$ to $\Gamma_w$, \textit{top right}:
 contact force $J_{P_{\gamma}}$, \textit{Bottom left}:
 Fulfilment of the contact condition
 $J_{\text{contact}}$. \textit{Bottom right}: 
 Continuity of velocities $J_{\text{vel,fsi}}$.}
 \end{figure}

\begin{figure}[bt]
   \begin{minipage}{0.5\textwidth}
  \includegraphics[width=\textwidth]{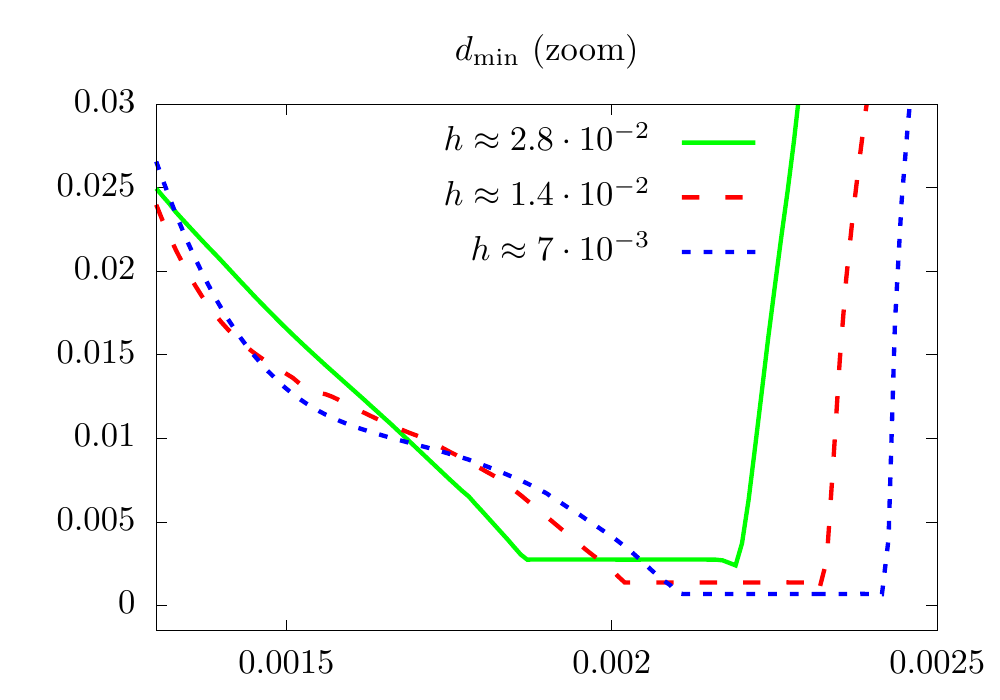}
 \end{minipage}
 \hspace{-0.5cm}
  \begin{minipage}{0.5\textwidth}
    \includegraphics[width=\textwidth]{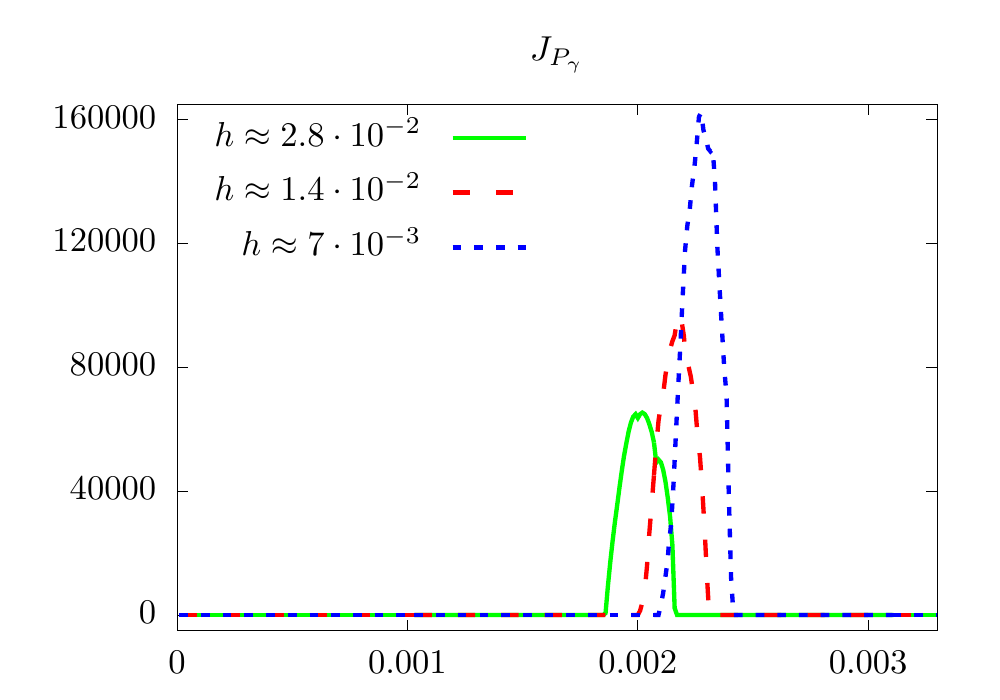}
 \end{minipage}
 
    \begin{minipage}{0.5\textwidth}
  \includegraphics[width=\textwidth]{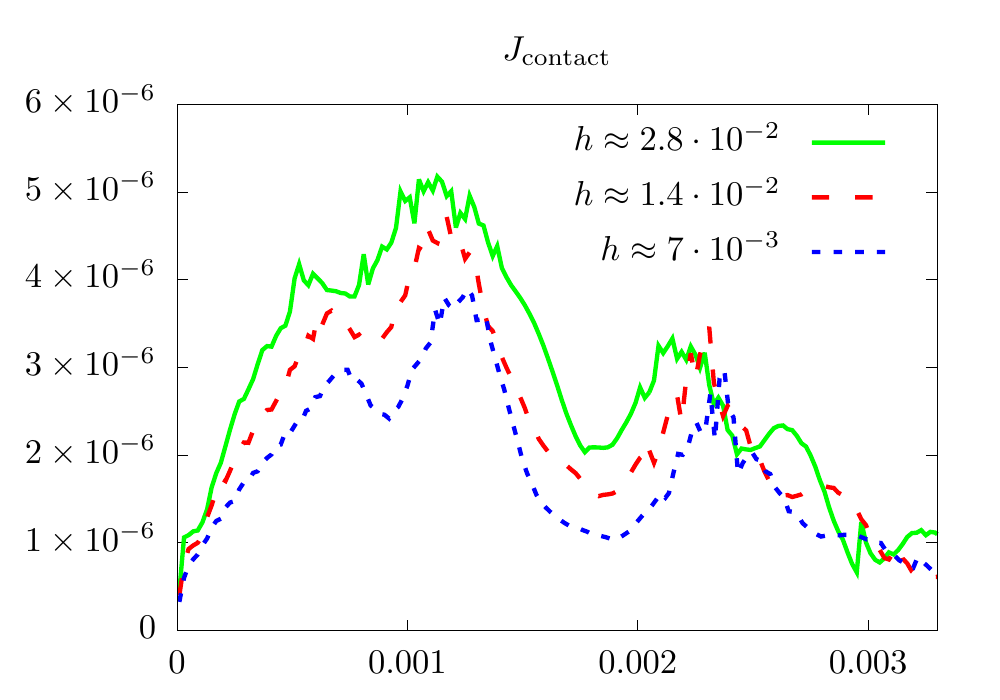}
 \end{minipage}
 \hspace{-0.5cm}
  \begin{minipage}{0.5\textwidth}
    \includegraphics[width=\textwidth]{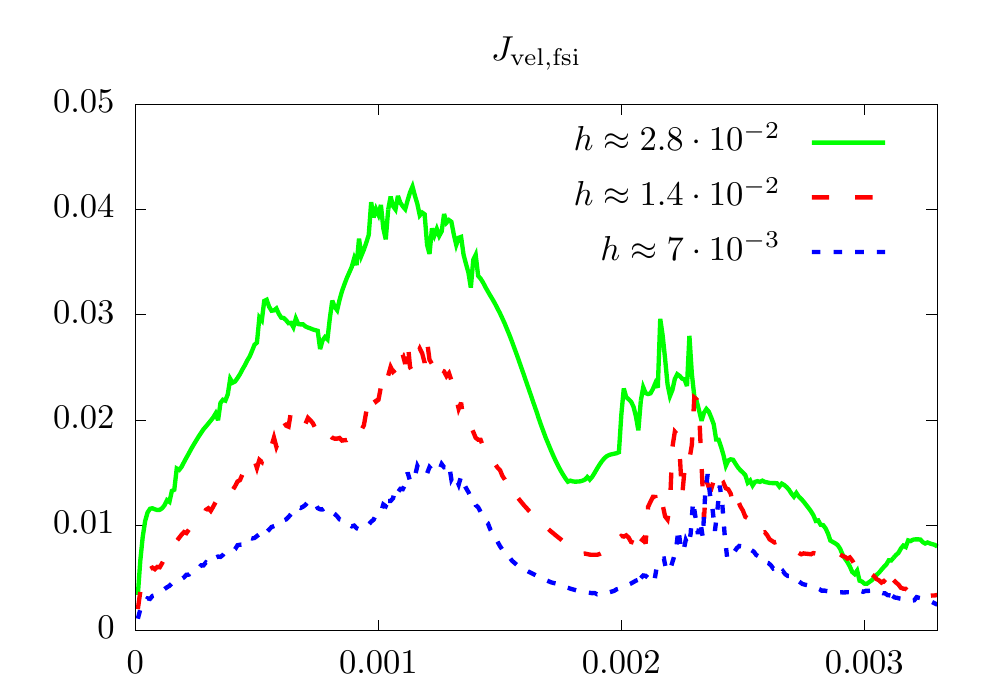}
 \end{minipage}
 \caption{\label{fig.convSpace_rel} Convergence studies under mesh refinement for the \textit{relaxed contact formulation} 
 by means of the following 
 functionals over time:
 \textit{Top left:} Minimal distance $d_{\min}$ of $\Gamma(t)$ to $\Gamma_w$, \textit{top right}:
 contact force $J_{P_{\gamma}}$, \textit{Bottom left}:
 Fulfilment of the contact condition
 $J_{\text{contact}}$. \textit{Bottom right}: 
 Continuity of velocities $J_{\text{vel,fsi}}$.}
 \end{figure}

We solve the same problem on three different meshes with 1.280, 5.120
and 20.480 mesh elements, where the finer
meshes are obtained from the coarsest one by global mesh refinement. The plots of the 
functionals $d_{\min}$ and $J_{P_{\gamma}}$
as well as the functionals
\begin{align*}
 J_{\text{contact}} &:= \frac{\left\| \gamma_C^{1/2} [P_{\gamma}(\lambda,d)]_+ 
 + \gamma_C^{-1/2} \lambda\right\|_{\Gamma(t)}}{\overline{\|\lambda\|}_{\Gamma}},\quad
 J_{\text{vel,fsi}}:=\frac{\|(\dot{d}-u)\cdot n\|_{\Gamma_{\text{fsi}}(t)}}
 {\overline{\|u\cdot n\|}_{\Gamma_{\text{fsi}}} + \overline{\|\dot{d}\cdot n\|}_{\Gamma_{\text{fsi}}}}
\end{align*}
measuring the fulfilment of the contact condition and 
the continuity of velocities on the part 
$\Gamma_{\text{fsi}}(t)$ of $\Gamma(t)$ that is not in contact with $\Gamma_w$
\begin{align*}
\Gamma_{\text{fsi}}(t) := \left\{x\in \Gamma(t) \, \big| \, P_\gamma(\lambda,d)(x) \leq 0\right\}.
\end{align*}
{are shown in Figure~\ref{fig.convSpace} for 
the artificial fluid formulation and in Figure~\ref{fig.convSpace_rel} for 
the relaxed contact formulation over time.
The quantities $\overline{\|\cdot\|}_*$ that are used to scale the functionals 
are temporal averages of the respective norms over the interval $I=[0,0.004]$, computed on the 
finest grid.}

First, we observe for both formulations in the plots on the top left that the contact happens later, the finer the 
discretisation is, as the fluid forces which act against the
closure of the fluid channel are better resolved on the finer meshes (see also Figure~\ref{fig.gap} and the related discussion above). 

{The curves for the contact force
$P_{\gamma}$ on the top right of both figures show significant differences between the two formulations.
While the functional values seem to converge for the artificial fluid formulation (if we neglect the time shift),
the contact force gets larger under mesh refinement for the relaxed formulation.
The larger values for the relaxed formulation are due to the presence of the fluid forces $\sigma_{f,n}$
during the whole contact interval, that are not penalised in this formulation. As a Lagrange multiplier for 
the incompressibility constraint, the continuous pressure $p$ gets singular
when it comes to contact. The discrete pressure $p_h$ gets 
larger and larger under mesh refinement in our computations.}

{On the other hand, the functional $J_{\text{contact}}$, that measures the difference between
$-\gamma_C^{1/2} [P_\gamma]_+$ and $\gamma_C^{-1/2} \lambda$ decreases under mesh refinement for both 
formulations. Besides the differences in the contact force $J_{P_\gamma}$,
the functional values on each of the mesh levels are actually very similar for the two formulations. 
The reason must be 
that the fluid forces $\sigma_{f,n}$ enter in both $\lambda$ and $[P_{\gamma}]_+$. We conclude that the 
increase in the functional $J_{P_\gamma}$ seems not to be an issue for the contact dynamics.

Both $J_{\text{contact}}$ and the functional $J_{\text{vel,fsi}}$ are controlled by the 
stability estimate in 
Theorem~\ref{theo.stability} for $\theta=1$. Although the parameter $\theta=0$ 
is used here, we observe that both
functionals decrease with mesh refinement before and during contact. While the 
convergence for the contact functional
is quite slow, the values of the velocity functional indicate a convergence order ${\cal O}(h^{\alpha})$
with $0.5\leq \alpha\leq 1$ for both formulations. Note that in contrast to the 
term $J_{\text{contact}}$, $ J_{\text{vel,fsi}}$
is controlled in the
stability estimate in Theorem~\ref{theo.stability} even with the pre-factor
$(\gamma_{\text{fsi}}^0 \mu_f)^{1/2} h^{-1/2}$.}

\paragraph{Flux formulations}
 
 {
 Next, we compare the different choices for $\lambda$. We show results exemplarily for the artificial fluid
 formulation with slip interface conditions. We will compare results using the jump of 
 stresses $\lambda=\jump{\sigma_{n,\text{slip}}}$~\eqref{LagrCompleteSlip},
 the jump of the numerical fluxes $\lambda=\jump{\widetilde{\sigma}_{n,\text{slip}}}$~\eqref{lambdaslip}
 and the extended fluxes
 $\lambda = \jump{\widetilde{\sigma}_{n,\text{slip}}} + \tau^T \sigma_s n (\tau \cdot n_w)$
 (see \eqref{eq:slip-perturbed}). As the results for the latter two choices are nearly 
 identical in this example, we show here only plots for the jump of stresses and the 
 (non-extended) numerical fluxes.
  We use the artificial fluid formulation (Variational Formulation~\ref{varForm_artF}) 
 and the previously used mesh with 5120 elements.} 
 
\begin{figure}[bt]
 \centering
   \begin{minipage}{0.5\textwidth}
  \includegraphics[width=\textwidth]{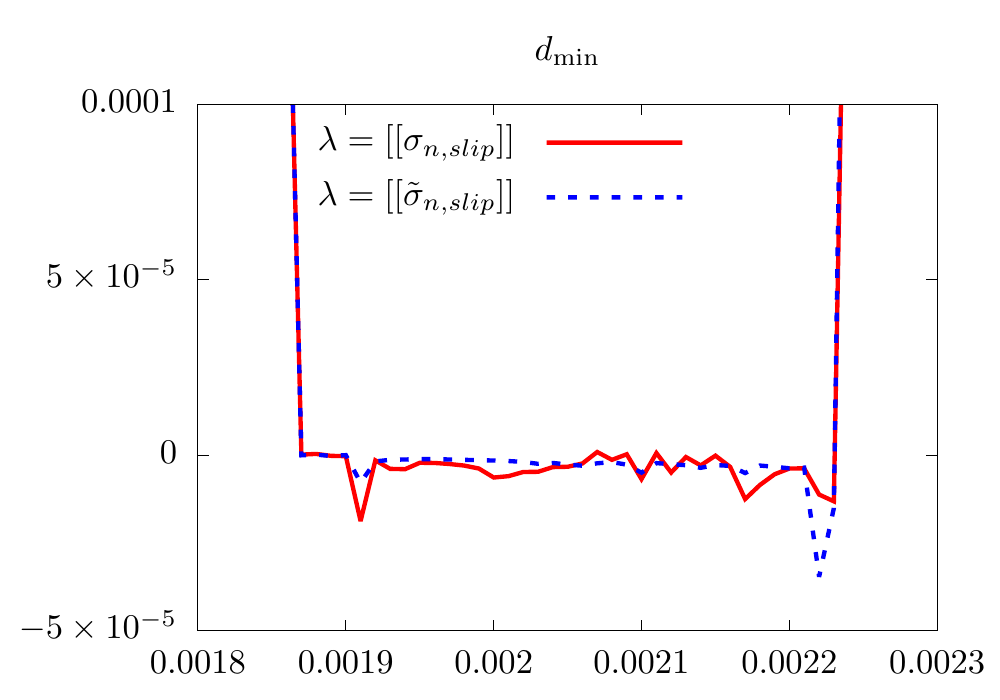}
 \end{minipage}
 \hspace{-0.5cm}
  \begin{minipage}{0.5\textwidth}
  \includegraphics[width=\textwidth]{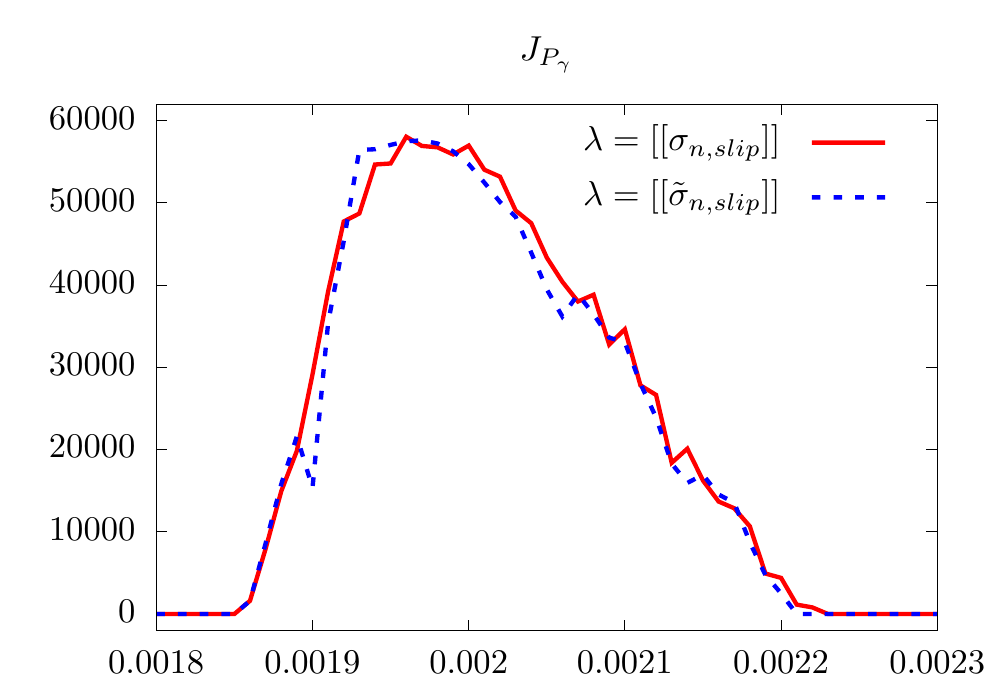}
 \end{minipage}
  \begin{minipage}{0.5\textwidth}
    \includegraphics[width=\textwidth]{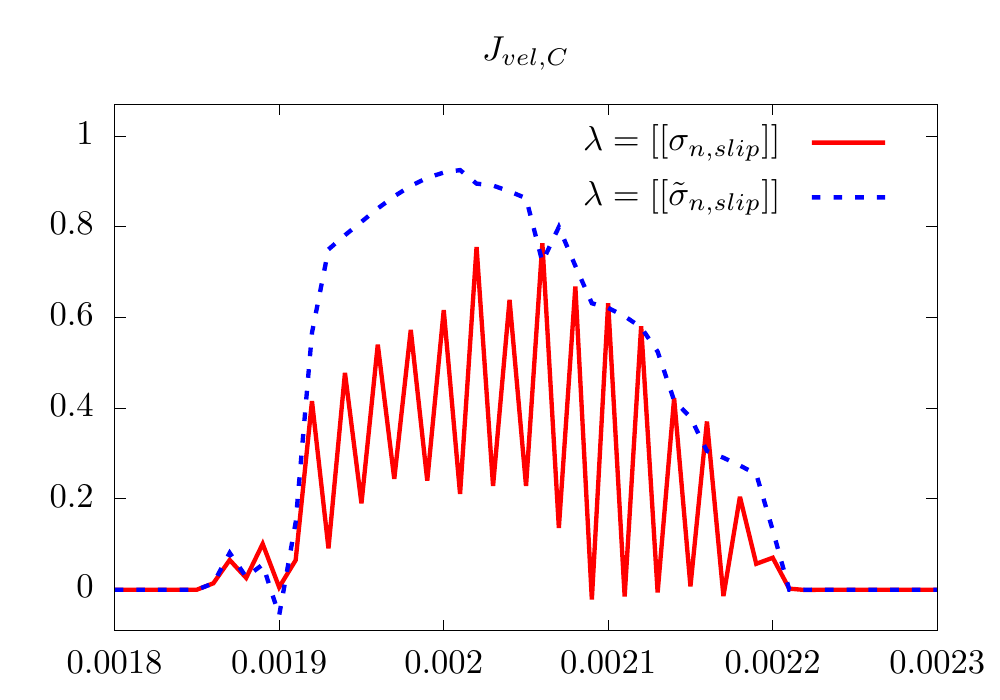}
 \end{minipage}
 \caption{\label{fig.dGfluxes} Comparison of the different possibilities to choose the fluxes $\lambda$. 
 Minimal distance $d_{\min}$ to $\Gamma_w$ (top left), contact force
 $J_{P_\gamma}$ (top right) and 
  velocity difference $J_{\text{vel},C}$ (bottom) integrated
 over the contact part $\Gamma_C(t)$ of $\Gamma(t)$
 over time for computations with {$\lambda=\jump{\sigma_{n,\text{slip}}}$ 
 and $\lambda=\jump{\widetilde{\sigma}_{n,\text{slip}}}$ for the artificial fluid formulation.}
 }
 \end{figure}

 In Figure~\ref{fig.dGfluxes}, we show the minimal distance $d_{\min}$ to $\Gamma_w$, the contact force
 $J_{P_{\gamma}}$
 and the integral over the velocity difference across the contact part $\Gamma_C(t)$ 
 of the interface over time
 \begin{align*}
 J_{\text{vel},C} := \int_{\Gamma_C(t)} (\dot{d}-u)\cdot n\, ds, \qquad
{\Gamma_C(t) := \left\{x\in \Gamma(t) \, \big| \, P_\gamma(\lambda,d)(x) \leq 0\right\}.}
\end{align*}
The fluid velocity $u$ is here artificial as it comes from $\Omega_f^C$. 
{When choosing $\lambda=\jump{\widetilde{\sigma}_{n,\text{slip}}}$, 
 we ensure that there is
 no feedback from this artificial velocity to the solid, see (\ref{dGflux_argumentation}). 
 For the jump of stresses $\lambda=\jump{\sigma_{n,\text{slip}}}$, we obtain a mixture of 
 the solid contact condition and the continuity of normal velocities and a feedback might result. This 
 follows analogously to the no-slip case, see~\eqref{mixture_wo_dGfluxes}.}
 
 In the left sketch of Figure~\ref{fig.dGfluxes}, we see that the minimal distance in the 
 {stress-based} formulation
 shows oscillations during the whole contact interval, especially in the second half. The interface jumps back
 and forth over the contact line many times. The curve corresponding to the formulation using discrete fluxes is 
 much smoother. Similarly,
 the contact force $J_{P_{\gamma}}$ looks smoother, when the flux formulation is used. 
 The reason
 for this behaviour is the mixture of the interface conditions during contact. On the bottom
 of Figure~\ref{fig.dGfluxes}, we see that
 the velocity difference $J_{\text{vel},C}$ shows wild oscillations for the stress formulation, while 
 it looks much smoother when using $\jump{\widetilde{\sigma}_{n,\text{slip}}}$. As the artificial 
 velocity in $\Omega_f^C$ 
 has no physical meaning, it is not a drawback that the absolute values of $J_{\text{vel},C}$ are larger. 
 {Due to the feedback of this velocity to the contact conditions, the oscillations 
 appear in the displacement as well.}

 On the other hand, we should mention that the 
 oscillations are relatively small. Especially those in $d_{\min}$ are almost by a factor $10^3$ smaller
 than the mesh size $h\approx 1.4\cdot 10^{-2}$ in vertical direction in this example and are
 therefore still acceptable.

 \paragraph{Influence of the contact parameter $\gamma_C^0$}

\begin{figure}[bt]
\centering
   \begin{minipage}{0.5\textwidth}
  \includegraphics[width=\textwidth]{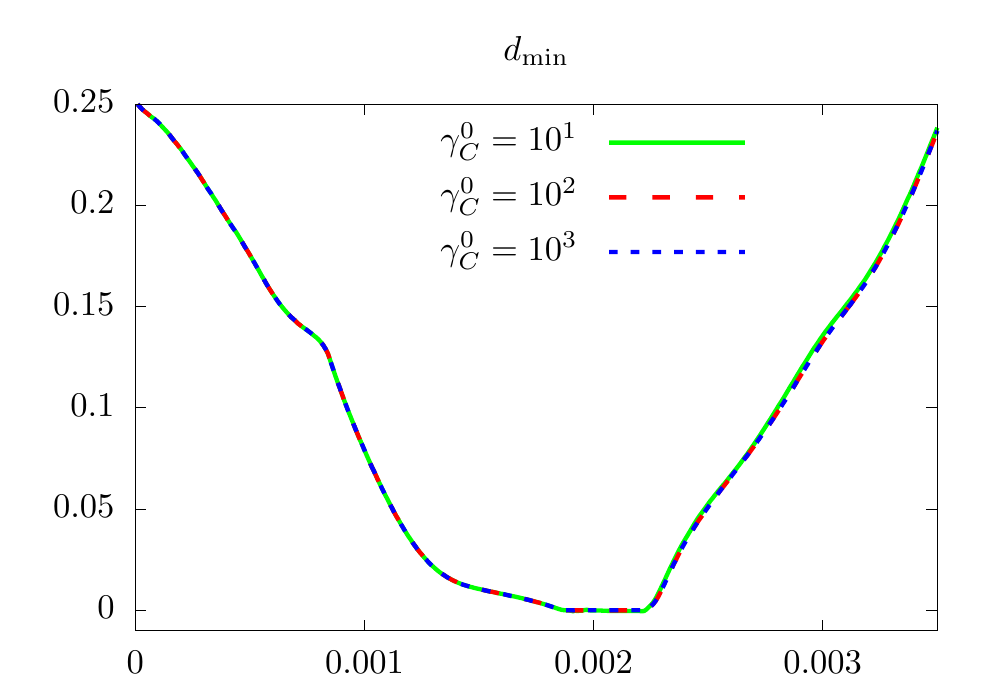}
 \end{minipage}
 \hspace{-0.5cm}
  \begin{minipage}{0.5\textwidth}
    \includegraphics[width=\textwidth]{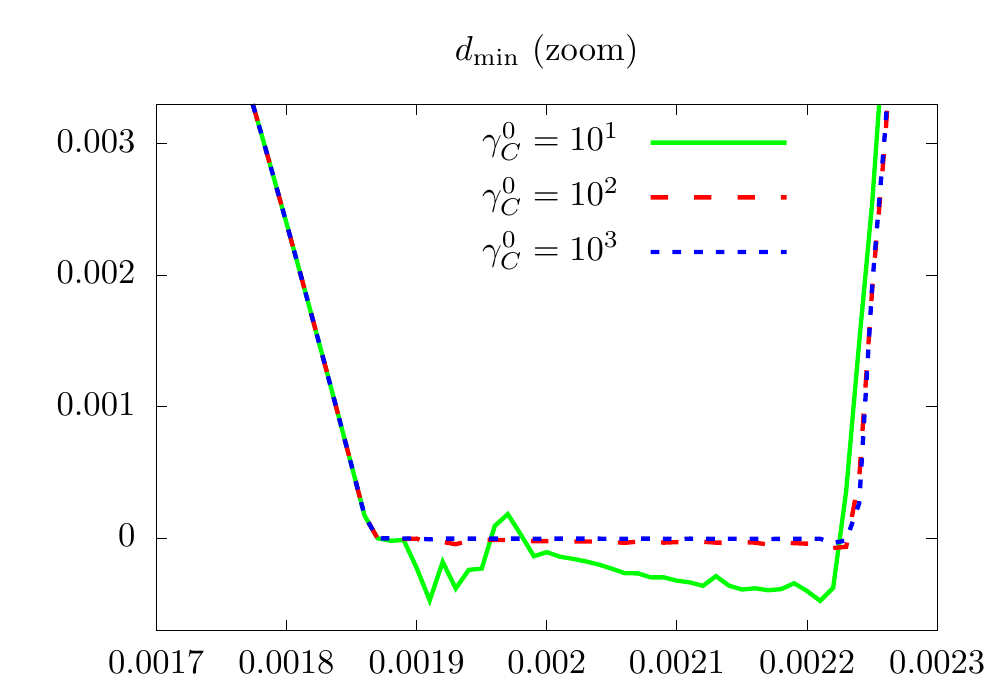}
 \end{minipage}
 
    \begin{minipage}{0.5\textwidth}
  \includegraphics[width=\textwidth]{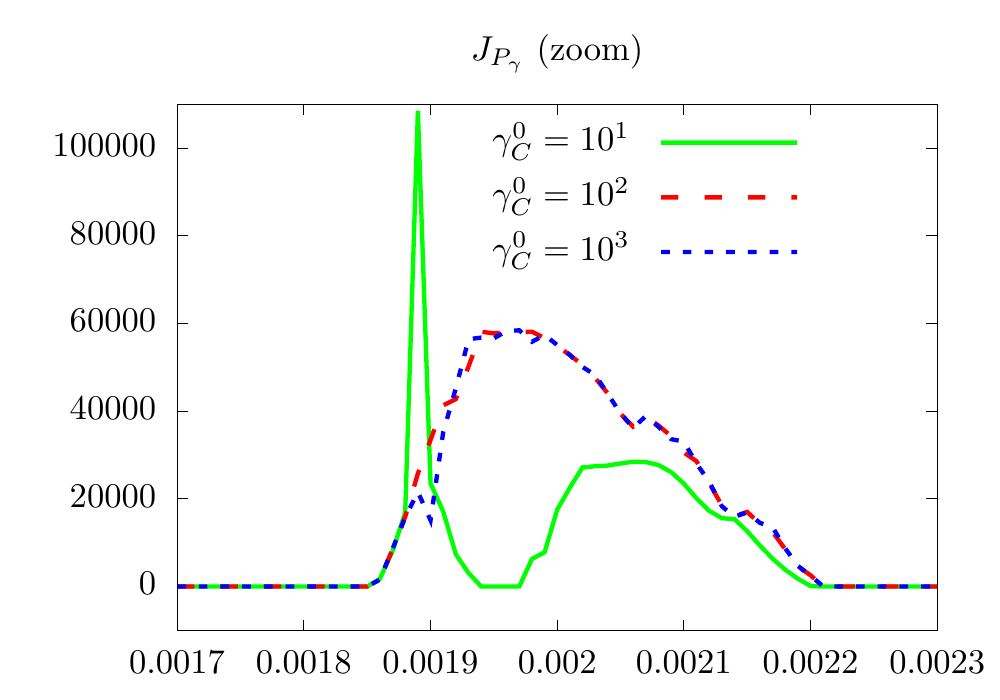}
 \end{minipage}
 \caption{\label{fig.contactpar} Parameter studies for the contact parameter $\gamma_C^0$: 
 Minimal distance of $\Gamma(t)$ to the {wall} $\Gamma_w$ (\textit{Top left}: Total time interval, 
 \textit{top right}: zoom around the contact interval)
 and contact force $J_{P_{\gamma}}$ (bottom) over time.}
 \end{figure}

{Next, we study the effect of different contact parameters $\gamma_C^0$ for the artificial
fluid formulation and $\lambda=\jump{\widetilde{\sigma}_{n,\text{slip}}}$ on the mesh with 5120 elements.
}
In Figure~\ref{fig.contactpar}, we show the 'minimal distance' $d_{\min}$ (top) and the {contact force
$J_{P_\gamma}$} over time for different contact parameters $\gamma_C^0$.
{The results are similar to the corresponding results for the virtual obstacle problem 
in Figure~\ref{fig.vertDefandpress}.
For} 
the smallest contact parameter
$\gamma_C^0=10$, the contact condition is violated throughout the contact 
interval ($d_{\min}<0$). The maximum overlap into the artificial fluid domain is again approximately by a factor
30 smaller than the mesh size $h\approx 1.4\cdot 10^{-2}$. 
This violation gets smaller, the larger the contact parameter is chosen. {The instabilities for the smallest 
parameter are still much better visible in the contact force $J_{P_{\gamma}}$. At time $t=1.89\cdot 10^{-3}$ the functional
shows a huge peak, as the contact condition $d\cdot n_w \leq g_0$ is severely violated and it vanishes 
from $t=1.94\cdot 10^{-3}$ to $t=1.97\cdot 10^{-3}$, when the contact is in fact shortly released.}

For the larger values {$\gamma_C^0\geq10^2$}, the curves are relatively smooth and very similar. {
Altogether, this
shows again that the assumption ``$\gamma_C^0$ sufficiently large'' in Theorem~\ref{theo.stability} is 
necessary in order to ensure stability.}

 \paragraph{Slip vs no-slip conditions}
 
\begin{figure}[bt]
   \begin{minipage}{0.5\textwidth}
  \includegraphics[width=\textwidth]{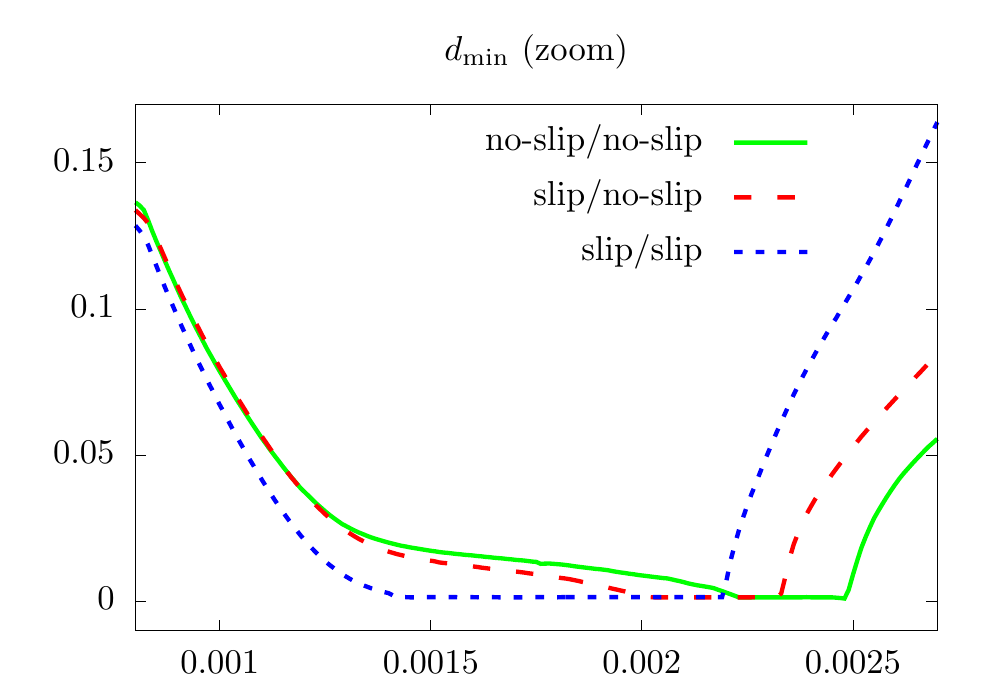}
 \end{minipage}
 \hspace{-0.5cm}
  \begin{minipage}{0.5\textwidth}
    \includegraphics[width=\textwidth]{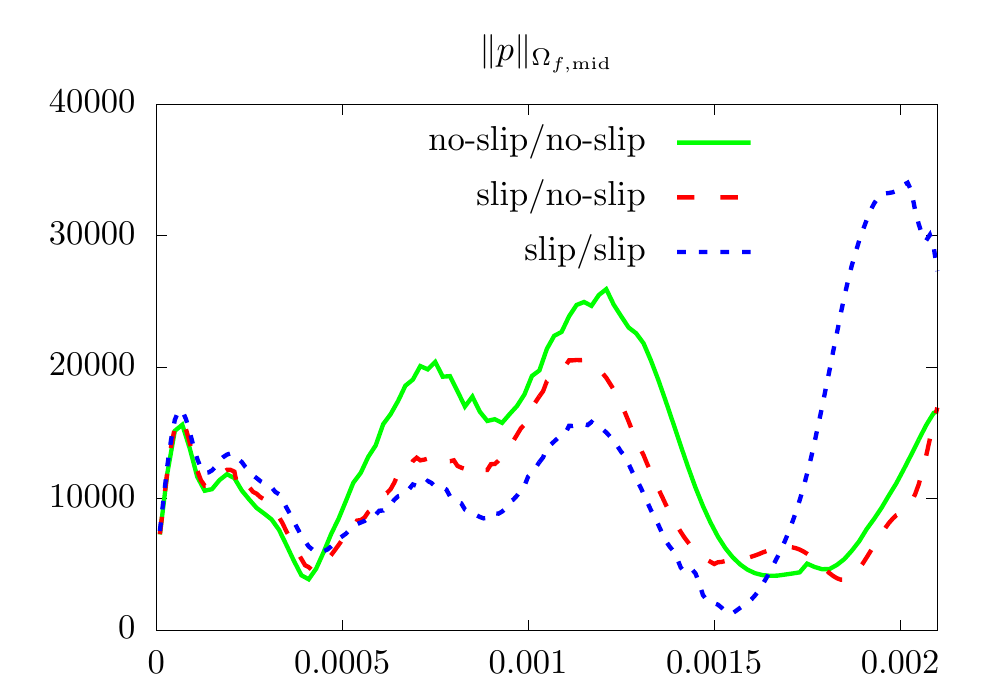}
 \end{minipage}
 \caption{\label{fig.slipvsnoslip} Comparison of slip- and no-slip interface/boundary conditions
 by means of the minimal distance $d_{\min}$ of the interface $\Gamma(t)$ to $\Gamma_w$ around 
 the contact interval (left) 
 and the $L^2$-norm of 
 the pressure over a region $\Omega_{f,{\rm mid}}(t)$ around the contact surface before contact (right) over time.
 Due to the larger pressure before contact, the impact happens later when using no-slip conditions. 
}
 \end{figure}

Next, we compare the effect of slip- and no-slip boundary and interface conditions in
Figure~\ref{fig.slipvsnoslip}. {Due to the difficulties associated
with the artificial fluid formulation and no-slip interface and boundary conditions
(see the discussion at the end of Section~\ref{sec.artf}), we use the relaxed contact formulation
on the mesh with 5120 elements in this paragraph.}

{We show results for 
\begin{itemize}
\item Slip conditions on the interface $\Gamma(t)$ and the lower wall $\Gamma_w$
\item A slip condition on $\Gamma(t)$ and a no-slip condition on $\Gamma_w$
\item No-slip conditions on $\Gamma(t)$ and $\Gamma_w$.
\end{itemize}
Note that the second option is possible, as for the relaxed contact formulation $\Gamma_w\cap \Gamma(t)=\emptyset$.

We observe that the contact condition (or more precisely the relaxed 
condition $d\cdot n_w \leq g_{\epsilon}$) is 
earlier active, when using slip-conditions: at $t_C= 1.42\cdot 10^{-3}$
for slip/slip conditions
compared to $t_C=2.02\cdot 10^{-3}$ for slip interface and no-slip boundary conditions and at 
$t_C=2.23 \cdot 10^{-3}$ for no-slip conditions on interface and boundary. The reason is that the fluid forces, and 
in particular the pressure,
that act against the contact are larger for no-slip conditions, as the fluid can not ``slip'' out
of the contact zone easily. This can be seen in the pressure plot on the right. The pressure is 
considerably larger from $t\approx 5\cdot 10^{-4}$ for the no-slip conditions until contact is reached
for the slip/slip case at $t_C= 1.42\cdot 10^{-3}$.

As we are allowing for a small gap between the solid and the ground, these results do not contradict 
the theoretical results by Gerard-Varet et al~\cite{GerardVaretetal2015}
discussed in Section~\ref{sec.slip}, who showed that (in their configuration with a rigid body) 
contact can not happen, 
when no-slip conditions are used on the interface and/or the boundary. 
As discussed in Section~\ref{sec.relaxed},
the basic assumption of the relaxed formulation is that a small or infinitesimal fluid layer remains during contact.
On the contrary, the results confirm
that contact is more likely to happen for slip-conditions, which is in agreement with the theoretical results.
}

 \paragraph{Comparison with an {explicit} ad hoc approach}
 
 {The probably simplest possibility to combine the FSI model introduced in Section~\ref{sec.FSI} and 
 the contact approach described in Section~\ref{sec.ContactModel} is to split $\Gamma$ explicitly in each time-step
into a fluid-structure interface $\Gamma_{\text{fsi}}(t_{m-1})$ and a contact surface $\Gamma_C(t_{m-1})$ 
based on the displacement
$d(t_{m-1})$ of the previous time-step
and to use the interface condition~\eqref{BalNormalForces} on $\Gamma_{\text{fsi}}(t_{m-1})$ 
and the contact condition~\eqref{ContactFormula} on $\Gamma_C(t_{m-1})$. 
A strategy of this type has been used by Hecht \& Pironneau~\cite{Pironneau2016}.}
The system of equations reads in the slip case:\\
\textit{Find $u\in {\cal V}, p \in {\cal Q}, d \in {\cal W}$ such that $\dot{d}= \partial_t d$ and}
\begin{eqnarray}\label{varForm_explicit}
 \begin{aligned}
 &\big(\partial_t u,v\big)_{\Omega_f(t)} +  \left(\sigma_f(u,p),\nabla v\right)_{\Omega_f(t)} 
 + \left({\rm div }\, u, q\right)_{\Omega_f(t)}
 +\left(\partial_t \dot{d},w\right)_{\Omega_s(t)}
 + \left(\sigma_s(d),\nabla w\right)_{\Omega_s(t)} \\ 
 %+ \left(\partial_t d,z\right)_{\Omega_s(t)} - \left(\dot{d},z\right)_{\Omega_s(t)}\\
 &\qquad- \left(n^T\sigma_f(u,p) n, (w-v)\cdot n\right)_{\Gamma_{\text{fsi}}(t_{m-1})} 
 +\gamma_{\text{fsi}} \left((\dot{d}-u)\cdot n, (w-v)\cdot n\right)_{\Gamma_{\text{fsi}}(t_{m-1})} \\
 &\qquad\qquad- \left((\dot{d}-u)\cdot n,n^T\sigma_f(v,-q) n\right)_{\Gamma_{\text{fsi}}(t_{m-1})}\\ 
 &\qquad\qquad\qquad+\gamma_C \left(P_{\gamma,s}(d), w\cdot n_w\right)_{\Gamma_C(t_{m-1})} 
 %\\ &\qquad\qquad\qquad\qquad\qquad
 =\left(f_f,v\right)_{\Omega_f(t)} + \left(f_s,w\right)_{\Omega_s(t)}\quad 
 \forall v,q,w \in {\cal V} \times {\cal Q} \times {\cal W},
\end{aligned}
\end{eqnarray}
{where $P_{\gamma,s}$ is defined in Section~\ref{sec.ContactModel} for the pure solid problem, i.e.$\,$without any fluid 
contributions.} 
We use the same numerical parameters as for the contact formulations presented in this work. 

\begin{figure}[bt]
 \centering
   \begin{minipage}{0.5\textwidth}
  \includegraphics[width=\textwidth]{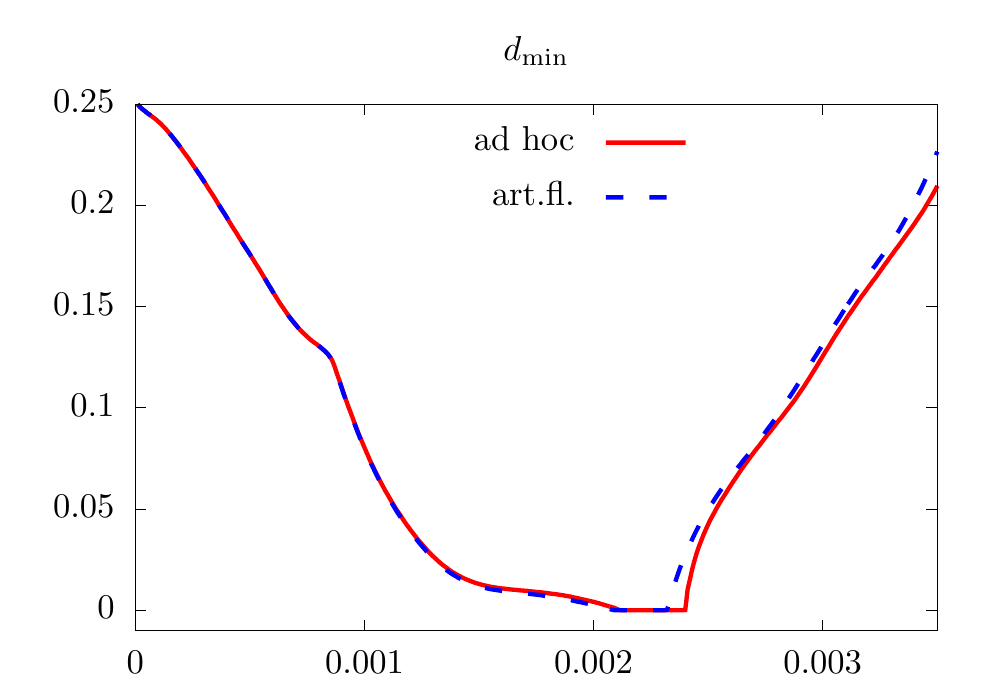}
 \end{minipage}
 \hspace{-0.5cm}
  \begin{minipage}{0.5\textwidth}
    \includegraphics[width=\textwidth]{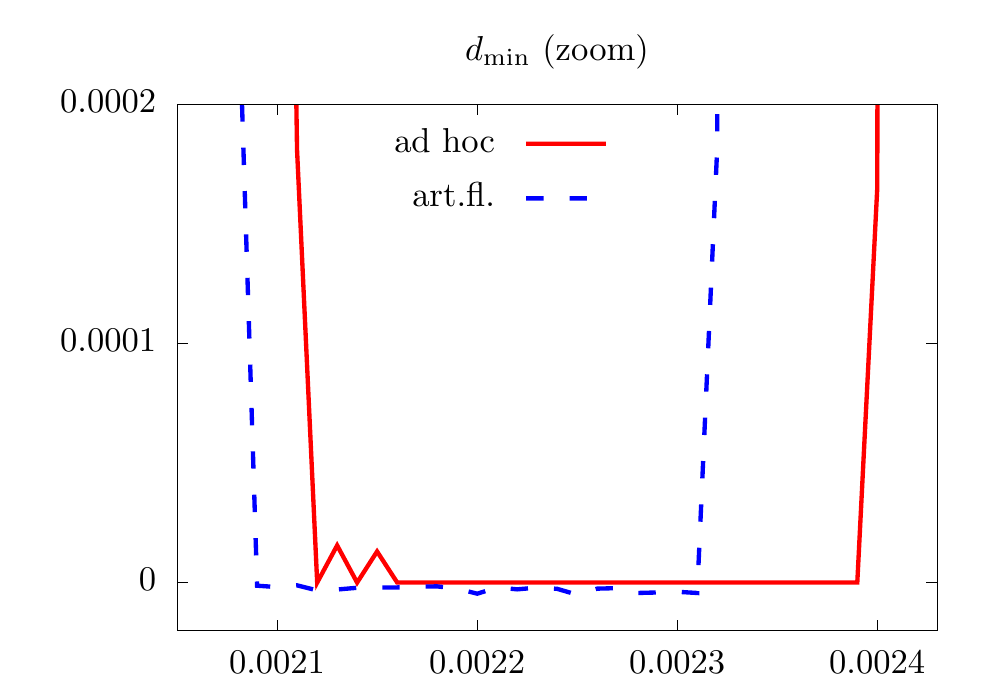}
 \end{minipage}
   \begin{minipage}{0.5\textwidth}
    \includegraphics[width=\textwidth]{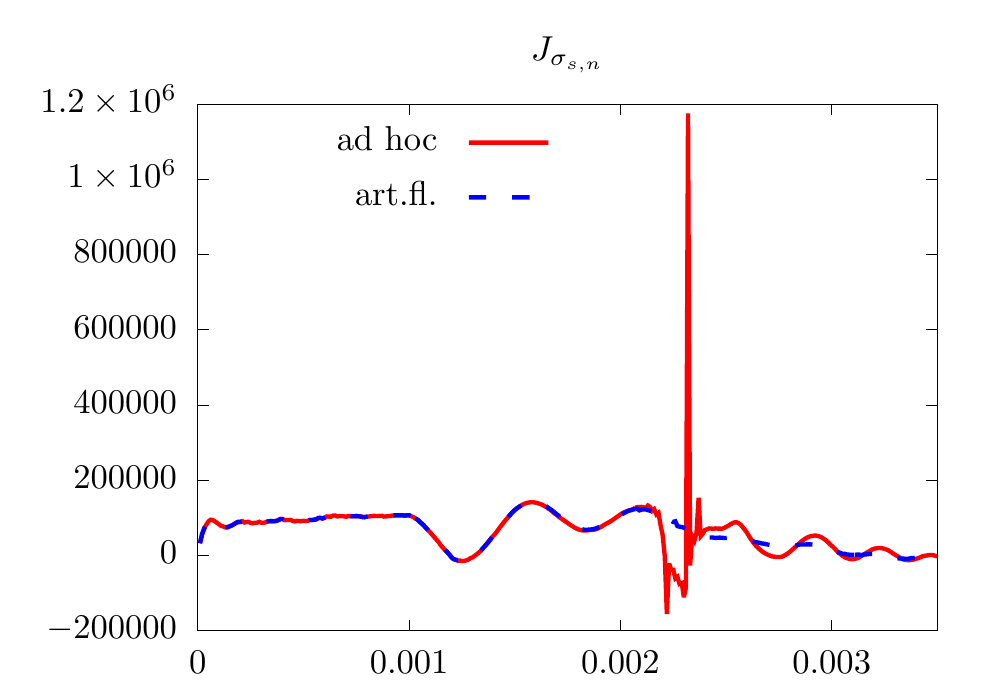}
 \end{minipage}
 \caption{\label{fig.chattering} Comparison of an explicit ad hoc approach to include the 
 contact and interface conditions with
 the approach using an artificial fluid described in Section~\ref{sec.artf}. We show the 
 minimal distance $d_{\min}$ of the 
 interface $\Gamma(t)$ to $\Gamma_w$ on the top left and a zoom-in on the top right. Chattering, i.e.$\,$an alteration between
 contact and no contact is visible for the ad hoc 
 approach at the beginning of the contact interval. The integral over the normal solid stresses
 $J_{\sigma_{s,n}}$ shown in the graph on the bottom shows large instabilities for the ad hoc approach. 
 }
 \end{figure}

To compare this approach with the artificial fluid formulation we show the minimal distance 
to the ground $d_{\min}$ and the integral over the normal solid stresses over 
$\Gamma(t) = \Gamma_{\text{fsi}}(t)\cup\Gamma_C(t)$
\begin{align*}
 J_{\sigma_{s,n}} = \int_{\Gamma(t)} \sigma_{s,n}\, ds
\end{align*}
on the finer mesh with 20480 elements in Figure~\ref{fig.chattering}. 
While the curves for $d_{\min}$ over 
the total time interval shown on the top left look similar, a zoom-in on the right shows again that the presence of the 
artificial fluid leads to an earlier time of impact. Moreover, we observe
chattering for 
the ad-hoc approach at the beginning of the 
contact interval, i.e.$\,$contact is released twice again before the solid stays in contact with $\Gamma_w$. 
The interface jumps back to the fluid domain,
with a (relatively small) minimal distance of approximately $1.6\cdot 10^{-5}$. 

{In fact, the functional $d_{\rm min}$ is not a good indicator to investigate stability for the ad hoc approach, as it 
is zero, as soon as one point of the interface lies on $\Gamma_w$. Note that this is different for the 
approaches presented in this work, where the interface can go beyond $\Gamma_w$ (or $\Gamma_{\epsilon}$ 
for the relaxed approach). In the actual computation, the interface oscillates considerably in each time-step and
contact is released and renewed frequently in different points. The functional $J_{\sigma_{s,n}}$ on the bottom left
of Figure~\ref{fig.chattering} serves 
to get a better impression of the instabilities during contact. It
oscillates throughout the contact interval including a huge peak at $t=2.32\cdot10^{-3}$. Moreover, we see that the
elastic dynamics after the contact are also significantly influenced by these instabilities. Compared to the artificial 
fluid approach the oscillations in the displacement are significantly larger after contact.}

 We have also tried to {iterate for the splitting into $\Gamma_C(t)$ and $\Gamma_{\text{fsi}}(t)$ 
 within each time-step} of 
 the ad-hoc approach, which can be seen as an active-set strategy. This did however not cure the problem, as cycling 
 between different active sets is not prevented.
 
  \section{Conclusions}
 \label{sec.conclusion}
 We have presented two consistent formulations for fluid-structure interactions with contact, both including a continuous switch
 between the FSI interface and the contact condition depending on the contact force $P_\gamma$. In contrast to certain penalty 
 approaches, the contact force is physically motivated and included in a consistent way in the variational formulations. 
 Our numerical results indicate 
 that the two proposed formulations have better stability properties than the usual ad hoc approaches and no chattering was observed in our computations.
 
 Moreover, we have derived analytically a stability result for 
 a generalised formulation including a parameter $\theta\in [0,1]$. As in the pure solid case (Chouly \& Hild~\cite{ChoulyHild2013}), this result 
 implies stability for $\theta=1$ and stability up to a term including the contact force for $\theta \neq 1$. In our computations,
 we have however not observed any stability issues for the choice $\theta=0$ either.
 
 The contact formulations were derived here for the simplified configuration of contact with a fixed and straight wall and using linear models for the fluid 
 and solid sub-problems. The algorithms can be applied to 
 more complex contact configurations by using approaches from the literature to compute the projection and the distances between 
 different surfaces~\cite{Wohlmuth2011,Yangetal2005,poulios-renard-15}. {In particular,} the extension to the incompressible Navier-Stokes equations in the fluid and to non-linear elasticity in the 
 solid {can be addressed by combining the proposed approach with the arguments recently reported in Mlika \emph{\emph{et al.}}\cite{milka-et-al-17}.}
Moreover, Coulumb or Tresca friction can {also be incorporated by} following Chouly \emph{et al.}~\cite{ChoulyCoulombFriction, ChoulyTrescaFriction}. 

\begin{acknowledgement}
The first author acknowledges support by the EPSRC grant EP/P01576X/1.
The third author was supported by the DFG Research Scholarship FR3935/1-1.
\end{acknowledgement}

\bibliographystyle{plainnat}

\begin{thebibliography}{10}

\bibitem{Wohlmuth2011}
Wohlmuth Barbara. Variationally consistent discretization schemes and numerical
  algorithms for contact problems.  {\it Acta Numerica. }2011;20:569-734.

\bibitem{Knaufetal}
Knauf Stefan, Frei Stefan, Richter Thomas, Rannacher Rolf. Towards a complete
  numerical description of lubricant film dynamics in ball bearings.  {\it
  Computational Mechanics. }2014;53(2):239--255.

\bibitem{Bruyere2012}
Bruyere Vincent, Fillot Nicolas, Morales-Espejel Guillermo~E., Vergne Philippe.
  Computational fluid dynamics and full elasticity model for sliding line
  thermal elastohydrodynamic contacts.  {\it Tribology International.
  }2012;46(1):3 - 13.

\bibitem{TezduyarSathe2007}
Tezduyar Tayfun~E., Sathe Sunil. Modeling of fluid-structure interactions with
  the space-time finite elements: solution techniques.  {\it International
  Journal for Numerical Methods in Fluids. }2007;54:855-900.

\bibitem{MayerWalletal2013}
Mayer Ursula~M, Popp Alexander, Gerstenberger Axel, Wall Wolfgang~A. 3D
  fluid--structure-contact interaction based on a combined {XFEM} {FSI} and
  dual mortar contact approach.  {\it Computational Mechanics.
  }2010;46(1):53--67.

\bibitem{DosSantosEtAl2008}
Santos N~Diniz, Gerbeau Jean-Fr{\'e}d{\'e}ric, Bourgat Jean-Fran\c{c}ois. A
  partitioned fluid--structure algorithm for elastic thin valves with contact.
  {\it Computer Methods in Applied Mechanics and Engineering.
  }2008;197(19):1750--1761.

\bibitem{AstorinoGerbeauetal2009}
Astorino Matteo, Gerbeau Jean-Fr\'ed\'eric, Pantz Olivier, Traor\'e
  Karim-Fr\'ed\'eric. Fluid-structure interaction and multi-body contact:
  Application to aortic valves.  {\it Computer Methods in Applied Mechanics and
  Engineering. }2009;198(45-46):3603 - 3612.

\bibitem{FreiRichter2017Sammelband}
Frei Stefan, Richter Thomas. An accurate {E}ulerian approach for
  fluid-structure interactions.  In:  Frei S., Holm B., Richter T., Wick T.,
  Yang H., eds. {\it Fluid-Structure Interaction: Modeling, Adaptive
  Discretization and Solvers}, Radon Series on Computational and Applied
  Mathematics. Walter de Gruyter, Berlin 2017.

\bibitem{Nitsche70}
Nitsche Joachim~A. {\"U}ber ein {V}ariationsprinzip zur {L}\"osung von
  {D}irichlet-{P}roblemen bei {V}erwendung von {T}eilr\"aumen, die keinen
  {R}andbedingungen unterworfen sind.  {\it Abhandlungen aus dem Mathematischen
  Seminar der Universit\"at Hamburg. }1970;36:9--15.

\bibitem{ChoulyHild2013}
Chouly Franz, Hild Patrick. A {N}itsche-based method for unilateral contact
  problems: numerical analysis.  {\it SIAM Journal on Numerical Analysis.
  }2013;51(2):1295--1307.

\bibitem{ChoulySymVsNonsym}
Chouly Franz, Hild Patrick, Renard Yves. Symmetric and non-symmetric variants
  of {N}itsche's method for contact problems in elasticity: {T}heory and
  numerical experiments.  {\it Mathematics of Computation.
  }2015;84(293):1089--1112.

\bibitem{Choulyetal2015}
Chouly Franz, Hild Patrick, Renard Yves. A {N}itsche finite element method for
  dynamic contact: 1. Space semi-discretization and time-marching schemes.
  {\it ESAIM: Mathematical Modelling and Numerical Analysis.
  }2015;49(2):481--502.

\bibitem{AlartCurnier91}
Alart Pierre, Curnier Alain. A mixed formulation for frictional contact
  problems prone to {N}ewton like solution methods.  {\it Computer Methods in
  Applied Mechanics and Engineering. }1991;92(3):353--375.

\bibitem{ChoulyTrescaFriction}
Chouly Franz. An adaptation of Nitsche's method to the Tresca friction problem.
   {\it Journal of Mathematical Analysis and Applications. }2014;411:329-339.

\bibitem{ChoulyCoulombFriction}
Chouly Franz, Hild Patrick, Lleras Vanessa, Renard Yves. {Nitsche-based finite
  element method for contact with Coulomb friction}
  https://hal.archives-ouvertes.fr/hal-01654487/file/authorv3.pdf,2017.

\bibitem{BurmanetalObstacle}
Burman Erik, Hansbo Peter, Larson Mats~G., Stenberg Rolf. Galerkin least
  squares finite element method for the obstacle problem.  {\it Computer
  Methods in Applied Mechanics and Engineering. }2017;313(Supplement C):362 -
  374.

\bibitem{BurmanHansboLarsonMembrane}
{Burman} Erik, {Hansbo} Peter, {Larson} Mats~G.. {Augmented Lagrangian and
  Galerkin least squares methods for membrane contact}.  {\it ArXiv e-prints.
  }2017;.
\newblock http://adsabs.harvard.edu/abs/2017arXiv171104494B.

\bibitem{Annavarapuetal2014}
Annavarapu Chandrasekhar, Hautefeuille Martin, Dolbow John~E.. A {N}itsche
  stabilized finite element method for frictional sliding on embedded
  interfaces. Part {I}: Single interface.  {\it Computer Methods in Applied
  Mechanics and Engineering. }2014;268:417 - 436.

\bibitem{BurmanHansbo2018}
Burman Erik, Hansbo Peter. Deriving robust unfitted finite element methods from
  augmented Lagrangian formulations.  In:  Bordas S.P.A., Burman E.N., Larson
  M.G., Olshanskii M.A., eds. {\it Geometrically Unfitted Finite Element
  Methods and Applications - Proceedings of the UCL-workshop 2016}, Springer
  2017 (pp. 1--24).

\bibitem{Hillairet2d}
Hillairet Matthieu. Lack of collision between solid bodies in a 2D
  incompressible viscous flow.  {\it Communications in Partial Differential
  Equations. }2007;32(9):1345-1371.

\bibitem{HeslaPhD}
{Hesla} Todd~I.. Collisions of smooth bodies in viscous fluids: A mathematical
  investigation.
\newblock PhD thesis.{ }University of Minnesota.{ }2004.

\bibitem{HillairetTakahashi3d}
Hillairet Matthieu, Takahashi Tak\'eo. Collisions in three-dimensional fluid
  structure interaction problems.  {\it SIAM Journal on Mathematical Analysis.
  }2009;40(6):2451-2477.

\bibitem{GerardVaretetal2015}
Gerard-Varet David, Hillairet Matthieu, Wang Chao. The influence of boundary
  conditions on the contact problem in a 3d {N}avier-{S}tokes flow.  {\it
  Journal de Math{\'e}matiques Pures et Appliqu{\'e}es. }2015;103:1--38.

\bibitem{GerardVaretHillairet}
G{\'e}rard-Varet David, Hillairet Matthieu. Regularity Issues in the Problem of
  Fluid Structure Interaction.  {\it Archive for Rational Mechanics and
  Analysis. }2010;195(2):375--407.

\bibitem{Wang2014}
Wang Chao. Strong solutions for the fluid--solid systems in a 2-D domain.  {\it
  Asymptotic Analysis. }2014;89(3-4):263--306.

\bibitem{GrandmontHillairet}
Grandmont C{\'e}line, Hillairet Matthieu. Existence of global strong solutions
  to a beam--fluid interaction system.  {\it Archive for Rational Mechanics and
  Analysis. }2016;220(3):1283--1333.

\bibitem{GrandmontEtAlInBook}
Grandmont C{\'e}line, Luk{\'a}{\v{c}}ov{\'a}-Medvid{\'o}v{\'a} M{\'a}ria,
  Ne{\v{c}}asov{\'a} {\v{S}}{\'a}rka. Mathematical and numerical analysis of
  some FSI problems:1--77.
\newblock Basel: Springer Basel 2014.

\bibitem{MuhaCanic}
Muha Boris, \v{C}ani\'c Sun\v{c}ica. Existence of a weak solution to a
  fluid-elastic structure interaction problem with the {N}avier slip boundary
  condition.  {\it Journal of Differential Equations. }2016;260(12):8550 -
  8589.

\bibitem{Puso2004}
Puso Michael~A. A 3D mortar method for solid mechanics.  {\it International
  Journal for Numerical Methods in Engineering. }2004;59(3):315--336.

\bibitem{Yangetal2005}
Yang Bin, Laursen Tod~A, Meng Xiaonong. Two dimensional mortar contact methods
  for large deformation frictional sliding.  {\it International Journal for
  Numerical Methods in Engineering. }2005;62(9):1183--1225.

\bibitem{ChoulyMlikaRenard18}
Chouly Franz, Mlika Rabii, Renard Yves. An unbiased Nitsche's approximation of
  the frictional contact between two elastic structures.  {\it Numerische
  Mathematik. }2018;139(3):593--631.

\bibitem{milka-et-al-17}
Mlika Rabii, Renard Yves, Chouly Franz. An unbiased {N}itsche's formulation of
  large deformation frictional contact and self-contact.  {\it Computer Methods
  in Applied Mechanics and Engineering. }2017;325:265--288.

\bibitem{DunneRannacher}
Dunne Thomas, Rannacher Rolf. Adaptive finite element approximation of
  fluid-structure interaction based on an {E}ulerian variational formulation.
  In:  Bungartz H.-J., Sch{\"a}fer M., eds. {\it Fluid-Structure Interaction:
  Modeling, Simulation, Optimization}, Lecture Notes in Computational Science
  and Engineering. Springer 2006 (pp. 110-145).

\bibitem{Cottetetal}
Cottet Georges-Henri, Maitre Emmanuel, Milcent Thomas. {E}ulerian formulation
  and level set models for incompressible fluid-structure interaction.  {\it
  ESAIM: Mathematical Modelling and Numerical Analysis. }2008;42(3):471--492.

\bibitem{Richter2012b}
Richter Thomas. A Fully {E}ulerian Formulation for Fluid-Structure
  Interactions.  {\it Journal of Computational Physics. }2013;233:227-240.

\bibitem{Pironneau2016}
Hecht Fr\'ed\'eric, Pironneau Olivier. An energy stable monolithic Eulerian
  fluid-structure finite element method.  {\it International Journal for
  Numerical Methods in Fluids. }2017;85(7):430--446.

\bibitem{Peskin1972}
Peskin Charles~S. Flow patterns around heart valves: a numerical method.  {\it
  Journal of Computational Physics. }1972;10(2):252--271.

\bibitem{BoffiGastaldi2003}
Boffi Daniele, Gastaldi Lucia. A finite element approach for the immersed
  boundary method.  {\it Computers \& Structures. }2003;81(8):491--501.

\bibitem{ZhangGerstenberger}
Zhang Lucy, Gerstenberger Axel, Wang Xiaodong, Liu Wing~Kam. Immersed finite
  element method.  {\it Computer Methods in Applied Mechanics and Engineering.
  }2004;193(21):2051--2067.

\bibitem{LegayChessaBelytschko2006}
Legay Antoine, Chessa Jack, Belytschko Ted. An {E}ulerian-{L}agrangian method
  for fluid-structure interaction based on level sets.  {\it Computer Methods
  in Applied Mechanics and Engineering. }2006;195:2070-2087.

\bibitem{GerstenbergerWall}
Gerstenberger Axel, Wall Wolfgang~A. An extended finite element
  method/{L}agrange multiplier based approach for fluid--structure interaction.
   {\it Computer Methods in Applied Mechanics and Engineering.
  }2008;197(19):1699--1714.

\bibitem{BurmanFernandez2014}
Burman Erik, Fern{\'a}ndez Miguel~A. An unfitted {N}itsche method for
  incompressible fluid--structure interaction using overlapping meshes.  {\it
  Computer Methods in Applied Mechanics and Engineering. }2014;279:497--514.

\bibitem{alauzet-et-al-15}
Alauzet Fr{\'e}d{\'e}ric, Fabr{\'e}ges Benoit, Fern{\'a}ndez Miguel~Angel,
  Landajuela Mikel. {Nitsche-XFEM for the coupling of an incompressible fluid
  with immersed thin-walled structures}.  {\it Computer Methods in Applied
  Mechanics and Engineering. }2016;301:300--335.

\bibitem{MassingLarsonetal}
Massing Andr{\'e}, Larson Mats, Logg Anders, Rognes Marie. A Nitsche-based cut
  finite element method for a fluid-structure interaction problem.  {\it
  Communications in Applied Mathematics and Computational Science.
  }2015;10(2):97--120.

\bibitem{HansboHermanssonSvedberg2004}
Hansbo Peter, Hermansson Joakim, Svedberg Thomas. Nitsche's method combined
  with space--time finite elements for ALE fluid--structure interaction
  problems.  {\it Computer Methods in Applied Mechanics and Engineering.
  }2004;193(39-41):4195--4206.

\bibitem{kamensky-et-al-15}
Kamensky David, Hsu Ming~Chen, Schillinger Dominik, et al. An immersogeometric
  variational framework for fluid--structure interaction: {A}pplication to
  bioprosthetic heart valves.  {\it Computer Methods in Applied Mechanics and
  Engineering. }2015;284:1005--1053.

\bibitem{boilevinkayl:hal-01704575}
Boilevin-Kayl Ludovic, Fern{\'a}ndez Miguel~A., Gerbeau Jean-Fr{\'e}d{\'e}ric.
  {\it {Numerical methods for immersed FSI with thin-walled structures}.
  }Research Report RR-9151: {Inria Paris}; 2018.
\newblock \url{https://hal.inria.fr/hal-01704575}.

\bibitem{DunnePhD}
Dunne Thomas. Adaptive finite element approximation of fluid-structure
  interaction based on {E}ulerian and Arbitrary Lagrangian-{E}ulerian
  variational formulations.
\newblock PhD thesis{ }Heidelberg University{ }2007.

\bibitem{FreiRichter2014}
Frei Stefan, Richter Thomas. A locally modified parametric finite element
  method for interface problems.  {\it SIAM Journal on Numerical Analysis.
  }2014;52(5):2315-2334.

\bibitem{FreiRichter2017}
Frei Stefan, Richter Thomas. A second order time-stepping scheme for parabolic
  interface problems with moving interfaces.  {\it ESAIM: Mathematical
  Modelling and Numerical Analysis. }2017;51(4):1539--1560.

\bibitem{Burman2010}
Burman Erik. Ghost penalty.  {\it Comptes Rendus Mathematique.
  }2010;348(21-22):1217--1220.

\bibitem{RichterBuch}
Richter Thomas. {\it Finite Elements for Fluid-Structure Interactions. Models,
  Analysis and Finite Elements.} Lecture Notes in Computational Science and
  Engineering, vol. 118: .
\newblock Springer; 2017.

\bibitem{FreiPhD}
Frei Stefan. Eulerian finite element methods for interface problems and
  fluid-structure interactions.
\newblock PhD thesis.{ }Heidelberg University{ }2016.
\newblock http://www.ub.uni-heidelberg.de/archiv/21590.

\bibitem{CimolinDiscacciati}
Cimolin Flavio, Discacciati Marco. Navier--Stokes/Forchheimer models for
  filtration through porous media.  {\it Applied Numerical Mathematics.
  }2013;72:205--224.

\bibitem{IlievLaptev}
Iliev Oleg, Laptev Vsevolod. On numerical simulation of flow through oil
  filters.  {\it Computing and Visualization in Science. }2004;6(2):139--146.

\bibitem{Angot1999}
Angot Philippe. Analysis of singular perturbations on the Brinkman problem for
  fictitious domain models of viscous flows.  {\it Mathematical Methods in the
  Applied Sciences. };22(16):1395-1412.

\bibitem{Kamenskyetal2017}
Kamensky David, Evans John~A., Hsu Ming-Chen, Bazilevs Yuri. Projection-based
  stabilization of interface Lagrange multipliers in immersogeometric
  fluid-thin structure interaction analysis, with application to heart valve
  modeling.  {\it Computers and Mathematics with Applications. }2017;74(9):2068
  - 2088.

\bibitem{BrezziPitkaeranta1984}
Brezzi Franco, Pitk{\"a}ranta Juhani. On the stabilization of finite element
  approximations of the Stokes equations.  In:  W.~Hackbusch, ed. {\it
  Efficient solutions of elliptic systems}, Springer 1984 (pp. 11--19).

\bibitem{BeckerBraack2001}
Becker Roland, Braack Malte. A finite element pressure gradient stabilization
  for the {S}tokes equations based on local projections.  {\it Calcolo.
  }2001;38(4):173-199.

\bibitem{HughesFrancaBalestra1986}
Hughes Thomas~JR, Franca Leopoldo~P, Balestra Marc. A new finite element
  formulation for computational fluid dynamics: V. Circumventing the
  {B}abu{\v{s}}ka-{B}rezzi condition: a stable {P}etrov-{G}alerkin formulation
  of the {S}tokes problem accommodating equal-order interpolations.  {\it
  Computer Methods in Applied Mechanics and Engineering. }1986;59(1):85--99.

\bibitem{BurmanHansbo2006}
Burman Erik, Hansbo Peter. Edge stabilization for the generalized {S}tokes
  problem: a continuous interior penalty method.  {\it Computer Methods in
  Applied Mechanics and Engineering. }2006;195(19):2393--2410.

\bibitem{BurmanHansbo2012}
Burman Erik, Hansbo Peter. Fictitious domain finite element methods using cut
  elements: {II}. {A} stabilized {N}itsche method.  {\it Applied Numerical
  Mathematics. }2012;62(4):328--341.

\bibitem{TureketalNewton}
Mandal Saptarshi, Ouazzi Abderrahim, Turek Stefan. Modified Newton Solver for
  Yield Stress Fluids.  In:  Karas{\"o}zen B{\"u}lent, Manguo{\u{g}}lu Murat,
  Tezer-Sezgin M{\"u}nevver, G{\"o}ktepe Serdar, U{\u{g}}ur {\"O}m{\"u}r, eds.
  {\it Numerical Mathematics and Advanced Applications ENUMATH 2015},
  :481--490{ }Springer International Publishing; 2016.

\bibitem{FreiToAppear}
Frei Stefan. An edge-based pressure stabilisation technique for finite elements
  on arbitrarily anisotropic meshes.
  http://www.homepages.ucl.ac.uk/\~{}ucahfre/Edgestab\_aniso.pdf,
  submitted(2017).

\bibitem{BesierWollner}
Besier Michael, Wollner Winnifried. On the pressure approximation in
  nonstationary incompressible flow simulations on dynamically varying spatial
  meshes.  {\it International Journal for Numerical Methods in Fluids.
  }2012;69(6):1045--1064.

\bibitem{Gascoigne}
Becker Roland, Braack Malte, Meidner Dominik, Richter Thomas, Vexler Boris. The
  finite element toolkit {G}ascoigne3d  http://www.gascoigne.uni-hd.de; .

\bibitem{poulios-renard-15}
Poulios Konstantinos, Renard Yves. An unconstrained integral approximation of
  large sliding frictional contact between deformable solids.  {\it Computers
  \& Structures. }2015;153:75--90.

\end{thebibliography}

\end{document}